%% file: F.tex
\documentstyle[11pt]{article}

\input{macros}

\begin{document}

\input {title} 
\input {sect1}

\input {sect2}

\input {sect3}

\input {sect4}
\input {sect5}

\input biblio\end{document}

%% file: macros
\newtheorem{theo}{Theorem}[section]
\newtheorem{prop}{Proposition}[section]
\newtheorem{cor}{Corollary}[section]
\newtheorem{defi}{Definition}[section]
\newtheorem{rem}{Remark}[section]
\newtheorem{lem}{Lemma}[section]

   {
   \catcode`\<=\active
   \catcode`\>=\active
   \gdef\bra<#1|{\left<\,#1\>\right|}
   \gdef\ket|#1>{\left|\>#1\,\right>}
   \gdef\braket<#1|#2>{\left<#1\>\mid\>#2\right>}
   \gdef\scalar<#1.#2>{\left<#1\cdot#2\right>}
   }

\def\Proof{\medskip\noindent {\it Proof --- \ }}
\def\QED{\hfill\nobreak\hbox{${}_\Box$}\par\medbreak}


%% file: title
\begin{titlepage}
\begin{flushright}
ENSLAPP-L-601/96\\
US-FT/47-96\\
\end{flushright}
\par \vskip .1in \noindent

\begin{center}
{\LARGE Drinfel'd Twists and Algebraic Bethe Ansatz}\\
 \end{center}
  \par \vskip .3in \noindent

\begin{center}

      {\bf J. M. MAILLET $^{(1)}$ and J. SANCHEZ de SANTOS $^{(2)}$}
  \par \vskip .1in \noindent

{\sl $^{(1)}$ Laboratoire de Physique Th\'eorique ENSLAPP $^{*}$\\
       ENS Lyon, 46 all\'ee d'Italie 69364 Lyon CEDEX 07
       France\\
and\\
$^{(2)}$ Departamento de Fisica de Particulas\\
                     Universidad de Santiago
                     E-15706, Santiago de Compostela, Spain.}\\[0.6in]
\end{center}

\par \vskip .10in
\begin{center}
{\bf Abstract}\\
\end{center}

\begin{quote}
We study representation theory of Drinfel'd twists, in terms of what we call    $F$-matrices, associated to finite dimensional irreducible modules of quantum  affine algebras, and which factorize the corresponding (unitary) $R$-matrices. We construct explicitly such {\em factorizing $F$-matrices} for irreducible finite $n^{th}$ tensor products of the fundamental evaluation representations of the quantum affine algebra ${\cal U}_q (\hat{sl}_2)$ and the Yangian $Y (sl_2)$. We then apply these  constructions to the $XXX$-$1 \over 2$ and $XXZ$-$1 \over 2$ Heisenberg (inhomogeneous) quantum spin chains of finite lenght $n$ in the framework of the Algebraic Bethe Ansatz. In particular, we show that these factorizing $F$-matrices diagonalize the generating matrix of scalar products of quantum states of these models. They also diagonalize the diagonal (operator) entries of the quantum monodromy matrix. Due to their algebraic properties, these $F$-matrices are shown to obey simple difference equations. This leads to a natural $F$-basis for the quantum space of states of the inhomogeneous $XXX$-$1 \over 2$ quantum spin chain of finite lenght in which this model can be interpreted as a {\em diagonal dressing} of the corresponding Gaudin model.  
\end{quote}
\par \vskip 0.5in

\begin{flushleft}
\rule{5.1 in}{.007 in}\\
$^{*}${\small URA 1436 ENSLAPP du CNRS, associ\'ee \`a  l'Ecole
Normale Sup\'erieure de Lyon, \`a l'Universit\'e de Savoie et au Laboratoire d'Annecy de Physique des  Particules.}\\
{\small This work is supported by CNRS (France), the EC contract ERBCHRXCT920069} and grant PB93-0344 from DGICYT (Spain).\\
{\small { email: maillet@enslapp.ens-lyon.fr}}\\[0.2 in]

December 1996
\end{flushleft}

\end{titlepage}

%% file: sect1
\section{Introduction}

Quasi-triangular Hopf (and quasi-Hopf) Quantum 
Universal Envelopping (QUE) algebras \cite{Drin1,Drin2,FRT,Jim1,Jim2,Jim3,KR1,Skl1} provide the natural framework for two-dimensional integrable models solvable through the 
Algebraic Bethe Ansatz (or Quantum Inverse Scattering Method (QISM))  
\cite{Bax1,F1,FST,Gaud1,JM1,KBI,KS1,Skl2}. At the root of this approach to quantum integrable models 
are two basic formulas. The first expresses commutation relations for the 
operator valued entries  of a matrix $T$ (the quantum monodromy matrix) as,
\begin{equation}
R_{12}\ T_{1}\  T_{2}\ =\ T_{2}\ T_{1}\ R_{12}\ , 
\label{eq:rtt}
\end{equation}
where $T_{1} \in End\ (V_{1}) \otimes {\cal A}_{R}$ and $T_{2} \in End\ (V_{2}) \otimes {\cal A}_{R}$, ${\cal A}_{R}$ is an associative algebra, $V_{i}$, $i\ =\ 1, 2, \dots$ are finite dimensional vector spaces, and $R_{12} \in End\ (V_{1} \otimes V_{2})$ is the so-called numerical $R$-matrix. The second formula is 
the   Yang-Baxter equation in $V_1 \otimes V_2 \otimes V_3$ for the $R$-matrix itself,
\begin{equation}
R_{12} \ R_{13} \ R_{23} \ =\ R_{23} \ R_{13} \ R_{12}\ ,
\label{eq:rrr}
\end{equation}
where $R_{12}$, $R_{13}$, $R_{23}$ act in the tensor product $V_1 \otimes V_2 \otimes V_3$, their action being equal to $R$ in the two spaces with corresponding indices and being the identity in the omitted third space.\\

The matrix elements of $T$ generates an associative algebra ${\cal A}_{R}$ with 
quadratic defining relations (\ref{eq:rtt}). ${\cal A}_{R}$ inherits a Hopf algebra structure with co-multiplication $\delta$ defined by,
\begin{eqnarray*}
{\delta}\ ( T_{ik} )\ &=&\ \sum_{k}\ T_{ik} \otimes T_{kj}\ , \\
{\delta}\ ({\mathbf 1})\ &=&\ {\mathbf 1} \otimes {\mathbf 1}\ ,
\label{eq:delta}
\end{eqnarray*}
or in compact notations, ${\delta}\ ( T_{1} )\ =\ T_{1} \otimes T_{1}$.
As stressed in \cite{Drin2,FRT} this algebraic picture is dual to the one of 
quasi-triangular Hopf QUE-algebras.\\

In \cite{Drin2,Drin3}, but also in a more elementary way in the beautiful paper \cite{Drin4}, Drinfel'd introduced the notion of twisting (or deformation) for quasi-Hopf 
algebras and quasi-triangular quasi-Hopf algebras.\\
Let (${\cal A}, \Delta, \epsilon, S, R, \Phi$) be a quasi-triangular 
quasi-Hopf algebra over ${\mathbf C}$ and  $F \in {\cal A} \otimes {\cal A}$ an invertible element such that $(\epsilon \otimes id)\ F\ =\ (id \otimes \epsilon)\  F\ =\ \mathbf 1$. The twist of $\cal A$ generated by $F$ defines a new  quasi-triangular 
quasi-Hopf algebra ($\tilde{\cal A}, \tilde{\Delta}, \tilde{\epsilon}, \tilde{S},  
\tilde{R}, \tilde{\Phi}$), having the same elements and product law as 
$\cal A$,  and such that in particular,
\begin{eqnarray*}
\tilde{R}_{12}\ &=&\ F_{21}\ R_{12}\ F^{-1}_{12}\ , \\
\tilde{\Phi}_{123}\ &=&\ F_{23}\ F_{1,23}\ {\Phi}_{123}\ F^{-1}_{12,3}\ 
F^{-1}_{12}\ , \\
\tilde{\Delta} (x)\ &=&\ F_{12}\ {\Delta} (x)\ F^{-1}_{12}\ ,\ x \in {\cal A}\ ,
\label{eq:tilde}
\end{eqnarray*}
where we adopted the following compact notations for the co-multiplication actions on $F$'s, 
\begin{eqnarray*}
({\mathbf 1} \otimes \Delta)\ F\ &=&\ F_{1,23}\ , \\ 
(\Delta \otimes {\mathbf 1})\ F\ &=&\ F_{12,3}\ .
\label{equation:compactd}
\end{eqnarray*}
The Hopf algebra case is given by putting the value of $\Phi$ equal to the identity.\\

An important class of twists is given by elements $F \in {\cal A} \otimes {\cal A}$ relating two Hopf algebraic structures on $\cal A$, namely twists such that $\Phi\ =\ \tilde{\Phi}\ =\ {\mathbf 1}$. Such twists have already been used in the context of quantum algebras \cite{Drin3,Drin4,E1,KT1,JMM1,R1}, in particular in \cite{KT1} to relate two different presentations of quantum affine algebras.\\
Another important class of twists (which is central in \cite{Drin3}) is given by elements $F \in {\cal A} \otimes {\cal A}$ relating a (non-co-commutative) co-multiplication $\Delta$ to a co-commutative one $\Delta_0$.\\
It is shown in \cite{Drin3} (Proposition 3.16) that if $\cal G$ is a finite dimensional simple Lie algebra over $\mathbf C$  with fixed invariant scalar product associated to a given element $t \in {\cal G} \otimes {\cal G}$, the standard quasi-triangular quantization of the universal envelopping algebra (${\cal U} {\cal G}, \Delta_0$) as defined in \cite{Drin2}, namely the  quasi-triangular Hopf QUE-algebra over ${\mathbf C}[[h]]$, $({\cal U}({\cal G})_{h}, {\Delta}, R, S, \epsilon)$, is related, up to isomorphism, by a twist $F$, to the quasi-triangular quasi-Hopf algebra $({\cal U}({\cal G})_{h}, {\Delta}_0, R_0, S, \epsilon, \Phi_0)$, where $\Delta_0$ is the standard co-commutative co-multiplication of ${\cal U} {\cal G}$, $R_0\ =\ e^{ht \over 2}$, and $\Phi_{0}$ is defined by means of the Knizhnik-Zamolodchikov system of equations \cite{KZ}. Similar results for the case of quantum affine algebras are announced in \cite{EK1}. However, explicit universal formulas for such twists are not available yet.\\
In the triangular  case where $R_{21}\ R_{12}\ =\ {\mathbf 1}$, it is shown in \cite{Drin3} (Proposition 3.6), that any triangular quasi-Hopf 
QUE-algebra over ${\mathbf C}[[h]]$ can, by a suitable twist, be brought into the form  $R\ =\ \mathbf 1$ and $\Phi\ =\ \mathbf 1$.\\
If applied to the triangular {\em Hopf} QUE-algebra situation, it means that the corresponding twist $F$ belongs to both of the classes defined above, namely, it relates the quantum co-multiplication $\Delta$ to some co-commutative one $\Delta_0$ (since $\tilde{R}\ =\ {\mathbf 1}$) and both $\Phi$ and $\tilde{\Phi}$ are equal to the identity. Such twists will be called {\em factorizing twists}, referring in particular to the fact that they factorize the $R$-matrix as,
\begin{equation}
R_{12}\ =\ F_{21}^{-1}\ F_{12}\ .
\label{eq:rff}
\end{equation}
In the dual picture, for a Yang-Baxter type algebra ${\cal A}_{R}$, a twist $F$ will leave unchanged the co-multiplication $\delta$ 
while deforming the product law between the (quantum operator) entries of the 
$T$-matrices as,
\begin{equation}
T_{1}\ {\ast}_{F}\ T_{2}\ =\ F_{12}\ T_{1} \cdot T_{2}\ F^{-1}_{12}\ ,
\label{eq:starf}
\end{equation} 
with $\tilde{R}_{12}\ =\ F_{21}\ R_{12}\ F^{-1}_{12}$, the relation 
(\ref{eq:rtt}) being covariantly changed as,
\begin{equation}
\tilde{R}_{12}\ T_{1}\ {\ast}_{F}\ T_{2}\ =\ T_{2}\ {\ast}_{F}\ T_{1}\ \tilde{R}_{12}\ .
\label{eq:rttf}
\end{equation}
Note that, $R_{12}$ stands for the representation of some universal $R$-matrix as, $R_{12}\ =\ (\rho_1 \otimes \rho_2) R$, while $F_{12}\ =\ (\rho_1 \otimes \rho_2) F$, $\rho_1$ and $\rho_2$ being two representations of ${\cal A}$.\\
The elements of ${\cal A}_{R}$ being generated by 
the matrix elements of arbitrary (finite) products of the $T$-matrices as 
$T_{q}\ =\ T_{1}\ T_{2} \cdots T_{n}$ where $q\ =\ (1,2,...,n)$,  
defining the operation ${\ast}_{F}$ on ${\cal A}_{R}$ amounts to give the value of the twist $F$ on such tensor products of representations. 
In other words, to achieve this definition, one
needs  to construct arbitrary left and right co-multiplication actions on $F \in {\cal A} \otimes {\cal A}$ leading to representations of $F$ on arbitrary products of the $T$-matrices $T_{q}$.\\
For triangular Hopf algebras $\cal A$, a factorizing twist $F$ corresponding to the above Drinfel'd's theorem, namely such that $\tilde{R}\ =\ {\mathbf 1}$ and $\tilde{\Phi}\ =\ {\mathbf 1}$, gives the following relation in the algebra ${\cal A}_R$,
\begin{equation}
T_{1}\ {\ast}_{F}\ T_{2}\ =\ T_{2}\ {\ast}_{F}\ T_{1}\ , 
\label{eq:starff}
\end{equation}
where we used the decomposition (\ref{eq:rff}) of the $R$-matrix in terms of the twist $F$. This is the dual counterpart of $\Delta_0$ being co-commutative on ${\cal A}$. The corresponding $\Phi$ and $\tilde{\Phi}$ being equal to the identity, $F$ satisfies the following cocycle relation,
\begin{equation}
F_{23}\ F_{1,23}\ =\ F_{12}\ F_{12,3}\ .
\label{eq:ff}
\end{equation}
More generally, a twist $F$, realizing the deformation of a non-commutative co-multiplication $\Delta$ to a co-commutative one $\Delta_0$, will put into correspondence a non-co-commutative (or non-commutative) algebraic context with a co-commutative (or commutative) one.\\

In our opinion, this simple observation has not been used up to now in its full strenght in the investigation of quantum integrable models associated to quantum algebras, in particular in the framework of the Algebraic Bethe Ansatz. This is probably mainly due to the fact that the knowledge of the corresponding  twists $F$, namely, in particular their representation theory relevant to these models, is far from being as well known 
as for the corresponding $R$-matrices. 
Indeed, the co-multiplication action on $F$ is not in general of multiplicative type 
as it is the case for $R$-matrices of quasi-triangular Hopf algebras.\\
This feature makes the problem of computing representation theory of twists $F$ quite involved, not speaking about universal formulas in ${\cal A} \otimes {\cal A}$.\\

Typical examples of Hopf QUE-algebras arising in the context of Algebraic Bethe Ansatz are quantum affine algebras and their associated Yangians. They lead in particular to trigonometric and rational spectral parameter dependent solutions of the Yang-Baxter equation \cite{Drin1,Jim1}. Such solutions are for example associated to intertwining operators between irreducible tensor products of finite dimensional irreducible modules $V_i\ (\lambda_i)$ of quantum affine algebras \cite{CP1}.\\
In such representations, the $R$-matrix $R_{12}\ (\lambda_1,\ \lambda_2)$ is an element of $End\ (V_1 (\lambda_1) \otimes V_2 (\lambda_2))$ depending on two spectral parameters $\lambda_1$ and $\lambda_2$. In the irreducible case, this $R$-matrix is unitary \cite{CP1,FR1}, namely,
\begin{equation}
R_{21}\ (\lambda_2,\ \lambda_1)\ R_{12}\ (\lambda_1,\ \lambda_2)\ =\ {\mathbf 1}\ ,
\label{eq:url}
\end{equation}
and satisfies the Yang-Baxter equation,
\begin{equation}
R_{12}(\lambda_1, \lambda_2 )\ R_{13}(\lambda_1, \lambda_3)\ R_{23}(\lambda_2, \lambda_3)\ =\ 
R_{23}(\lambda_2, \lambda_3)\ R_{13}(\lambda_1, \lambda_3)\ R_{12}(\lambda_1, \lambda_2) \ ,
\label{eq:qr}
\end{equation}
where $\lambda_1$, $\lambda_2$, $\lambda_2$ are spectral parameters 
attached to the vector spaces $V_1$, $V_2$ and $V_3$ respectively.\\

Although quantum affine algebras are  quasi-triangular rather than triangular Hopf QUE-algebras, for such unitary $R$-matrices, associated to finite dimensional irreducible representations, it is possible to define (see section 2.4) the notion of {\em factorizing $F$-matrices} in much the same way as the notion of factorizing twists arises for triangular Hopf algebras. As part of their definition, such factorizing $F$-matrices decompose the associated 
unitary  $R$-matrix as,
\begin{equation}
R_{12}(\lambda_1,\ \lambda_2)\ =\ F_{21}^{-1}\ (\lambda_2,\ \lambda_1)\ F_{12}\ (\lambda_1,\ \lambda_2)\ ,
\label{eq:rffl}
\end{equation}
with $F_{12}\ (\lambda_1,\ \lambda_2) \in End\ (V_1 (\lambda_1) \otimes V_2 (\lambda_2))$.\\
Moreover, for irreducible finite tensor product of finite dimensional irreducible  modules $V_{q}\ (\lambda_q) \equiv V_{1}\ (\lambda_1) \otimes \cdots \otimes V_{N}\ (\lambda_n)$, $(q \equiv 1, 2, \dots , n)$ and $\lambda_q \equiv (\lambda_1,\ \dots\ ,\ \lambda_n)$, one can define the action of the permutation group $S_n$ induced by the $R$-matrix. Let $R_q^{\sigma}\ (\lambda_q)$ be the $R$-matrix associated to the element $\sigma \in S_n$. A factorizing $F$-matrix acting in the tensor product $V_q\ (\lambda_q)$ gives the following decomposition of the $R$-matrix,
\begin{equation}
R_q^{\sigma}\ (\lambda_q)\ =\ F_{\sigma(q)}^{-1}\ (\lambda_{\sigma(q)})\ F_q\ (\lambda_q)\ ,
\label{eq:rffq}
\end{equation}
where $\lambda_{\sigma(q)} \equiv (\lambda_{\sigma(1)},\ \dots\ ,\ \lambda_{\sigma(n)})$, which generalizes (\ref{eq:rff}).\\

If explicitly constructed, such factorizing $F$-matrices are expected to provide new insights in the investigation of the associated quantum integrable models, more specifically in the context of the Algebraic Bethe Ansatz. Indeed, these factorizing $F$-matrices relate non co-commutative (or non commutative) calculus to co-commutative (or commutative) one.\\

The purpose of this article is to show that in some elementary interesting cases like the Yang-Baxter type algebras associated to the irreducible $N^{th}$ tensor products of the fundamental  representations of the Yangians $Y(sl_2)$ and the quantum affine algebra $U_{q}({\hat {sl_2}})$, i.e., associated to the inhomogeneous $XXX$ and $XXZ$ spin $1 \over 2$ quantum chains of lenght $N$, it is possible to construct explicitly the relevant $F$-matrices satisfying eq. (\ref{eq:rffq}) and to use them in the context of Algebraic Bethe Ansatz.\\

As first applications, we will show that for the inhomogeneous $XXX$ and $XXZ$ spin $1 \over 2$ quantum chains of lenght $N$, the obtained  $F$-matrices can be used to diagonalize the generating matrix of scalar products  
of quantum states of these models and to construct the complete set of eigenvectors of 
the diagonal matrix entries of the $T$-matrix $T(\lambda)$ in the finite 
chain situation.\\ 
It will also appear that these $F$-matrices have very simple algebraic properties leading to a natural basis for the $XXX-{1 \over 2}$ quantum spin chain of lenght $N$ corresponding to co-commutative representations of the monodromy matrix $T(\lambda)$. In such a basis, the operator entries of the $T$-matrix take a form that can be identified as a "diagonal dressing"of the corresponding operators for the Gaudin model \cite{Gaud1}.\\
As a consequence, we will argue that these results should be relevant for  
the computation of correlation functions in the framework developped by Izergin and Korepin \cite{IK1,KBI}.\\

This article is organized as follows :\\

In section $2$ we recall some basic definitions and properties of Drinfel'd twists for quasi-bialgebras, quasi-Hopf algebras, and quasi-triangular quasi-Hopf and Hopf algebras. The case of triangular Hopf algebras leads to the notion of {\em factorizing twists} in section 2.3. In section 2.4 we define the notion of {\em factorizing $F$-matrices} associated to irreducible finite dimensional modules of quantum affine algebras.\\ 
Using these results, we construct, in section $3$,  explicit and compact expressions of such $F$-matrices associated to the irreducible $N^{th}$ tensor products of the fundamental  representations of the Yangians $Y(sl_2)$ and the quantum affine algebra $U_q ({\hat {sl_2}})$, i.e., corresponding to the inhomogeneous XXX and XXZ spin $1 \over 2$ quantum chains of lenght $N$.\\ 
In section $4$, the main algebraic properties of the above constructed $F$-matrices are given, including their inversion, unitary and crossing symmetries. Moreover it is shown that the (operator) diagonal entries of the quantum monodromy matrix $T(\lambda)$ can be diagonalized by means of such $F$-matrices. Similarly, we obtain that the pseudo-vaccuum expectation value of any $T$-matrix product, $T_{q}$, can be diagonalized by using a specific $F$-matrix, namely,
\begin{equation}
< T_{1}\ T_{2}\ \dots\ T_{n} >\ =\ F_{12 \dots n}^{-1}\ {\otimes}_{i = 1}^{n} < T_{i} >\ F_{12 \dots n}\ ,
\label{eq:Ti}
\end{equation}
where $< T_{i} >$ are diagonal matrices depending on the model.\\
Note that $< T_{1}\ T_{2}\ \dots\ T_{n} >$ contains all possible scalar products of quantum states of the model. This formula separates the contributions to these scalar products coming on the one hand from the $R$-matrix (contained in $F$), and on the other hand coming from the specific model (contained in $< T_{i} >$'s).
Finally, very simple difference equations are obtained for these $F$-matrices.\\
In section $5$, we study the inhomogeneous $XXX-{1 \over 2}$ quantum spin chain 
of lenght $N$ in the $F$-basis, i.e., in the basis generated by the column vectors of the inverse of the $F$-matrices. It is shown that, in this basis, the operator entries of the quantum monodromy matrix take  very simple forms, in contrast to the usual basis. In particular, the off-diagonal entries of the quantum monodromy matrix, i.e., the  
$B\ (\lambda)\ =\ T_{12}\ (\lambda)$ and $C\ (\lambda)\ =\ T_{21}\ (\lambda)$ operators, that can be considered as creation and anihilation operators of quantum states of the model, are given as sums of only $N$ elementary quasi-local operators, instead of their usual expressions which are given as sums containing up to $2^{N-1}$ such  terms. Namely, in this basis, we have for example,
\begin{equation}
B\ (\lambda ; {\xi}_{1} \dots {\xi}_{N})\ =\ \sum_{i = 1}^{N} S_{i}^{-}\  {\Delta}_{i} (\lambda ; {\xi}_{1}, \dots, {\xi}_{N})\ ,
\label{eq:BF}
\end{equation} 
where the variables ${\xi}_{i}$ are the inhomogeneity parameters on the spin chain, $S_{i}^{-}$ are usual Pauli matrices ${\sigma}^{-}$ acting at site $(i)$,  and  ${\Delta}_{i}$'s are {\em diagonal} matrices acting on all sites but $(i)$.\\
This leads to a remarkable interpretation of the $XXX-{1 \over 2}$ model as a {\em ``diagonal dressing''} of the corresponding Gaudin model.\\
The use of this $F$-basis also leads to a simplified formula for the $F$-matrix itself.\\
Some conclusions and perspectives, in particular with regards to the computation of correlation functions of these models are also given in the last section.

%% file: sect2
\section{Drinfel'd twists}

The purpose of this section is to give the caracterizing properties of representation theory of Drinfel'd twists associated to finite dimensional irreducible representations of (quasi-triangular) quantum affine algebras and which factorize the corresponding unitary $R$-matrices in these representations as in eqs. (\ref{eq:rffl}, \ref{eq:rffq}). They are given in section 2.4.  Matrices acting in these modules and satisfying this set of properties will be called {\em factorizing $F$-matrices}.\\
Before doing this, we first recall, for completness and also to fix our notations, the basic definitions and properties of Drinfel'd twists for quasi-bialgebras and quasi-Hopf algebras in section 2.1. The case of quasi-triangular Hopf algebras is treated in section 2.2. Then in section 2.3, we describe the case of triangular Hopf algebras. There we define the notion of {\em factorizing twists} relating the quantum non co-commutative co-multiplication $\Delta$ to a co-commutative one $\Delta_0$ and we give their main properties. Then, this is used in section 2.4 as a guideline to define the notion of {\em factorizing $F$-matrices} associated to  finite dimensional irreducible modules of (quasi-triangular) quantum affine algebras whose $R$-matrices are unitary in these representations \cite{CP1,FR1}.

\subsection{General Properties of Drinfel'd twists}

All definitions and results (and their proofs) given in this subsection are taken from the Drinfel'd article \cite{Drin2}. We recall them here only for completness and to fix our notations.

\begin{defi}
By a quasi-bialgebra over $k$ ($k$ being a field of caracteristic zero) is meant a set  ${\cal B}_{\cal A}\ =\ ({\cal A},\ \Delta,\ \epsilon,\ \Phi)$ where ${\cal A}$ is an associative $k$-algebra with unity, $\Delta$ a homomorphism ${\cal A} \rightarrow {\cal A} \otimes {\cal A}$, $\epsilon$ a homomorphism ${\cal A} \rightarrow k$ and $\Phi$ an invertible element of ${\cal A} \otimes {\cal A} \otimes {\cal A}$ such that,
\begin{eqnarray*}
(id \otimes \Delta) \Delta (a)\ &=&\ \Phi\ (\Delta \otimes id) \Delta (a)\ {\Phi}^{-1},\ \ \ a \in {\cal A},\\
(id \otimes id \otimes \Delta) \Phi\  (\Delta \otimes id \otimes id) \Phi\ &=& \ ({\mathbf 1} \otimes \Phi)\  (id \otimes \Delta \otimes id) \Phi\ (\Phi \otimes {\mathbf 1}), \\
(\epsilon \otimes id) \circ \Delta\ &=&\ id\ =\ (id \otimes \epsilon) \circ 
\Delta, \\
(id \otimes \epsilon \otimes id) \Phi\ &=&\ {\mathbf 1}\ .
\end{eqnarray*}
\label{defi:qba}
\end{defi}
We have the following proposition,
\begin{prop}
Suppose given a quasi-bialgebra ${\cal B}_{\cal A}$ and an invertible element $F \in {\cal A} \otimes {\cal A}$ such that $(\epsilon \otimes id)\ F\ =\ (id \otimes \epsilon)\  F\ =\ \mathbf 1$. Put, 
\begin{eqnarray*}
\tilde{\Delta} (a)\ &=&\ F\ {\Delta} (a)\ F^{-1}\ ,\ \ a \in {\cal A}\ , \\
\tilde{\Phi}\ &=&\ ({\mathbf 1} \otimes F)\ (id \otimes \Delta) F\ \Phi\ (\Delta \otimes id) F^{-1}\ (F^{-1} \otimes {\mathbf 1})\ .
\end{eqnarray*}
Then the set $\tilde{\cal B}_{\cal A}\ =\ ({\cal A},\ \tilde{\Delta},\ \epsilon,\ \tilde{\Phi})$ is also a quasi-bialgebra obtained by twisting the set ${\cal B}_{\cal A}$ via the element $F$. The element $F$ itself is called a twist. Twisting via $F_{1}$ and then via $F_{2}$ is equivalent to twisting via $F_{2}\ F_{1}$.
\label{prop:qbaf}
\end{prop}

\begin{defi}
By a quasi-Hopf algebra is meant a quasi-bialgebra ${\cal B}_{\cal A}$ for which there exists $\alpha,\ \beta \in {\cal A}$ and a bijective anti-homomorphism 
$S$ (antipode) of ${\cal A}$ such that,
\begin{eqnarray*}
\sum_i S(b_i)\ \alpha\ c_i\ &=&\ \epsilon (a)\ \alpha\ , \\
\sum_i b_i\ \beta\ S(c_i)\ &=&\ \epsilon (a)\ \beta\ ,
\end{eqnarray*}
for $a \in {\cal A}$ and $\Delta (a)\ =\ \sum_i b_i \otimes c_i$, and, 
\begin{eqnarray*}
\sum_i X_i\ \beta\ S(Y_i)\ \alpha\ Z_i\ &=&\ {\mathbf 1}\ ,\ \ \ \Phi\ =\ \sum_i X_i \otimes Y_i \otimes Z_i\ , \\
\sum_j S(P_j)\ \alpha\ Q_j\ \beta\ S(R_j)\ &=&\ {\mathbf 1}\ ,\ \ {\Phi}^{-1}\ =\ \sum_j P_j \otimes Q_j \otimes R_j\ .
\end{eqnarray*}
\label{defi:qha}
\end{defi}

\begin{prop}
The property of being quasi-Hopf is preserved by twisting.
\label{prop:qhpt}
\end{prop}
Further important properties of quasi-Hopf algebras are given in \cite{Drin2}. Let us now turn to the cases of quasi-triangular and triangular quasi-Hopf algebras, for which is introduced the notion of $R$-matrix.
\begin{defi}
By a quasi-triangular quasi-Hopf algebra is meant a set ${\cal H}_{\cal A}\ =\ ({\cal A},\ R,\ {\Delta},\ \epsilon,\ {\Phi})$ where ${\cal B}_{\cal A}\ =\ ({\cal A},\ {\Delta},\ \epsilon,\ {\Phi})$ is a quasi-hopf algebra and $R \in {\cal A} \otimes {\cal A}$ is an invertible element such that,
\begin{eqnarray*}
R\ \Delta (a)\ &=&\ \sigma \circ \Delta (a)\ R\ ,\ \ a \in {\cal A}\ ,  \\
\sigma\ &:&\ {\cal A} \otimes {\cal A} \rightarrow {\cal A} \otimes {\cal A}\ , \\
\sigma\ &:&\ x \otimes y \mapsto y \otimes x\ ,
\label{eq:rdelta}
\end{eqnarray*}
together with,
\begin{eqnarray*}
(\Delta \otimes id)\ R\ &=&\ {\Phi}_{321}\ R_{13}\ {\Phi}_{132}^{-1}\ R_{23}\ {\Phi}_{123}\ , \\
(id \otimes \Delta)\ R\ &=&\ {\Phi}_{231}^{-1}\ R_{13}\ {\Phi}_{213}\ R_{12}\ {\Phi}_{123}^{-1}\ .
\end{eqnarray*}
It is said triangular if in addition $R_{21}\ R_{12}\ =\ {\mathbf 1}$.\\
The twisting of ${\cal H}_{\cal A}$ is defined similarly as for quasi-Hopf algebras, supplemented with the definition of the twisted $R$-matrix as,
\begin{equation}
\tilde{R}\ =\ F_{21}\ R_{12}\ F_{12}^{-1}\ .
\end{equation}
\label{defi:qtqh}
\end{defi}
Here, for simplicity, we have introduced standard tensorial notations that we will use in all the rest of this article. They are defined as follows.
\begin{defi}
Let ${\cal A}$ be an associative $k$-algebra with unity, and  $X \in {\cal A}^{\otimes\ n}$ be denoted by 
\begin{eqnarray*}
X_{12 \dots n}\ &=&\ \sum_i x^{(1)\ i} \otimes x^{(2)\ i} \otimes \dots \otimes x^{(n)\ i}\ , \\
&=&\ \sum_i x_{1}^{(1)\ i}\ x_{2}^{(2)\ i}\ \dots\ x_{n}^{(n)\ i}\ ,
\end{eqnarray*}
where each element $x_{k}^{(k)\ i}$ acts in the $k^{th}$ copy of ${\cal A}$ of the tensor product ${\cal A}^{\otimes\ n}$ as $x^{(k)\ i}$ and as the identity in all other copies of $\cal A$ in the tensor product. Then we have for any set $(\alpha_1 ,\ \alpha_2 , \dots\ ,\ \alpha_n)$ which is a permutation of the set $(1,\ 2,\ \dots\ ,\ n)$,
\begin{eqnarray*}
X_{\alpha_1 \ \alpha_2 \ \dots\ \alpha_n}\ =\ \sum_i x_{\alpha_1}^{(1)\ i}\ x_{\alpha_2}^{(2)\ i}\ \dots\ x_{\alpha_n}^{(n)\ i}\ ,
\end{eqnarray*}
where each element $x_{\alpha_k}^{(k)\ i}$ acts in the $\alpha_k^{th}$ copy of ${\cal A}$ in the tensor product ${\cal A}^{\otimes\ n}$ as $x^{(k)\ i}$ and as the identity in all other copies of ${\cal A}$ in the tensor product.
\label{defi:xtensor}
\end{defi}
Now, let $\sigma \in S_n$ be an arbitrary element of the permutation group of $n$ elements $S_n$, and $\sigma \ =\ \sigma_{\beta_1}\ \sigma_{\beta_2}\ \dots\ \sigma_{\beta_p}$ a decomposition of it in terms of simple transpositions $\sigma_{\beta_k}$ of $(\beta_k,\ \beta_k + 1)$. We have,
\begin{defi}
Let $\sigma \in S_n$ and arbitrary element $X \in {\cal A}^{\otimes\ n}$. We define the action of $S_n$ on ${\cal A}^{\otimes\ n}$ as,
\begin{eqnarray*}
\sigma\ (X_{12\ \dots\ n})\ =\ X_{\sigma (1)\ \sigma (2)\ \dots\ \sigma (n)}\ .
\end{eqnarray*}
This action is equivalent to the adjoint action of the operator $P_{12\ \dots\ n}^{\sigma}\ =\ P_{\beta_1 ,\ \beta_1 +1}\  \dots\ P_{\beta_p ,\ \beta_p +1}$, $P_{i,\ i+1}$ being the transposition of the $i^{th}$ and $(i + 1)^{th}$ copies in the tensor product ${\cal A}^{\otimes\ n}$, namely,
\begin{eqnarray*}
\sigma\ (X_{12\ \dots\ n})\ &=&\ P_{12\ \dots\ n}^{\sigma}\ X_{12\ \dots\ n}\   (P_{12\ \dots\ n}^{\sigma})^{-1}\ , \\
(P_{12\ \dots\ n}^{\sigma})^{-1}\ &=&\ (P_{12\ \dots\ n}^{{\sigma}^{-1}})\ ,\\
\sigma_i\ (X_{12\ \dots\ n})\ &=&\ P_{i,\ i+1}\ X_{12\ \dots\ n}\ P_{i,\ i+1}\ .
\end{eqnarray*}
\label{defi:sigmax}
\end{defi}
Now, as in the quasi-Hopf case we have the property,
\begin{prop}
The property of being quasi-triangular or triangular quasi-Hopf algebra is preserved by twisting.
\label{prop:qthpt}
\end{prop}
Moreover, the following Theorem shows that any triangular quasi-Hopf Quantum Universal Envelopping Algebra over $k[[h]]$ is a deformation of the trivial Hopf Quantum Universal envelopping algebra over $k[[h]]$.
\begin{theo}[Drinfel'd]
Any triangular quasi-Hopf Quantum Universal Envelopping Algebra over $k[[h]]$ can, by a suitable twist, be brought into the form $R\ = {\mathbf 1}$ and $\Phi\ = {\mathbf 1}$
\label{theo:tqh}
\end{theo}
This ends the list of general definitions and properties of twists we will use in the subsequent sections. We now turn to the Hopf algebra situation which will be relevant to quantum integrable models solvable by the Algebraic Bethe Ansatz.

\subsection{Drinfel'd Twists and Quasi-Triangular Hopf Algebras}

\begin{defi}
A quasi-Hopf algebra ${\cal B}_{\cal A}\ =\ ({\cal A},\ {\Delta},\ \epsilon,\ {\Phi})$ for which $\Phi\ =\ {\mathbf 1}$ is called a Hopf algebra.\\
Similarly, a quasi-triangular (resp. triangular) quasi-hopf algebra ${\cal H}_{\cal A}\ =\ ({\cal A},\ R,\ {\Delta},\ S,\ \epsilon,\ {\Phi})$ with $\Phi\ =\ {\mathbf 1}$ is called a quasi-triangular (resp. triangular) Hopf algebra.
\label{defi:ha}
\end{defi}
Note that if $\Phi\ =\ {\mathbf 1}$, then the co-multiplication $\Delta$ is co-associative, namely,
\begin{equation}
(\Delta \otimes id)\ \Delta (a)\ =\ (id \otimes \Delta)\ \Delta (a)\ ,\ \ \ a \in {\cal A}\ .
\label{eq:del}
\end{equation}
Moreover for quasi-triangular Hopf algebra, the $R$ operator satisfies the Yang-Baxter equation,
\begin{equation}
R_{12} \ R_{13} \ R_{23} \ =\ R_{23} \ R_{13} \ R_{12}\ .
\label{eq:r3}
\end{equation}
The following proposition gives the behaviour of quasi-triangular Hopf algebra under twisting.
\begin{prop}[Drinfel'd]
The property of being quasi-triangular Hopf algebra is preserved by twisting via an invertible element $F\ =\ \sum_i f^i \otimes f_i \in {\cal A} \otimes {\cal A}$ such that $(\epsilon \otimes id)\ F\ =\ (id \otimes \epsilon)\  F\ =\ \mathbf 1$ and satisfying the cocycle condition,
\begin{equation}
(F \otimes {\mathbf 1})\ (\Delta \otimes id)\ F\ =\ ({\mathbf 1} \otimes F)\ (id \otimes \Delta)\ F\ .
\end{equation}
Moreover, $u\ =\ \sum_i f^i\ S(f_i)$ is invertible and the twisted antipode is given by $\tilde{S} (a)\ =\ u\ S(a)\ u^{-1}$, while the co-multiplication, $R$-matrix and co-unit change according to the same rule as in the quasi-triangular quasi-Hopf situation. Such a twist $F$ will be called an admissible twist.\\
Moreover, if $F$ is an admissible twist, then $G_{12}\ =\ F_{21}\ R_{12}$ is also an admissible twist with $\tilde{\Delta}^{F} (a)\ =\ \sigma \circ \tilde{\Delta}^{G} (a)$.
\label{prop:qthf}
\end{prop} 
We now introduce some convenient notations for the co-multiplication action for any element $X \in {\cal A} \otimes {\cal A}$, $X\ =\ X_{12}\ =\ \sum_i a_i \otimes b_i$,  as,
\begin{eqnarray}
({\Delta}^{(n-2)} \otimes id)\ X\ &=&\ X_{1 \dots n-1,\ n}\ ,\nonumber \\
(id \otimes {\Delta}^{(n-2)})\ X\ &=&\ X_{1,\ 2 \dots n}\ ,
\label{eq:notd}
\end{eqnarray}
where we have for any integer $p$,
\begin{eqnarray*}
{\Delta}^{(p)} (a)\ =\ ({\Delta}^{(p-1)} \otimes id)\ \Delta (a)\ =\ (id \otimes {\Delta}^{(p-1)})\ \Delta (a),\ \ a \in {\cal A}\ .
\end{eqnarray*}
For later use, consider now the element $\tilde{X}_{12}\ =\ X_{21}\ =\ \sum_i b_i \otimes a_i$. We have the following co-multiplication action on it,
\begin{eqnarray}
({\Delta}^{(n-2)} \otimes id)\ \tilde{X}\ &=&\ \tilde{X}_{1 \dots n-1,\ n}\ , \nonumber \\
&=&\ X_{n,\ 1 \dots n-1}\ , 
\label{eq:tdlx}
\end{eqnarray}
and similarly,
\begin{eqnarray}
(id \otimes {\Delta}^{(n-2)})\ \tilde{X}\ &=&\ \tilde{X}_{1,\ 2 \dots n}\ , \nonumber \\
&=&\ X_{2 \dots n,\ 1}\ .
\label{eq:tdrx}
\end{eqnarray}
We have the following proposition,
\begin{prop}
Let $F$ be an admissible twist for a given Hopf algebra ${\cal H}_{\cal A}\ =\ ({\cal A},\ R,\ {\Delta},\ \epsilon)$, and let the twisted co-multiplication be given as $\tilde{\Delta} (a)\ =\ F\ {\Delta} (a)\ F^{-1}\ ,\ \ a \in {\cal A}$. Define for any integer $n$ the ``total'' twist $F_{12 \dots n} \in {\cal A}^{\otimes n}$ as,
\begin{equation}
F_{12 \dots n}\ =\ F_{12}\ F_{12,\ 3}\ \dots\ F_{12 \dots n-1,\ n}\ =\ F_{12 \dots n-1}\ F_{12 \dots n-1,\ n}\ .
\label{eq:tfl}
\end{equation}
It satisfies,
\begin{equation}
F_{12 \dots n}\ {\Delta}^{(n-1)} (a)\ F_{12 \dots n}^{-1}\ =\ \tilde{\Delta}^{(n-1)} (a)\ ,\ \  a \in {\cal A}\ .
\label{eq:ddt}
\end{equation}
Moreover we have the cocycle relation,
\begin{equation}
F_{1,\ 2 \dots n-1}\ F_{12 \dots n-1,\ n}\ =\ F_{2\ \dots\ n-1,\ n}\ F_{1,\ 23 \dots n}\ .
\label{eq:crf}
\end{equation}
Hence we obtain the following alternative expression for the total twist $F_{12 \dots n}$ ,
\begin{equation}
F_{12 \dots n}\ =\ F_{n-1\ n}\ F_{n-2,\ n-1\ n}\ \dots\ F_{1,\ 23 \dots n}\ =\ F_{23 \dots n}\ F_{1,\ 23 \dots n}\ .
\label{eq:tfr}
\end{equation}
\label{prop:tf}
\end{prop}
\Proof
First we prove (\ref{eq:ddt}) by induction on $n$. It is true by definition for $n\ =\ 2$. Suppose it is true up to some  $n \geq 2$. Then we have for any $a \in {\cal A}$,
\begin{eqnarray*}
{\Delta}^{(n)} (a)\ &=&\ ({\Delta}^{(n-1)} \otimes id)\ \Delta (a)\ , \\
&=&\ ({\Delta}^{(n-1)} \otimes id)\ (F^{-1}\ \tilde{\Delta} (a)\ F)\ , \\
&=&\ ({\Delta}^{(n-1)} \otimes id)\ (F^{-1})\ ({\Delta}^{(n-1)} \otimes id)\  \tilde{\Delta} (a)\ ({\Delta}^{(n-1)} \otimes id)\ (F)\ .
\end{eqnarray*}
Then applying the relation (\ref{eq:ddt}) for $n$ we have,
\begin{equation}
{\Delta}^{(n-1)} (a)\ =\ F_{12 \dots n}^{-1}\ \tilde{\Delta}^{(n-1)} (a)\ F_{12 \dots n}\ ,
\end{equation}
with,
\begin{equation}
F_{12 \dots n}\ =\ F_{12}\ F_{12,\ 3}\ \dots\ F_{12 \dots n-1,\ n}\ =\ F_{12 \dots n-1}\ F_{12 \dots n-1,\ n}\ .
\end{equation}
So we have,
\begin{eqnarray*}
{\Delta}^{(n)} (a)\ &=&\ F_{12 \dots n,\ n+1}^{-1}\ F_{12 \dots n}^{-1}\ \tilde{\Delta}^{(n)} (a)\ F_{12 \dots n}\ F_{12 \dots n,\ n+1}\ , \\
&=&\ F_{12 \dots n+1}^{-1}\ \tilde{\Delta}^{(n)} (a)\ F_{12 \dots n+1}\ ,
\end{eqnarray*}
which prove the first part of the proposition.
Then recall that admissible twists satisfy the cocycle relation,
\begin{equation}
F_{23}\ F_{1,23}\ =\ F_{12}\ F_{12,3}. 
\end{equation}
Applying the operation $(id \otimes {\Delta}^{(n-3)} \otimes id)$ to it we get exactly (\ref{eq:crf}).\\
Let 
\begin{equation}
G_{12 \dots n}\ =\ F_{n-1\ n}\ F_{n-2,\ n-1\ n}\ \dots\ F_{1,\ 23 \dots n}\ =\ G_{23 \dots n}\ F_{1,\ 23 \dots n}\ ,
\end{equation}
then obviously, $G_{123}\ =\ F_{123}$ by the elementary cocycle relation. Suppose it is so up to $(n-1)$ spaces, namely, $G_{12 \dots p}\ =\ F_{12 \dots p},\ \ p \leq n-1$. Then we have the following chain of identities,
\begin{eqnarray*}
G_{12 \dots n}\ &=&\ G_{23 \dots n}\ F_{1,\ 23 \dots n}\ =\ F_{23 \dots n}\ F_{1,\ 23 \dots n}\ , \\
&=&\ F_{23 \dots n-1}\  F_{23 \dots n-1,\ n}\ F_{1,\ 23 \dots n}\ , \\
&=&\ F_{23 \dots n-1}\  F_{1,\ 23 \dots n-1}\ F_{12 \dots n-1,\  n}\ , \\
&=&\ G_{23 \dots n-1}\  F_{1,\ 23 \dots n-1}\ F_{12 \dots n-1,\  n}\ , \\
&=&\ G_{12 \dots n-1}\ F_{12 \dots n-1,\  n}\ , \\
&=&\ F_{12 \dots n-1}\ F_{12 \dots n-1,\  n}\ , \\
&=&\ F_{12 \dots n}\ ,
\end{eqnarray*}
where we used first the induction hypothesis, then the definition of the total $F$, then the cocycle relation, then again the induction hypothesis and finally the definition of the total $F$.
\QED

\subsection{Drinfel'd Twists and Triangular Hopf Algebras}

We are now interested in the case of triangular Hopf algebras. We have the following proposition,
\begin{prop}
Let ${\cal H}_{\cal A}\ =\ ({\cal A},\ R,\ S,\ {\Delta},\ \epsilon)$ be a triangular Hopf algebra. Then for any integer $n$ and for any element $\sigma$ in the permutation group $S_n$ decomposed in a minimal way in terms of simple transpositions as $\sigma \ =\ \sigma_{\beta_1}\ \sigma_{\beta_2}\ \dots\ \sigma_{\beta_p}$, where $\sigma_k$ interchange the $k^{th}$ and the $(k + 1)^{th}$ copies of ${\cal A}$ in the tensor product ${\cal A}^{\otimes\ n}$, we can define a map from $S_n$ to ${\cal A}^{\otimes\ n}$ which associate in a unique way an  element $R_q^{\sigma}\ =\ R_{12 \dots n}^{\sigma} \in {\cal A}^{\otimes n}$ to any element $\sigma \in S_n$, $q\ =\ (12\ \dots\ n)$, and  satisfying,
\begin{equation}
R_{12 \dots n}^{\sigma}\ {\Delta}^{(n-1)} (a)\ =\ {\sigma} \circ {\Delta}^{(n-1)} (a)\ R_{12 \dots n}^{\sigma}\ .
\label{eq:rs}
\end{equation}
It is given as,
\begin{eqnarray*}
R_q^{\sigma}\ &=&\ P_q^{\sigma}\ \hat{R}_q^{{\sigma}^{-1}}\ , \\
\hat{R}_q^{{\sigma}}\ &=&\ \hat{R}_q^{{\sigma}_{\beta_1}}\ \hat{R}_q^{{\sigma}_{\beta_2}}\ \dots\ \hat{R}_q^{{\sigma}_{\beta_p}}\ , \\
P_q^{\sigma}\ &=&\ P_{\beta_1 ,\ \beta_1 +1}\ P_{\beta_2 ,\ \beta_2 +1}\ \dots\ P_{\beta_p ,\ \beta_p +1}\ , 
\end{eqnarray*}
where $\hat{R}_q^{{\sigma}_{i}}\ =\ P_{i, i +1}\ R_{i, i+1}$, $R_{i, i+1}$ being  the $R$ operator and $P_{i,\ i+1}$ being the transposition operator of the $i^{th}$ and $(i + 1)^{th}$ copies in the tensor product ${\cal A}^{\otimes\ n}$.\\
Moreover $\hat{R}_q^{{\sigma}}$ does not depend on the chosen decomposition of the element $\sigma \in S_n$. It gives a representation of the permutation group $S_n$ in ${\cal A}^{\otimes\ n}$ and the composition law for the operators $R$  are given by,
\begin{eqnarray}
\hat{R}_q^{{\sigma}_1 {\sigma}_2 }\ &=&\ \hat{R}_q^{{\sigma}_1}\ \hat{R}_q^{{\sigma}_2}\ , \label{eq:hatrss}\\
R_q^{{\sigma}_1 {\sigma}_2 }\ &=&\ R_{{\sigma}_2 (q)}^{\sigma_1}\  R_q^{{\sigma}_2}\ .\label{eq:rss}
\end{eqnarray}
In particular let us define, 
\begin{eqnarray*}
({\Delta}^{(n-2)} \otimes id)\ R\ &=&\ R_{1\ \dots\ n-1,\ n}\ , \\
(id \otimes {\Delta}^{(n-2)})\ R\ &=&\ R_{1,\ 2\ \dots\ n}\ ,
\end{eqnarray*}
we have,
\begin{eqnarray}
R_{1\ \dots\ n}^{\sigma}\ R_{1\ \dots\ n,\ 0}\ =\ R_{\sigma(1)\ \dots\ \sigma(n),\ 0}\ R_{1\ \dots\ n}^{\sigma}\ , 
\label{eq:rsr0q1}\\
R_{1\ \dots\ n}^{\sigma}\ R_{0,\ 1\ \dots\ n}\ =\ R_{0,\ \sigma(1)\ \dots\ \sigma(n)}\ R_{1\ \dots\ n}^{\sigma}\ .
\label{eq:rsr0q2}
\end{eqnarray}
\label{prop:rsigma}
\end{prop}
\Proof
Recall first that $\hat{R}_q^{{\sigma}}$ as defined above provides a representation of the permutation group $S_n$ due to the Yang-Baxter equation and unitarity for the operator  $R$. Hence to any element $\sigma \in S_n$ we can associate a unique operator $\hat{R}_q^{{\sigma}}$. In particular it doesn't depend on the particular decomposition of the element $\sigma$ in terms of simple transpositions. The same is true for the permutation operator $P_q^{{\sigma}}$. 
Hence $R_q^{\sigma}$ is also uniquely defined.\\ 
Recall also the fact that $\hat{R}_q^{{\sigma}_{i}}\ =\ P_{i, i +1}\ R_{i, i+1}$ commutes with ${\Delta}^{(n-1)}$. This follows from the Definition (\ref{defi:qtqh}) rewritten as,
\begin{eqnarray*}
P_{12}\ R_{12}\ \Delta (a)\ =\ \Delta (a)\ P_{12}\ R_{12}\ ,\ \ a \in 
{\cal A}\ ,
\end{eqnarray*}
and the co-associativity of the co-multiplication $\Delta$.\\
Hence we have for any permutation $\sigma \in S_n$, with $q \equiv (1\ \dots\ n)$,
\begin{eqnarray*}
\hat{R}_q^{{\sigma}}\ \Delta^{(n-1)} (a)\ =\ \Delta^{(n-1)} (a)\ \hat{R}_q^{{\sigma}}\ ,\ \ a \in {\cal A}\ .
\end{eqnarray*}
It leads to the desired relation for $R_q^{\sigma}\ =\ P_q^{\sigma}\ \hat{R}_q^{{\sigma}^{-1}}$. Equations (\ref{eq:rsr0q1}, \ref{eq:rsr0q2}) follows from the definition. Equation (\ref{eq:rss}) follows from eq. (\ref{eq:hatrss}). \QED
For example let us consider the cyclic permutation ${\sigma}_c$,
\begin{equation}
{\sigma}_c (X_{12\ \dots\ n})\ =\ X_{2\ \dots\ n1}\ ,
\label{eq:sc}
\end{equation}
we have, $q\ =\ (12\ \dots\ n)$,
\begin{eqnarray}
R_q^{{\sigma}_c}\ &=&\ R_{1n}\ \dots\ R_{12}\ =\ R_{1,\ 23\ \dots\ n}\ , \nonumber\\ 
&=&\ (id \otimes {\Delta}^{(n-2)})\ R\ ,
\label{eq:rsc1}
\end{eqnarray}
and similarly,
\begin{eqnarray}
R_q^{{\sigma}_c^{-1}}\ &=&\ R_{1\ n}\ \dots\ R_{n-1\ n}\ =\ R_{12\ \dots\ n-1,\ n}\ , \nonumber\\
&=&\ ({\Delta}^{(n-2)} \otimes id)\ R\ .
\label{eq:rsc2}
\end{eqnarray}
We now define the notion of factorizing twist $F$.
\begin{defi}
Let ${\cal H}_{\cal A}\ =\ ({\cal A},\ R,\ S,\ {\Delta},\ \epsilon)$ be a triangular Hopf algebra, and $F$ an admissible twist such that $\tilde{R}\ = {\mathbf 1}$. Then $R\ =\ F_{21}^{-1}\ F_{12}$ and the deformed co-multiplication $\tilde{\Delta}$ is co-commutative.\\
Such a twist will be called ``factorizing'' twist.
\label{defi:fact}
\end{defi}
We have the following fundamental property of ``factorizing''twists.
\begin{prop}
Let ${\cal H}_{\cal A}\ =\ ({\cal A},\ R,\ S,\ {\Delta},\ \epsilon)$ be a triangular Hopf algebra, and $F$ a ``factorizing'' twist for it. For any integer $n$, let us consider the total twist $F_q\ =\ F_{12\ \dots\ n}$,  $q\ =\ (12\ \dots\ n)$. For any element $\sigma \in S_n$ it satisfies,
\begin{equation}
F_{{\sigma}(q)}\ R_q^{\sigma}\ =\ F_q\ ,
\label{eq:frfs}
\end{equation}
where $F_{{\sigma}(q)}\ =\ F_{{\sigma}(1)\ \dots\ {\sigma}(n)}$. Equivalently, we also have,
\begin{equation}
\hat{R}_q^{\sigma}\ =\ F_q^{-1}\ P_q^{\sigma}\ F_q\ .
\label{eq:rfpf}
\end{equation}
This relation is directly related to the co-commutativity of the deformed co-multiplication by $F$.
In particular we have,
\begin{eqnarray}
R_q^{\sigma}\ F_{0,\ q}\ &=&\ F_{0,\ \sigma (q)}\ R_q^{\sigma}\ , 
\label{eq:rpfa} \\
R_q^{\sigma}\ F_{q,\ 0}\ &=&\ F_{\sigma (q),\ 0}\ R_q^{\sigma}\ , 
\label{eq:rpfb} \\
F_{q,\ 0}\ R_{0,\ q}\ &=&\ F_{0,\ q}\ , 
\label{eq:rpfc} \\
F_{0,\ q}\ R_{q,\ 0}\ &=&\ F_{q,\ 0}\ ,
\label{eq:rpfd}
\end{eqnarray}
where we have used the notations (\ref{eq:notd}) for the co-multiplication actions on $F$ and $R$, and $(0)$ is an additional copy of ${\cal A}$. Moreover, using the co-multiplication $\tilde{\Delta} (a)\ =\ F_{12}\ \Delta (a)\ F_{12}^{-1}$, we can define $\tilde{F}_{q,\ 0}\ =\ F_q\ F_{q,\ 0}\ F_q^{-1}$ and similarly, $\tilde{F}_{0,\ q}\ =\ F_q\ F_{0,\ q}\ F_q^{-1}$. They verify,
\begin{eqnarray*}
F_{12\ \dots\ n}\ &=&\ \tilde{F}_{1\ \dots\ n-1,\ n}\ F_{12\ \dots\ n-1}\ , \\
&=&\ \tilde{F}_{1\ \dots\ n-1,\ n}\ \tilde{F}_{1\ \dots\ n-2,\ n-1}\ \dots\ \tilde{F}_{12,\ 3}\ F_{12}\ ,
\end{eqnarray*}
and,
\begin{eqnarray*}
F_{12\ \dots\ n}\ &=&\ \tilde{F}_{1,\ 23\ \dots\ n}\ F_{23\ \dots\ n}\ , \\ 
&=&\ \tilde{F}_{1,\ 23\ \dots\ n}\ \tilde{F}_{2,\ 3\ \dots\ n}\ \dots\ \tilde{F}_{n-2,\ n-1n}\ F_{n-1 n}\ ,
\end{eqnarray*}
and hence we have for any permutation $\sigma \in S_n$, $q\ =\ (12\ \dots\ n)$,  \begin{eqnarray}
\tilde{F}_{0,\ q}\ &=&\ \tilde{F}_{0,\ {\sigma}(q)}\ , \label{eq:fts1}\\
\tilde{F}_{q,\ 0}\ &=&\ \tilde{F}_{{\sigma}(q),\ 0}\ .
\label{eq:fts2}
\end{eqnarray}
\label{prop:fundf}
\end{prop}
\Proof
We first prove eqs. (\ref{eq:rpfa} - \ref{eq:rpfd}). 
The proof is straightforward from the formula $R\ =\ F_{21}^{-1}\ F_{12}$ by applying to it various co-multiplication action on the right or on the left to the required power, and the use of the Proposition (\ref{prop:rsigma}) and relations (\ref{eq:tdlx}, \ref{eq:tdrx}).\\
Let us now prove equation (\ref{eq:frfs}) in the case of the transposition of the first two elements ${\sigma}_{12}$ and for the cyclic permutation ${\sigma}_c$. This follows from eqs. (\ref{eq:rsc1}, \ref{eq:rpfb}) and the definition of the total twist $F_{12\ \dots\ n}$ (\ref{eq:tfl}, \ref{eq:tfr}). Indeed, we have for the transposition ${\sigma}_{12}$, using eq. (\ref{eq:rpfb}) we just proved,
\begin{eqnarray*}
F_{213\ \dots\ n}\ R_{12}\ &=&\ F_{21}\ F_{21,\ 3}\ \dots\ F_{21\ \dots\ n-1,\ n}\ R_{12}\ , \\
&=&\ F_{21}\ F_{21,\ 3}\ \dots\ F_{21\ \dots\ n-2,\ n-1}\ R_{12}\ F_{12\ \dots\ n-1,\ n}\ , \\
&=&\ F_{21}\ R_{12}\ F_{12,\ 3}\ \dots\ F_{12\ \dots\ n-1,\ n}\ , \\
&=&\ F_{12}\ F_{12,\ 3}\ \dots\ F_{12\ \dots\ n-1,\ n}\ =\ F_{123\ \dots\ n}\ .
\end{eqnarray*}
For the cyclic permutation ${\sigma}_c$ we have similarly from eqs. (\ref{eq:rsc1}, \ref{eq:tdrx}), and the formula $R\ =\ F_{21}^{-1}\ F_{12}$, 
\begin{eqnarray*}
F_{1,\ 23\ \dots\ n}\ &=&\ F_{23\ \dots\ n,\ 1}\ R_{1,\ 23\ \dots\ n}\ , \\ 
&=&\ F_{23\ \dots\ n,\ 1}\ R_{12\ \dots\ n}^{{\sigma}_c}\ . 
\end{eqnarray*}
Now multiplying by $F_{23\ \dots\ n}$ on the left we obtain,
\begin{eqnarray*}
F_{23\ \dots\ n}\ F_{1,\ 23\ \dots\ n}\ &=&\ F_{23\ \dots\ n}\ F_{23\ \dots\ n,\ 1}\ R_{12\ \dots\ n}^{{\sigma}_c} \ .
\end{eqnarray*}
Hence leading to,
\begin{eqnarray*}
F_{12\ \dots\ n}\ &=&\ F_{23\ \dots\ n1}\ R_{12\ \dots\ n}^{{\sigma}_c}\ , \\
&=&\ {\sigma}_c (F_{12\ \dots\ n})\ R_{12\ \dots\ n}^{{\sigma}_c}\ ,
\end{eqnarray*}
which is the desired relation.\\
So, for $\hat{R}_q^{\sigma}$, associated to the transposition of the first two elements ${\sigma}_{12}$ and for the cyclic permutation ${\sigma}_c$,
\begin{eqnarray*}
\hat{R}_q^{\sigma}\ &=&\ P_q^{\sigma}\ R_q^{{\sigma}^{-1}}\ , \\
&=&\ P_q^{\sigma}\ F_{{\sigma}^{-1}(q)}^{-1}\ F_q\ , \\
&=&\ F_q^{-1}\ P_q^{\sigma}\ F_q\ .
\end{eqnarray*}
Now we use the well known fact that the permutation group $S_n$ is generated by the cyclic permutation ${\sigma}_c$ and the transposition of the first two elements ${\sigma}_{12}$. Using this, together with the fact that $\hat{R}_q^{\sigma}$ gives a representation of $S_n$, shows that any such $\hat{R}_q^{\sigma}$ operator can be decomposed as an ordered product of $\hat{R}_q^{\sigma}$'s associated to ${\sigma}_{12}$'s and ${\sigma}_c$'s. It follows, $P_q^{\sigma}$ being also a representation of $S_n$, that the 
relations (\ref{eq:rfpf}) and (\ref{eq:frfs}) are true for any permutation $\sigma \in S_n$. \\
Finally, eqs. (\ref{eq:fts1}, \ref{eq:fts2}) follow from eqs. (\ref{eq:rpfa}, \ref{eq:rpfb}, \ref{eq:frfs}).
\QED

\subsection{Representation theory for admissible twists : F-matrices}

We will now give the caracterizing properties of representation theory of a particular class of Drinfel'd twists associated to finite dimensional irreducible evaluation representations of (quasi-triangular) quantum affine algebras and which factorize the corresponding (unitary) $R$-matrices in these representations as in eqs. (\ref{eq:rffl}, \ref{eq:rffq}). Matrices acting in these modules and satisfying this set of properties will be called {\em factorizing $F$-matrices}.\\
It is the unitarity of the $R$-matrix which allows to introduce the notion of factorizing $F$-matrices in much the same way as the notion of factorizing twists is defined for triangular Hopf algebras in section 2.3.\\

Let us now fix some notations. Following \cite{FR1}, let $V$ be a finite dimensional irreducible ${\cal U}_q (\hat{\cal G})$ module with ${\pi}_V : {\cal U}_q (\hat{\cal G}) \rightarrow End(V)$. Using the element $d$ of ${\cal U}_q (\hat{\cal G})$ it is possible to define an automorphism $D_z$ of ${\cal U}_q (\hat{\cal G}) \otimes {\mathbf C}[z,\ z^{-1}]$ and a new representation, 
$ {\pi}_{V (z)} : {\cal U}_q (\hat{\cal G}) \rightarrow End(V) \otimes {\mathbf C}[z,\ z^{-1}]$ 
by the formula, 
${\pi}_{V (z)} (a)\ =\ {\pi}_V (D_z (a)),\ \ \ a \in {\cal U}_q (\hat{\cal G})$.\\
Hence, from the finite dimensional module $V$, we get a one-parameter family of finite dimensional modules $V(z)$ connected by the action of the automorphism $D_z$ (with $V\ =\ V(1)$).\\
Let $R(z)$ be the corresponding $R$-matrix satisfying the Yang-Baxter equation,
\begin{equation}
R_{12} (z)\ R_{13} (zw)\ R_{23} (w)\ =\ R_{23} (w)\ R_{13} (zw)\ R_{12} (z)\ ,
\label{eq:rrrzw}
\end{equation}
and the co-multiplication relation,
\begin{equation}
R^{VW} (z)\ {\pi}_{V(zw) \otimes W(w)} (\Delta (a))\ =\ {\pi}_{V(zw) \otimes W(w)} (\sigma \circ \Delta (a))\ R^{VW} (z)\ .
\label{eq:dz}
\end{equation}
It is shown in \cite{CP1,FR1} that, in the irreducible finite dimensional case, this $R$-matrix is unitary,
\begin{equation}
R_{21} (z^{-1})\ R_{12} (z)\ =\ {\mathbf 1}\ .
\label{eq:ur}
\end{equation}
For any integer $n$ let us consider an irreducible tensor product of $n$ finite dimensional modules $V_i (z_i )$, all being connected by the automorphism $D_z$,  $V_q\ =\ V_1 (z_1 ) \otimes V_2 (z_2 ) \otimes \dots \otimes V_n (z_n )$. Elements of $End (V_q)$ are denoted as $X_{12 \dots n} (z_{1} ,\ z_{2} ,\ \dots,\ z_n )$. Moreover, to simplify formulas, we will not write explicitly the $z_i$'s dependence everywhere. So, unless another dependence on the $z_i$'s is specified, we  assume the following identification $ X_{12 \dots n} (z_{1} ,\ z_{2} ,\ \dots,\ z_n ) \equiv X_{12 \dots n}$. It means that to any module $V_i$ is associated the parameter $z_i$.
\begin{defi}
Let $n$ be an integer and $\sigma \in S_n$ an arbitrary element of the permutation group with $n$ elements. Let $X \in End(V_q )$, $V_q\ =\ V_1 (z_1 ) \otimes V_2 (z_2 ) \otimes \dots \otimes V_n (z_n )$. We define the extended  action of the symmetry group $S_n$ on $End (V_q )$ as,
\begin{equation}
\sigma ( X_{12\ \dots\ n}\ (z_{1} ,\ z_{2} ,\ \dots,\ z_n ) )\ =\ X_{\sigma (1) \sigma (2) \dots \sigma (n)}\ (z_{\sigma (1)} ,\ z_{\sigma (2)} ,\ \dots,\ z_{\sigma (n)} )\ .
\label{eq:sigmaz}
\end{equation}
Let $\sigma \in S_n$ be decomposed in a minimal way into an ordered product of simple transpositions as, $\sigma \ =\ {\sigma}_{{\beta}_1}\ \dots\ {\sigma}_{{\beta}_p}$, then we can define, 
${\cal P}_q^{\sigma}\ =\ {\cal P}_{\beta_1 ,\ \beta_1 +1}\  \dots\ {\cal P}_{\beta_p ,\ \beta_p +1}$ where ${\cal P}_{i\ i+1}$ 
is the permutation operator in $V_i (z_i ) \otimes V_{i+1} (z_{i+1} )$ acting on any element $X \in End (V_q )$ as, 
\begin{eqnarray*}
{\cal P}_{i\ i+1}\ X_{1\ \dots\ i\ i+1\ \dots\ n} (z_{1} ,\ \dots,\ z_{i},\ z_{i+1},\ \dots,\ \ z_n ) \ {\cal P}_{i\ i+1}\ =\  \\
=\ X_{1\ \dots\ i+1\ i\ \dots\ n} (z_{1} ,\ \dots,\ z_{i+1},\ z_{i},\ \dots,\ \ z_n )\ .
\end{eqnarray*}
\label{defi:sigmaxz}
\end{defi}
Here is now the well known  analogue of Proposition (\ref{prop:rsigma}) (see \cite{CP1,FR1}).
\begin{prop}
Let $R$ be an $R$-matrix, associated to irreducible finite dimensional evaluation representation of a quantum affine algebra, satisfying Yang-Baxter equation and unitarity relation. Then, for any integer $n$ and to any element $\sigma$ in the permutation group $S_n$ decomposed in a minimal way in terms of simple transpositions as $\sigma \ =\ \sigma_{\beta_1}\ \sigma_{\beta_2}\ \dots\ \sigma_{\beta_p}$, we can associate in a unique way an element $R_q^{\sigma}\ =\ R_{12 \dots n}^{\sigma} \in End (V_q )$, $V_q\ =\ V_1 (z_1 ) \otimes V_2 (z_2 ) \otimes \dots V_n (z_n )$. It is given as,
\begin{eqnarray*}
R_q^{\sigma}\ &=&\ {\cal P}_q^{\sigma}\ \hat{\cal R}_q^{{\sigma}^{-1}}\ , \\
\hat{\cal R}_q^{{\sigma}}\ &=&\ \hat{\cal R}_q^{{\sigma}_{\beta_1}}\ \hat{\cal R}_q^{{\sigma}_{\beta_2}}\ \dots\ \hat{\cal R}_q^{{\sigma}_{\beta_p}}\ , \\
{\cal P}_q^{\sigma}\ &=&\ {\cal P}_{\beta_1 ,\ \beta_1 +1}\ {\cal P}_{\beta_2 ,\ \beta_2 +1}\ \dots\ {\cal P}_{\beta_p ,\ \beta_p +1}\ , 
\end{eqnarray*}
where $\hat{\cal R}_q^{{\sigma}_{i}}\ =\ {\cal P}_{i, i +1}\ R_{i, i+1} (z_i, z_{i+1})$, $R_{i, i+1} (z_i, z_{i+1})$ being  the $R$-matrix acting in $V_i (z_i ) \otimes V_{i+1} (z_{i+1} )$ as $R$ and as the identity in all other modules in the tensor product $V_q$.\\
Moreover $\hat{\cal R}_q^{{\sigma}}$ gives a representation of the extended permutation group action $S_n$ in $End (V_q )$. In particular for the cyclic permutaion  ${\sigma}_c$,
\begin{equation}
{\sigma}_c ( X_{12\ \dots\ n} ( z_{1} ,\ \dots,\ z_n ) )\ =\ X_{2\ \dots\ n1} ( z_{2} ,\ \dots,\ z_{n},\ z_{1} )\ ,
\end{equation}
we have,
\begin{eqnarray}
R_{1\ \dots\ n}^{{\sigma}_c} ( z_{1} ,\ \dots,\ z_n )\ =\ R_{1n} (z_1 ,\ z_n )\ \dots\ R_{12} ( z_1 ,\ z_2 )\ =\ R_{1,\ 23\ \dots\ n} ( z_{1} ;\ z_{2} \dots,\ z_{n})\ ,  
\label{eq:rsc1z}
\end{eqnarray}
and similarly,
\begin{eqnarray}
R_{1\ \dots\ n}^{{\sigma}_c^{-1}}\ =\ R_{1\ n} (z_1 ,\ z_n )\ \dots\ R_{n-1\ n} (z_{n-1} ,\ z_n )\ =\ R_{12\ \dots\ n-1,\ n} ( z_{1},\ \dots,\ z_{n-1};\ z_{n} )\ .
\label{eq:rsc2z}
\end{eqnarray}
Moreover, it satisfies, for any element $\sigma \in S_n$,
\begin{eqnarray}
R_q^{\sigma}\ R_{0,\ q}\ &=&\ R_{0,\ \sigma(q)}\ R_q^{\sigma}\ ,  \label{eq:r0qsa}\\
R_q^{\sigma}\ R_{q,\ 0}\ &=&\ R_{\sigma(q),\ 0}\ R_q^{\sigma}\ .
\label{eq:r0qsb}
\end{eqnarray}
\label{prop:rsigmaz}
\end{prop}
\begin{defi}
By factorizing $F$-matrices is meant a set of invertible elements $F \in End (V_q )$ , associated to irreducible finite dimensional modules $V_q\ =\ V_1 (z_1 ) \otimes V_2 (z_2 ) \otimes \dots \otimes V_n (z_n )$, defined for any integer $n$ and such that for any element $\sigma \in S_n$,
\begin{equation}
F_{{\sigma}(1)\ \dots\ {\sigma}(n)}\ (z_{\sigma (1)} ,\ \dots,\ z_{\sigma (n)} )\ R_{1\ \dots\ n}^{\sigma}\ (z_{1} ,\ \dots,\ z_n )\ =\ F_{1\ \dots\ n}\ (z_{1} ,\ \dots,\ z_n )\ ,
\label{eq:frz}
\end{equation}
or in compact notations, $F_{{\sigma}(q)}\ R_q^{\sigma}\ =\ F_q$, with $q\ =\ (12\ \dots\ n)$.
\label{defi:factfm}
\end{defi}
Now the caracterizing properties of these ``factorizing'' $F$-matrices is given in the following proposition to be used when constructing explicit examples in the next section.
\begin{prop}
For any integer $n$, let $F_{1\ \dots\ n} (z_{1} ,\ \dots,\ z_n )$ be factorizing $F$-matrices. We define ``partial'' $F$-matrices through the formulas,
\begin{eqnarray*}
F_{1,\ 2\ \dots\ n}\ (z_{1} ;\ z_{2} \dots,\ z_{n})\ &=&\ F_{2\ \dots\ n}^{-1}\  (z_{2},\ \dots,\ z_{n})\ F_{1\ \dots\ n}\ (z_{1} ,\ \dots,\ z_n )\ , \\
F_{1\ \dots\ n-1,\ n}\ (z_{1},\ \dots,\ z_{n-1};\ z_{n})\ &=&\ F_{1\ \dots\ n-1}^{-1}\ (z_{1},\ \dots,\ z_{n-1})\ F_{1\ \dots\ n}\ (z_{1} ,\ \dots,\ z_n )\ .
\end{eqnarray*}
They satisfy,
\begin{eqnarray}
F_{2\ \dots\ n-1,\ n}\ (z_{2} ,\ \dots,\ z_{n-1} ;\ z_{n})\ F_{1,\ 2\ \dots\ n}\  (z_{1} ;\ z_{2} \dots,\ z_{n})\ = \nonumber \\ 
=\ F_{1,\ 2\ \dots\ n-1}\ (z_{1} ;\ z_{2} \dots,\ z_{n-1})\ F_{1\ \dots\ n-1,\ n}\ (z_{1},\ \dots,\ z_{n-1};\ z_{n})\ , 
\label{eq:cocyclez}
\end{eqnarray}
\begin{eqnarray}
R_{1\ \dots\ n}^{\sigma}\ (z_{1} ,\ \dots,\ z_n )\ F_{0,\ 1\ \dots\ n}\ (z_{0};\ z_{1} ,\ \dots,\ z_n )\ = \nonumber \\ 
=\ F_{0,\ {\sigma}(1)\ \dots\ {\sigma}(n)}\ (z_{0};\ z_{\sigma (1)} ,\ \dots,\ z_{\sigma (n)} )\ R_{1\ \dots\ n}^{\sigma}\ (z_{1} ,\ \dots,\ z_n )\ ,
\label{eq:rpfaz} 
\end{eqnarray}
\begin{eqnarray}
R_{1\ \dots\ n}^{\sigma}\ (z_{1} ,\ \dots,\ z_n )\ F_{1\ \dots\ n,\ 0}\ (z_{1} ,\ \dots,\ z_n ;\ z_{0})\ = \nonumber \\
=\ F_{{\sigma}(1)\ \dots\ {\sigma}(n),\ 0}\ (z_{\sigma (1)} ,\ \dots,\ z_{\sigma (n)} ;\ z_{0} )\  R_{1\ \dots\ n}^{\sigma}\ (z_{1} ,\ \dots,\ z_n )\ ,
\label{eq:rpfbz} 
\end{eqnarray}
\begin{eqnarray}
F_{1\ \dots\ n,\ 0}\ (z_{1} ,\ \dots,\ z_n ;\ z_{0})\ R_{0,\ 1\ \dots\ n}\ (z_{0};\ z_{1} ,\ \dots,\ z_n )\ = \nonumber \nonumber \\
=\ F_{0,\ 1\ \dots\ n}\ (z_{0};\ z_{1} ,\ \dots,\ z_n )\ ,
\label{eq:rpfcz} 
\end{eqnarray}
\begin{eqnarray}
F_{0,\ 1\ \dots\ n}\ (z_{0};\ z_{1} ,\ \dots,\ z_n )\  R_{1\ \dots\ n,\ 0}\ (z_{1} ,\ \dots,\ z_n ;\ z_{0})\ = \nonumber \\ 
=\ F_{1\ \dots\ n,\ 0}\ (z_{1} ,\ \dots,\ z_n ;\ z_{0})\ .
\label{eq:rpfdz}
\end{eqnarray}
Conversely, suppose we have defined for any $n$ sets of matrices $F_{1,\ 2\ \dots\ n}\ (z_{1} ;\ z_{2} \dots,\ z_{n})$ and $F_{1\ \dots\ n-1,\ n}\ (z_{1},\ \dots,\ z_{n-1};\ z_{n})$ such that they coincide for $n\ =\ 2$ and such that they satisfy the above properties (\ref{eq:cocyclez} - \ref{eq:rpfdz}). Then it is posssible to define recursively a set of factorizing $F$-matrices as,
\begin{equation}
F_{1\ \dots\ n}\ (z_{1} ,\ \dots,\ z_n )\ =\ F_{2\ \dots\ n}\ (z_{2},\ \dots,\ z_{n})\ F_{1,\ 2\ \dots\ n}\ (z_{1} ;\ z_{2} \dots,\ z_{n})\ ,
\label{eq:deftf1}
\end{equation}
or equivalently, due to the cocycle relation (\ref{eq:cocyclez}),
\begin{equation}
F_{1\ \dots\ n}\ (z_{1} ,\ \dots,\ z_n )\ =\ F_{1\ \dots\ n-1}\ (z_{1},\ \dots,\ z_{n-1})\ F_{1\ \dots\ n-1,\ n}\ (z_{1},\ \dots,\ z_{n-1};\ z_{n})\ .
\label{eq:deftf2}
\end{equation}
\label{prop:factfz}
\end{prop}
\Proof
We use simplified notations, namely, without writing the explicit $z_i$'s dependence. First, the cocycle relation (\ref{eq:cocyclez}) follows from the definition as,
\begin{eqnarray*}
F_{2\ \dots\ n-1,\ n}\ F_{1,\ 2\ \dots\ n}\ &=&\ F_{2\ \dots\ n-1}^{-1}\ F_{2\ \dots\ n}\ F_{2\ \dots\ n}^{-1}\ F_{12\ \dots\ n}\ , \\
&=&\ F_{2\ \dots\ n-1}^{-1}\ F_{12\ \dots\ n}\ , \\
&=&\ F_{2\ \dots\ n-1}^{-1}\ F_{1\ \dots\ n-1}\ F_{1\ \dots\ n-1}^{-1}\ F_{12\ \dots\ n}\ , \\
&=&\ F_{1,\ 2\ \dots\ n-1}\ F_{1\ \dots\ n-1,\ n}\ .
\end{eqnarray*}
Equation (\ref{eq:rpfaz}) follows from the fundamental relation (\ref{eq:frz}) applied to the tensor product of $(n+1)$ modules $V(z_i)$. Let $q\ =\ (012\ \dots\ n)$, and $\sigma \in S_n$ acting only in the subset $q'\ =(12\ \dots\ n)$. Then we obtain,
\begin{eqnarray*}
F_{0{\sigma}(1)\ \dots\ {\sigma}(n)}\ R_{1\ \dots\ n}^{\sigma}\ =\ F_{01\ \dots\ n}\ ,
\end{eqnarray*}
hence leading to,
\begin{eqnarray*}
R_{q'}^{\sigma}\ F_{q'}^{-1}\ F_{q}\ &=&\ F_{{\sigma}(q')}^{-1}\ F_{q}\ =\ F_{{\sigma}(q')}^{-1}\ F_{0{\sigma}(q')}\ R_{q'}^{\sigma}\ ,  \\
&=&\ F_{0,\ {\sigma}(q')}\ R_{q'}^{\sigma}\ .
\end{eqnarray*}
Equation (\ref{eq:rpfbz}) follows from a similar argument.\\
To obtain (\ref{eq:rpfcz}), let us write again the fundamental relation (\ref{eq:frz}) for the element ${\sigma}_c \in S_{n+1}$ acting in $(01\ \dots\ n)$ as, ${\sigma}_c (X_{01\ \dots\ n})\ =\ X_{1\ \dots\ n0}$. We have,
\begin{eqnarray*}
F_{1\ \dots\ n0}\ R_{0,\ 1\ \dots\ n}\ =\ F_{01\ \dots\ n}\ .
\end{eqnarray*}
Then multiplying on the left by $F_{1\ \dots\ n}^{-1}$, we get the desired relation. Equation (\ref{eq:rpfdz}) is obtained in the same way.\\

Conversely, let eqs. (\ref{eq:cocyclez} - \ref{eq:rpfdz}) be satisfied, then by induction on $n$ one can prove in the same way as in Proposition (\ref{prop:tf}) that the total $F$-matrix $F_{12\ \dots\ n}$ can be well defined by the two equivalent  ordered products,
\begin{eqnarray*}
F_{12 \dots n}\ &=&\ F_{n-1\ n}\ F_{n-2,\ n-1\ n}\ \dots\ F_{1,\ 23 \dots n}\ , \\ 
&=&\ F_{12}\ F_{12,\ 3}\ \dots\ F_{12 \dots n-1,\ n}\ .
\end{eqnarray*}
To prove  the fundamental relation (\ref{eq:frz}) it is enough to prove it for the permutation of the two first elements ${\sigma}_{12}$ and for the cyclic permutation ${\sigma}_c$. This can be done as in the proof of Proposition (\ref{prop:fundf}). We have,
\begin{eqnarray*}
F_{213\ \dots\ n}\ R_{12}\ &=&\ F_{21}\ F_{21,\ 3}\ \dots\ F_{21\ \dots\ n-1,\ n}\ R_{12}\ , \\
&=&\ F_{21}\ F_{21,\ 3}\ \dots\ F_{21\ \dots\ n-2,\ n-1}\ R_{12}\ F_{12\ \dots\ n-1,\ n}\ , \\
&=&\ F_{21}\ R_{12}\ F_{12,\ 3}\ \dots\ F_{12\ \dots\ n-1,\ n}\ , \\
&=&\ F_{12}\ F_{12,\ 3}\ \dots\ F_{12\ \dots\ n-1,\ n}\ =\ F_{123\ \dots\ n}\ .
\end{eqnarray*}
For the cyclic permutation ${\sigma}_c$ we have similarly from eq. (\ref{eq:rpfdz}),
\begin{eqnarray*}
F_{1,\ 23\ \dots\ n}\ &=&\ F_{23\ \dots\ n,\ 1}\ R_{1,\ 23\ \dots\ n}\ , \\ 
&=&\ F_{23\ \dots\ n,\ 1}\ R_{12\ \dots\ n}^{{\sigma}_c}\ .
\end{eqnarray*}
Now multiplying by $F_{23\ \dots\ n}$ on the left we obtain,
\begin{eqnarray*}
F_{123\ \dots\ n}\ &=&\ F_{23\ \dots\ n}\ F_{23\ \dots\ n,\ 1}\ R_{12\ \dots\ n}^{{\sigma}_c}\ , \\
&=&\  F_{23\ \dots\ n1}\ R_{12\ \dots\ n}^{{\sigma}_c}\ .
\end{eqnarray*}
\QED

\begin{cor}
Let $F_{1\ \dots\ n}$ be factorizing $F$-matrices and $\tilde{F}_{1,\ 2\ \dots\ n}\ =\ F_{1\ \dots\ n}\ F_{2\ \dots\ n}^{-1}$ and similarly $\tilde{F}_{1\ \dots\ n-1,\ n}\ =\ F_{1\ \dots\ n}\ F_{1\ \dots\ n-1}^{-1}$. They satisfy,
\begin{eqnarray*}
\tilde{F}_{1,\ 2\ \dots\ n}\ \tilde{F}_{2\ \dots\ n-1,\ n}\ &=&\ \tilde{F}_{1\ \dots\ n-1,\ n}\ \tilde{F}_{1,\ 2\ \dots\ n-1}\ , \\
\tilde{F}_{0,\ 1\ \dots\ n}\ &=&\ \tilde{F}_{0,\ \sigma (1)\ \dots\ \sigma (n)}
\ ,\ \forall \sigma \in S_n \ , \\
\tilde{F}_{1\ \dots\ n,\ 0}\ &=&\ \tilde{F}_{\sigma (1)\ \dots\ \sigma (n),\ 0}
\ ,\ \forall \sigma \in S_n\ ,  \\
\tilde{F}_{1\ \dots\ n,\ 0}\ F_{1\ \dots\ n}\ R_{0,\ 1\ \dots\ n}\ &=&\ \tilde{F}_{0,\ 1\ \dots\ n}\ F_{1\ \dots\ n}\ .
\end{eqnarray*}
Conversely, suppose we have defined sets of matrices $\tilde{F}_{1,\ 2\ \dots\ n}$ and $\tilde{F}_{1\ \dots\ n-1,\ n}$ for any integer $n$ satisfying the above properties with $F_{1\ \dots\ n}\ =\ \tilde{F}_{1\ \dots\ n-1,\ n}\ \dots\ \tilde{F}_{12,\ 3}\ F_{12}$ or equivalently due to the cocycle relation, $F_{1\ \dots\ n}\ =\ \tilde{F}_{1,\ 2\ \dots\ n}\ \dots\ F_{n-1 n}$, then this set of matrices define factorizing $F$-matrices.
\label{cor:fundtf}
\end{cor}
\begin{rem}
The above proposition and its corollary show that the construction of factorizing $F$-matrices is reduced to the construction of the partial $F$-matrices $F_{12\ \dots\ n-1,\ n}$ and $F_{1,\ 2\ \dots\ n}$ for any integer $n$ provided they satisfy the properties (\ref{eq:cocyclez} - \ref{eq:rpfdz}). We will use this point of view in the next section.
\end{rem}
\begin{rem}
Suppose given a set of factorizing $F$-matrices $F_{1\ \dots\ n}$. Let $S_{1\ \dots\ n}$ be a set of matrices defined for any integer $n \geq 2$ such that for any element $\sigma \in S_n$, $S_{1\ \dots\ n}\ =\ S_{\sigma (1)\ \dots\ \sigma (n)}$. Then, $F_{1\ \dots\ n}^S\ =\ S_{1\ \dots\ n}\ F_{1\ \dots\ n}$ are also factorizing $F$-matrices. Note however that the relation between $F_{1,\ 2\ \dots\ n}^S$, $F_{1\ \dots\ n-1,\ n}^S$ and $F_{1,\ 2\ \dots\ n}$, $F_{1\ \dots\ n-1,\ n}$ are not simple. For example,
\begin{equation}
F_{1,\ 2\ \dots\ n}^S\ =\ F_{2\ \dots\ n}^{-1}\ S_{2\ \dots\ n}^{-1}\ S_{1\ \dots\ n}\ F_{1\ \dots\ n}\ ,
\end{equation}
while,
\begin{equation}
F_{1,\ 2\ \dots\ n}\ =\ F_{2\ \dots\ n}^{-1}\ F_{1\ \dots\ n}\ .
\end{equation}
The relations for the matrices  $\tilde{F}$ is simpler. We have,
\begin{equation}
\tilde{F}_{1,\ 2\ \dots\ n}^S\ =\ S_{1\ \dots\ n}\ \tilde{F}_{1,\ 2\ \dots\ n}\ S_{2\ \dots\ n}^{-1}\ .
\end{equation}
\end{rem}
\begin{rem}
Two sets of factorizing $F$-matrices $F_{1\ \dots\ n}$ and $G_{1\ \dots\ n}$ are related by a set of completly symmetric $S$-matrices, $G_{1\ \dots\ n}\ =\ S_{1\ \dots\ n}\ F_{1\ \dots\ n}$; indeed, it follows from the fundamental relation (\ref{eq:frz}) that,
\begin{equation}
R_q^{\sigma}\ =\ G_{{\sigma}(q)}^{-1}\ G_q\ = \ F_{{\sigma}(q)}^{-1}\ F_q\ ,
\end{equation}
leading to the relation,
\begin{equation}
G_q\ F_q^{-1}\ =\ G_{{\sigma}(q)}\ F_{{\sigma}(q)}^{-1}\ =\ S_q\ =\ S_{{\sigma}(q)}\ .
\end{equation}
\end{rem}

%% file: sect3
\section{Factorizing $F$-matrices associated to $Y(sl_2 )$ and ${\cal U}_q (\hat{sl_2})$}

The purpose of this section is to construct explicitly factorizing $F$-matrices, in the sense of Definition (\ref{defi:factfm}), acting in the irreducible $n^{th}$ tensor products of the spin-$1 \over 2$ (finite dimensional) irreducible representation of the quantum affine Hopf algebra  ${\cal U}_q (\hat{sl_2})$ and of its associated Yangian $Y(sl_2 )$. To achieve this we will make use of the Proposition (\ref{prop:factfz}) of the last section. These representations correspond to the integrable Heisenberg $XXX$ and $XXZ$ spin-$1 \over 2$ inhomogeneous quantum chains of lenght $n$.\\
We first recall the basic properties of the $R$-matrices associated to these finite dimensional irreducible fundamental representations in section 3.1. Then four sets of factorizing $F$-matrices are constructed in section 3.2, together with their relations.\\

\subsection{$R$-matrices associated to spin-$1 \over 2$  module of $Y(sl_2 )$ and ${\cal U}_q (\hat{sl_2})$}

The fundamental (two-dimensional) spin-$1 \over 2$ representation of ${\cal U}_q (\hat{sl_2})$ is associated to modules $V(z)$ such that $dim\  V(z)\ =\ 2$. The $R$-matrix acting in $V(\lambda ) \otimes V(\mu )$ is given as the following $4 \times 4$-matrix (up to a scalar factor to be discussed later on),
\begin{equation}
R(\lambda, \mu)=
\left(
\begin{array}{cccc}
1 & 0 & 0  & 0 \\
0 & b\ (\lambda, \mu) & c\ (\lambda, \mu) & 0 \\
0 & c\ (\lambda, \mu) & b\ (\lambda, \mu) & 0 \\
0 & 0 & 0 & 1
\end{array}
\right)\ .
\label{eq:vi}
\end{equation}
Corresponding quantum integrable models related to this class of $R$-matrices are the spin-$1 \over 2$ $XXX$ model for which $b$ and $c$ are rational functions of $\lambda$ and $\mu$,
\begin{equation}
b\ ( \lambda, \mu)\ =\frac{\lambda - \mu}{\lambda - \mu + \eta}\ \ \ 
c\ (\lambda, \mu)=\frac{\eta}{\lambda - \mu + \eta}\ ,
\label{eq:bcr}
\end{equation}
or the spin-$1 \over 2$ $XXZ$ model for which they are trigonometric functions of the difference $\lambda - \mu$,
\begin{equation}
b\ ( \lambda, \mu)=\frac { \sinh ( \lambda - \mu)}{\sinh (\lambda - \mu + \eta)}\ \ \ c\ (\lambda, \mu)=\frac { \sinh (  \eta)}{\sinh (\lambda - \mu + \eta)}\ ,
\label{eq:bct}
\end{equation}
where $\eta \in {\mathbf C}$ is the quantum deformation parameter. For the trigonometric case, it is useful to introduce other notations as, $e^{\lambda}\ =\ u$, $e^{\mu}\ =\ v$ with $z\ =\ {u \over v}$ and $e^{\eta}\ =\ q$. The $R$-matrix is depending on $\lambda, \mu$ only through $z$, leading to,
\begin{equation}
b\ (z)\ =\ \frac {q (z^2 - 1)}{z^2 q^2 - 1}\ \ \ c\ ( z)=\frac {z (q^2 - 1)}{z^2 q^2 - 1}\ .
\label{eq:bcq}
\end{equation}
As it is well known, the rational $R$-matrix can be obtained as a special limit of the trigonometric one.\\
Let us now list briefly the main properties of these $R$-matrices that we will use in the following \cite{CP1,FR1}.\\
The above $R$-matrices satisfy the Yang-Baxter equation,
\begin{eqnarray*}
R_{12}\ (\lambda_1,\lambda_2)\ R_{13}\ (\lambda_1,\lambda_3)\ R_{23}\ (\lambda_2,\lambda_3)\ =\ 
R_{23}\ (\lambda_2,\lambda_3)\ R_{13}\ (\lambda_1,\lambda_3)\ R_{12}\ (\lambda_1,\lambda_2)\ ,
\end{eqnarray*}
the unitary condition (provided it has no singularity and is invertible, namely, if $b\ ( \lambda_1, \lambda_2) \ne \pm c\ (\lambda_1, \lambda_2)$),
\begin{equation}
R_{12}\ (\lambda_1,\lambda_2)\ R_{21}\ (\lambda_2,\lambda_1)\ =\ {\mathbf 1}\ ,
\label{eq:eu}
\end{equation}
and the crossing symmetry relation,
\begin{equation}
({\gamma} \otimes {\mathbf 1})\ R_{12}\ (\lambda_1^s,\lambda_2)\ ({\gamma} \otimes {\mathbf 1})\ =\ R_{21}^{t_1}\ (\lambda_2 ,\lambda_1 )\ \rho (\lambda_1,\lambda_2)\ ,
\label{eq:ecross}
\end{equation}
with $\rho (\lambda_1,\lambda_2)$ being a scalar function, and $\gamma$ a $2 \times 2$ matrix such that $\gamma^2\ =\ {\mathbf 1}$, $\gamma^t \ =\ \pm \gamma$, the uppersript $t_j$ meaning the usual transposition of matrices in the corresponding space $(j)$.\\
For the rational case,  
\begin{eqnarray*}
\lambda_1^s\ =\ \lambda_1 - \eta\ ,\ \gamma\ =\ \sigma^y\ ,\ \rho (\lambda_1,\lambda_2)\ =\ \frac{\lambda_1 - \lambda_2 - \eta}{\lambda_1 - \lambda_2}\ ,
\end{eqnarray*}
and for the trigonometric case, 
\begin{eqnarray*}
\lambda_1^s\ =\ \lambda_1 - \eta + i \pi\ ,\ \gamma\ =\ \sigma^x\ ,\ \rho (\lambda_1,\lambda_2)\ =\ \frac{sh(\lambda_1 - \lambda_2 - \eta)}{sh(\lambda_1 - \lambda_2)}\ ,
\end{eqnarray*}
where $\sigma^x$ and $\sigma^y$ are the standard Pauli matrices.
Finally let us note also some useful symmetries of the $R$-matrices. First we have,
\begin{eqnarray}
R_{12}^{t_1 t_2}\ (\lambda_1,\lambda_2)\ &=&\ R_{12}\ (\lambda_1,\lambda_2)\ , \\
R_{21}\ (\lambda_1,\lambda_2)\ &=&\ R_{12}\ (\lambda_1,\lambda_2)\ , \\
(\sigma^x \otimes \sigma^x ) \ R_{12}\ (\lambda_1,\lambda_2)\ &=&\  R_{12}\  (\lambda_1,\lambda_2)\ (\sigma^x \otimes \sigma^x )\ . 
\label{eq:tinv}
\end{eqnarray}
For the rational case this last equation is true for any Pauli matrix. Moreover, denoting by $e_{\alpha}^{(ii)}$ the elementary $2 \times 2$-matrices 
such that ${(e^{(ij)})}_{kl}={\delta}_{ik}{\delta}_{jl}$ acting non-trivially only in the $\alpha$-th space,
\begin{eqnarray}
R_{12}\ (\lambda_1,\lambda_2)\ e_1^{(ii)}\ e_2^{(ii)}\ &=&\ e_1^{(ii)} \ e_2^{(ii)}\ , \\
e_1^{(ii)} \ e_2^{(ii)}\ R_{12}\ (\lambda_1,\lambda_2)\ &=&\ e_1^{(ii)} \ e_2^{(ii)}\ .
\label{eq:ree}
\end{eqnarray}
Let us now translate this list of properties for the elementary $R$-matrices to ordered products of them (corresponding to the co-multiplication action).
From the Yang-Baxter equation for the $R$-matrix we get,
\begin{prop}
We define as in Proposition (\ref{prop:rsigmaz}) the $R$-matrices acting in $V(\lambda_1)\ \otimes \dots \otimes V(\lambda_n)$. 
\begin{eqnarray}
R_{1,\ 23\ \dots\ n}\ &=&\  R_{1,\ 23\ \dots\ n}\ ( z_{1} ;\ z_{2} \dots,\ z_{n})\ , \nonumber\\
&=&\ R_{1n}\ (z_1 ,\ z_n )\ \dots\ R_{12}\ ( z_1 ,\ z_2 )\ , \\  
R_{12\ \dots\ n-1,\ n}\ &=&\ R_{12\ \dots\ n-1,\ n}\ ( z_{1},\ \dots,\ z_{n-1};\ z_{n} )\ , \nonumber \\
&=&\ R_{1\ n}\ (z_1 ,\ z_n )\ \dots\ R_{n-1\ n}\ (z_{n-1} ,\ z_n )\ .
\end{eqnarray}
They satisfy the following cocycle relation,
\begin{equation}
R_{2\ \dots\ n-1,\ n}\ R_{1,\ 2\ \dots\ n}\ =\ R_{1,\ 2\ \dots\ n-1}\ R_{1\ \dots\  n-1,\ n}\ ,
\label{eq:cocr}
\end{equation}
and the unitarity relation,
\begin{equation}
R_{1,\ 23\ \dots\ n}\ R_{23\ \dots\ n,\ 1}\ =\ {\mathbf 1}\ .
\label{eq:urn}
\end{equation}
provided $b\ ( \lambda_1, \lambda_j) \ne \pm c\ (\lambda_1, \lambda_j)$, $j\ =\ 2,\ \dots\ n$.
\label{prop:cocyr}
\end{prop}
\Proof
Let us first prove the unitarity relation. We have, from the elementary relation $R_{12}\ R_{21}\ =\ {\mathbf 1}$,
\begin{eqnarray*}
R_{1,\ 23\ \dots\ n}\ R_{23\ \dots\ n,\ 1}\ &=&\ R_{1,\ 3\ \dots\ n}\ R_{12}\ R_{21}\ R_{3\ \dots\ n,\ 1}\ , \\
&=&\ R_{1,\ 3\ \dots\ n}\ R_{3\ \dots\ n,\ 1}\ ,
\end{eqnarray*}
leading to the relation by induction on $n$. Note that the condition $b\ ( \lambda_1, \lambda_j) \ne \pm c\ (\lambda_1, \lambda_j)$, $j\ =\ 2,\ \dots\ n$ is necessary to ensure invertibility of the $R$-matrix.\\
The cocycle relation is also proven by induction on $n$. First it reduces to the Yang-Baxter relation for $n\ =\ 3$. Suppose it is verified up to some integer $n$. Then we have, for the set of $n$ modules labelled by $(13\ \dots\ n+1)$, the cocycle relation is verified, namely,
\begin{eqnarray*}
R_{3\ \dots\ n,\ n+1}\ R_{1,\ 3\ \dots\ n+1}\ =\ R_{1,\ 3\ \dots\ n}\ R_{13\ \dots\ n,\ n+1}\ .
\end{eqnarray*}
Let us now consider $(n+1)$ modules labelled by the set $(123\ \dots\ n+1)$, we have,
\begin{eqnarray*}
R_{2\ \dots\ n,\ n+1}\ R_{1,\ 2\ \dots\ n+1}\ &=&\ R_{2\ n+1}\ R_{3\ \dots\ n,\ n+1}\ R_{1,\ 3\ \dots\ n+1}\ R_{12}\ , \\
&=&\ R_{2\ n+1}\ R_{1,\ 3\ \dots\ n}\ R_{13\ \dots\ n,\ n+1}\ R_{12}\ , \\
&=&\ R_{1,\ 3\ \dots\ n}\ R_{2\ n+1}\ R_{1\ n+1}\ R_{12}\ R_{3\ \dots\ n,\ n+1}\ , \\
&=&\ R_{1,\ 3\ \dots\ n}\ R_{12}\ R_{1\ n+1}\ R_{2\ n+1}\ R_{3\ \dots\ n,\ n+1}\ , \\
&=&\ R_{1,\ 2\ \dots\ n}\ R_{123\ \dots\ n+1}\ .
\end{eqnarray*}
\QED
\begin{rem}
For $R$-matrices which are representations of a universal $R$-operator associated to a quasi-triangular Hopf algebra, this cocycle relation follows from the corresponding universal cocycle relation,
\begin{equation}
R_{12}\ ({\Delta} \otimes id)\ R\ =\ R_{23}\ (id \otimes {\Delta})\ R\ .
\label{eq:ucocyr}
\end{equation}
\end{rem}
The crossing symmetry of the elementary $R$-matrix leads to the following Proposition.
\begin{prop}
\begin{equation}
\gamma_0 \ R_{0,\ 1\ \dots\ n}\ ({\lambda}_0^s\ ;\ \lambda_1 ,\ \dots\ ,\ \lambda_n)\ \gamma_0\ =\ \rho ({\lambda}_0\ ;\ \lambda_1 ,\ \dots\ ,\ \lambda_n)\ R_{1\ \dots\ n,\ 0}^{t_0}\ ({\lambda}_0\ ;\ \lambda_1 ,\ \dots\ ,\ \lambda_n)\ ,
\label{eq:crn}
\end{equation}
or equivalently,
\begin{equation}
R_{0,\ 1\ \dots\ n}\ ({\lambda}_0\ ;\ \lambda_1 ,\ \dots\ ,\ \lambda_n)\ \gamma_0\ R_{0,\ 1\ \dots\ n}^{t_0}\ ({\lambda}_0^s\ ;\ \lambda_1 ,\ \dots\ ,\ \lambda_n)\ \gamma_0\ =\ \rho ({\lambda}_0\ ;\ \lambda_1 ,\ \dots\ ,\ \lambda_n)\ {\mathbf 1}\ ,
\label{eq:qdet}
\end{equation}
where $\rho ({\lambda}_0\ ;\ \lambda_1 ,\ \dots\ ,\ \lambda_n)\ =\ \prod_{i = 1}^{n}\ \rho ({\lambda}_0 ,\ \lambda_i)$, and $\gamma_0$ is the matrix $\gamma$ of eq. (\ref{eq:ecross}) acting in space $0$.
\label{prop:crossrn}
\end{prop}
\Proof
We have, $R_{0,\ 1\ \dots\ n}\ =\ R_{0n}\ \dots\ R_{01}$, hence, using simplified notations,
\begin{eqnarray*}
\gamma_0 \ R_{0,\ 1\ \dots\ n}\ ({\lambda}_0^s\ ;\ \lambda_1 ,\ \dots\ ,\ \lambda_n)\ \gamma_0\ &=&\ \gamma_0 \ R_{0n}\ ({\lambda}_0^s,\ \lambda_n)\ \gamma_0\ \dots\ \gamma_0 \ R_{01}\ ({\lambda}_0^s,\ \lambda_1)\ \gamma_0\ , \\
&=&\ \rho ({\lambda}_0,\ \lambda_n)\ R_{n0}^{t_0} (\lambda_n ,\ \lambda_0)\ \dots\ \rho ({\lambda}_0,\ \lambda_1)\ R_{10}^{t_0} (\lambda_1 ,\ \lambda_0)\ , \\
&=&\ \rho ({\lambda}_0\ ;\ \lambda_1 ,\ \dots\ ,\ \lambda_n)\ (R_{10}\ \dots\ R_{n0})^{t_0}\ , \\
&=&\ \rho ({\lambda}_0\ ;\ \lambda_1 ,\ \dots\ ,\ \lambda_n)\ R_{1\ \dots\ n,\ 0}^{t_0}\ .
\end{eqnarray*}
\QED
\begin{prop}
Let us define the operation $\dagger$ on elements $X_{1\ \dots\ n}\ (\lambda_1 ,\ \dots\ ,\ \lambda_n)$ belonging to the space $End\ (V(\lambda_1)\ \otimes \dots \otimes V(\lambda_n))$ as,
\begin{equation}
X_{1\ \dots\ n}^{\dagger}\ (\lambda_1 ,\ \dots\ ,\ \lambda_n)\ =\ X_{1\ \dots\ n}^{t_1\ \dots\ t_n }\ (- \lambda_1 ,\ \dots\ ,\ - \lambda_n)\ .
\label{eq:dagger}
\end{equation}
We have,
\begin{equation}
R_{0,\ 1\ \dots\ n}^{\dagger}\ =\ R_{0,\ 1\ \dots\ n}^{-1}\ =\ R_{1\ \dots\ n,\ 0}\ .
\label{eq:dagr}
\end{equation}
\label{prop:daggerr}
\end{prop}
\Proof
First we remark that the relation is true for the elementary $R$-matrix, namely, $R_{0i}^{\dagger}\ =\ R_{0i}$. Then,
\begin{eqnarray*}
R_{0,\ 1\ \dots\ n}^{\dagger}\ &=& (R_{0n}\ \dots\ R_{01})^{\dagger}\ , \\
&=&\ R_{01}^{\dagger}\ \dots\ R_{0n}^{\dagger}\ , \\
&=&\ R_{10}\ \dots\ R_{n0}\ , \\
&=&\ R_{1\ \dots\ n,\ 0}\ =\ R_{0,\ 1\ \dots\ n}^{-1}\ .
\end{eqnarray*}
\QED
\begin{defi}
We now define the $2^{n+1} \times 2^{n+1}$ matrix $R_{0,\ 1\ \dots\ n}$ as a $2 \times 2$ matrix in module $V_0 (\lambda_0)$ whose entries are $2^n \times 2^n$ matrices acting in $V_1 (\lambda_1)\ \otimes \dots \otimes V_n (\lambda_n)$ as,
\begin{equation}
R_{0,\ 1\ \dots\ n}\ =\  \left(
\begin{array}{cc}
A_{1\ \dots\ n}\ ({\lambda}_0\ ;\ \lambda_1 ,\ \dots\ ,\ \lambda_n) &  B_{1\ \dots\ n}\ ({\lambda}_0\ ;\ \lambda_1 ,\ \dots\ ,\ \lambda_n)\\
C_{1\ \dots\ n}\ ({\lambda}_0\ ;\ \lambda_1 ,\ \dots\ ,\ \lambda_n)  & D_{1\ \dots\ n}\ ({\lambda}_0\ ;\ \lambda_1 ,\ \dots\ ,\ \lambda_n)
\end{array}
\right)_{[0]}\ .
\label{eq:RABCD}
\end{equation}
Whenever it is not ambiguous we will write,\\ 
$A_{1\ \dots\ n}\ ({\lambda}_0\ ;\ \lambda_1 ,\ \dots\ ,\ \lambda_n)\ =\ A_{1\ \dots\ n}\ ({\lambda}_0)\ =\ A\ ({\lambda}_0)$ ,\\ 
and so on for $B, C, D$. With these notations we have,
\begin{equation}
R_{0,\ 1\ \dots\ n}^{\dagger}\ =\  \left(
\begin{array}{cc}
A_{1\ \dots\ n}^{\dagger}\ ({\lambda}_0) &  C_{1\ \dots\ n}^{\dagger}\ ({\lambda}_0)\\
B_{1\ \dots\ n}^{\dagger}\ ({\lambda}_0)  & D_{1\ \dots\ n}^{\dagger}\ ({\lambda}_0)
\end{array}
\right)_{[0]}\ .
\label{eq:Rdag}
\end{equation}
\label{defi:RRdag}
\end{defi}
For later use we now give triangularity properties of the $2^n \times 2^n$ matrices $A,\ B,\ C,\ D$.
\begin{prop}
The $2^n \times 2^n$ matrices $A_{1\ \dots\ n}$ and $C_{1\ \dots\ n}$ are upper triangular matrices while $B_{1\ \dots\ n}$ and $D_{1\ \dots\ n}$ are lower triangular. Moreover, $B_{1\ \dots\ n}$ and $C_{1\ \dots\ n}$ have only zeros on their diagonal. We also have for the diagonal part of $A_{1\ \dots\ n}$ and $D_{1\ \dots\ n}$,
\begin{eqnarray*}
diag\ (A_{1\ \dots\ n}\ (\lambda_0 ))\ =\ \otimes_{i\ =\ 1}^{n}  \left( 
\begin{array}{cc}
1 & 0 \\
0 & b  (\lambda_0,\lambda_{i})
\end{array}
\right)_{[i]}\ ,
\end{eqnarray*}
and,
\begin{eqnarray*}
diag\ (D_{1\ \dots\ n}\ (\lambda_0 ))\ =\ \otimes_{i\ =\ 1}^{n} \left( 
\begin{array}{cc}
b  (\lambda_0,\lambda_{i})& 0 \\
0 & 1
\end{array}
\right)_{[i]}\ .
\end{eqnarray*}
\label{prop:triabcd}
\end{prop}
\Proof
By induction on $n$. The proposition is true for $n\ =\ 1$ from the definition of the $R$-matrix. Let it be true for $n-1$. Then using, 
\begin{eqnarray*}
R_{0,\ 1\ \dots\ n}\ =\ R_{0n}\ R_{0,\ 1\ \dots\ n-1}\ ,
\end{eqnarray*}
and writing this relation as a $2 \times 2$ matrix identity in space $0$,
\begin{eqnarray*}
B_{1\ \dots\ n}\ (\lambda_0)\ =\ B_{1\ \dots\ n-1}\ (\lambda_0)
\otimes
\left( 
\begin{array}{cc}
1 & 0 \\
0 & b  (\lambda_0,\lambda_{n})
\end{array}
\right)_{[n]}
+
D_{1\ \dots\ n-1}\ (\lambda_0)
\otimes
\left( 
\begin{array}{cc}
0 & 0 \\
c\ (\lambda_0,\lambda_{n}) & 0
\end{array}
\right)_{[n]}\ ,
\end{eqnarray*}
where the second terms in the above tensor product is a matrix in module $V_n (\lambda_n)$. Hence, if $B_{1\ \dots\ n-1}\ (\lambda_0)$ and $D_{1\ \dots\ n-1}\ (\lambda_0)$ are lower triangular matrices, so is $B_{1\ \dots\ n}\ (\lambda_0)$. Moreover, if $B_{1\ \dots\ n-1}\ (\lambda_0)$ has only zeros on the diagonal it is also the case for $B_{1\ \dots\ n}\ (\lambda_0)$. Similarly,
\begin{eqnarray*}
D_{1\ \dots\ n}\ (\lambda_0)\ =\ B_{1\ \dots\ n-1}\ (\lambda_0)
\otimes
\left( 
\begin{array}{cc}
0 & c\ (\lambda_0,\lambda_{n})\\
0 & 0
\end{array}
\right)_{[n]}
+
D_{1\ \dots\ n-1}\ (\lambda_0)
\otimes
\left( 
\begin{array}{cc}
b  (\lambda_0,\lambda_{n})& 0 \\
0 & 1
\end{array}
\right)_{[n]}\ .
\end{eqnarray*}
Then, since $B_{1\ \dots\ n-1}\ (\lambda_0)$ is strictly lower triangular, the first term in this sum is also lower triangular. $D_{1\ \dots\ n-1}\ (\lambda_0)$ being lower triangular, the second term is also lower triangular, leading to the desired result. Then the diagonal part of $D_{1\ \dots\ n}\ (\lambda_0)$ can be obtained from the above relation since we have,
\begin{eqnarray*}
diag\ (D_{1\ \dots\ n}\ (\lambda_0))\ =\ diag\ (D_{1\ \dots\ n-1}\ (\lambda_0))
\otimes
\left( 
\begin{array}{cc}
b  (\lambda_0,\lambda_{n})& 0 \\
0 & 1
\end{array}
\right)_{[n]}\ ,
\end{eqnarray*}
leading to the result by induction on $n$. The result for the matrices $A_{1\ \dots\ n}$ and $C_{1\ \dots\ n}$ are obtained in a similar way.
\QED
Finally, let us give the symmetry relations induced by the elementary relation (\ref{eq:tinv}).
\begin{prop}
Let ${\cal C}\ =\ {\cal C}_{1\ \dots\ n}\ =\ {\sigma}_1^x\ \dots\ {\sigma}_n^x$. Then ${\cal C}^2\ =\ {\mathbf 1}$ and $\left[ {\cal C}_{01\ \dots\ n},\ R_{0,\ 1\ \dots\ n} \right]\ =\ 0$. It follows that,
\begin{eqnarray*}
{\cal C}\ A\ (\lambda_0 )\ {\cal C}\ &=&\ D\ (\lambda_0 )\ , \\
{\cal C}\ B\ (\lambda_0 )\ {\cal C}\ &=&\ C\ (\lambda_0 )\ .
\end{eqnarray*}
\label{prop:calc}
\end{prop}
\begin{rem}
${\cal C}_{1\ \dots\ n}$ relates upper triangular matrices in $V_1 (\lambda_1 ) \otimes\ \dots\ \otimes V_n (\lambda_n )$ to lower triangular ones.
\end{rem}

\subsection{Factorizing $F$-matrices associated to irreducible tensor products of the spin-$1 \over 2$ module of $Y(sl_2 )$ and ${\cal U}_q (\hat{sl_2})$}

Section 3.1 provides us with all the necessary ingredients to give explicit factorizing $F$-matrices associated to the $R$-matrices (\ref{eq:vi}).\\

The first very simple fact we can show is the following,
\begin{lem}
Let $b\ ( \lambda_1 ,\ \lambda_2) \ne 0$, and also $b\ ( \lambda_1, \lambda_2) \ne \pm c\ (\lambda_1, \lambda_2)$, then the $R$-matrix (\ref{eq:vi}) is unitary if and only if there exists an invertible matrix $F_{12}\ (\lambda_1 ,\ \lambda_2)$ such that,
\begin{equation}
F_{21}\ (\lambda_2 ,\ \lambda_1)\ R_{12}\ (\lambda_1 ,\ \lambda_2)\ =\ F_{12}\ (\lambda_1 ,\ \lambda_2)\ .
\label{eq:frfz}
\end{equation}
Moreover a particular choice of $F_{12}\ (\lambda_1 ,\ \lambda_2)$ is a lower triangular matrix in $V_1 (\lambda_1) \otimes V_2 (\lambda_2)$, given as,
\begin{equation}
F(\lambda_1,\ \lambda_2)\ =\ 
\left(
\begin{array}{cccc}
1 & 0 & 0  & 0 \\
0 & 1 & 0 & 0 \\
0 & c\ (\lambda_1,\ \lambda_2) & b\ ( \lambda_1,\ \lambda_2) & 0 \\
0 & 0 & 0 & 1
\end{array}
\right)_{[12]}\ ,
\end{equation}
or equivalently,
\begin{equation}
F_{12}=e_1^{(11)} + e_1^{(22)}R_{12}=e_2^{(22)} + e_2^{(11)}R_{12}\ .
\label{eq:f12}
\end{equation}
\end{lem}
\Proof
If an invertible $F_{12}$ exists, then $R_{12}$ is unitary,
\begin{eqnarray*}
R_{21}\ R_{12}\ =\ F_{12}^{-1}\ F_{21}\ F_{21}^{-1}\ F_{12}\ =\ {\mathbf 1}\ .
\end{eqnarray*}
Conversely, if $R_{12}$ is unitary with $b\ ( \lambda_1, \lambda_2) \ne \pm c\ (\lambda_1, \lambda_2)$ and $b\ ( \lambda_1 ,\ \lambda_2) \ne 0$, then we have the relations,
\begin{eqnarray}
b\ ( \lambda_1,\ \lambda_2)\ b\ (  \lambda_2,\ \lambda_1)\ +\ c\ (\lambda_1,\ \lambda_2)\ c\ ( \lambda_2,\ \lambda_1)\ &=&\ 1\ ,  \nonumber\\
b\ ( \lambda_1,\ \lambda_2)\ c\ (\lambda_2,\ \lambda_1)\ +\ c\ (\lambda_1,\ \lambda_2)\ b\ (  \lambda_2,\ \lambda_1)\ &=&\ 0 \ .
\label{eq:bccb0}
\end{eqnarray}
Using these relations, we directly obtain the equation (\ref{eq:frfz}). Moreover the determinant of the matrix  $F_{12}\ (\lambda_1 ,\ \lambda_2)$ is equal to 
$b\ ( \lambda_1 ,\ \lambda_2)$. Hence it is invertible if and only if $b\ ( \lambda_1 ,\ \lambda_2) \ne 0$.
\QED
\begin{rem}
If $R$ is some unitary $R$-matrix which admits a decomposition by factorizing $F$-matrices, then,
\begin{eqnarray*}
P_{12}\ R_{12}\ (\lambda_1,\ \lambda_2)\ =\ F_{12}^{-1}\ (\lambda_2,\ \lambda_1)\ P_{12}\ F_{12}\ (\lambda_1,\ \lambda_2)\ .
\end{eqnarray*}
It means that if the point $\lambda_1\ =\ \lambda_2$ is regular, the matrix $P_{12}\ R_{12}$ has the same spectrum as $P_{12}$. Note however that at different spectral parameters, the above factorization {\em is not a similarity transformation} for $P_{12}\ R_{12}$. Hence, the spectrum of $P_{12}\ R_{12}$ is not restricted in general to the one of the permutation operator $P_{12}$. This is due to the non-trivial spectral parameter dependence.\\
In our example (\ref{eq:vi}), note that at the point $\lambda_1\ =\ \lambda_2$, $R_{12}\ =\ P_{12}$ which cannot be decomposed by an invertible $F_{12}$.
\end{rem}
\begin{rem}
$F_{12}\ (\lambda_1 ,\ \lambda_2)\ =\ F_{12}\ (\lambda_1 - \lambda_2)$ is not unique since we can obtain another solution to the above factorization problem by multiplying it on the left by a completly symmetric $4 \times 4$ matrix $S_{12}\ (\lambda_1 ,\ \lambda_2)\ =\ S_{21}\ (\lambda_2 ,\ \lambda_1)$. Moreover two different solutions to this factorisation problem for the $R$-matrix (\ref{eq:vi}) are necessarily related by such an $S_{12}$ matrix.
\end{rem}
\begin{rem}
Since $\sigma^x \otimes \sigma^x$ commutes with the $R$-matrix, the matrix $\overline{F_{12}}\ =\ \sigma^x \otimes \sigma^x\ F_{12}\ \sigma^x \otimes \sigma^x$ also factorizes the $R$-matrix (\ref{eq:vi}). It is given as an upper triangular matrix in $V_1 (\lambda_1 ) \otimes V_2 (\lambda_2 )$,
\begin{equation}
\overline{F_{12}}\ =\ e_1^{(22)}\ +\ e_1^{(11)}\ R_{12}\ =\ e_2^{(11)}\ +\ e_2^{(22)}\ R_{12}\ .
\label{eq:fbar}
\end{equation}
\end{rem}
Similarly, using the fact that $R_{12}^{t_1 t_2}\ =\ R_{12}$, we obtain,
\begin{lem}
\begin{equation}
G_{12}\ (\lambda_1 ,\ \lambda_2)\ =\ ( F_{21}^{t_1 t_2}\ (\lambda_2 ,\ \lambda_1))^{-1}\ ,
\label{eq:g}
\end{equation}
is also factorizing the $R$-matrix (\ref{eq:vi}). It satisfies, $F_{12}\ G_{12}^{-1}\ =\ \Delta_{12}^{-1}$ with $\Delta_{12}\ =\ \Delta_{21}$ being the following diagonal matrix,
\begin{equation}
\Delta(\lambda_1 ,\ \lambda_2)=
\left(
\begin{array}{cccc}
1 &  &   &  \\
& b^{-1}(\lambda_2 ,\ \lambda_1) &  &  \\
 &  &b^{-1}(\lambda_1,\ \lambda_2) &  \\
 &  &  & 1
\end{array}
\right)_{[12]}\ .
\label{eq:Delta}
\end{equation}
\end{lem}
\Proof
By direct computation using the unitarity of the $R$-matrix written as in equations  (\ref{eq:bccb0}).
\QED
\begin{rem}
The formula (\ref{eq:f12}) for $F_{12}$ involves the $R$-matrix and for example the projectors $e_1^{(11)}$ and $e_1^{(22)}$ acting in the vector space $V_1$ only. Hence one can expect the representation of the  partial $F$-matrix $F_{1,\ 2\ \dots\ n}$ (which corresponds to the comultiplication action in second space which here acts only on $R$) to be given as,
\begin{eqnarray*}
F_{1,\ 2\ \dots\ n}\ =\ e_1^{(11)}\ +\ e_1^{(22)}\ R_{1,\ 2\ \dots\ n}\ . 
\end{eqnarray*}
Similarly, one can expect also from the second part of the formula (\ref{eq:f12}), 
\begin{eqnarray*}
F_{2\ \dots\ n,\ 1}\ =\ e_1^{(22)}\ +\ e_1^{(11)}\ R_{2\ \dots\ n,\ 1}\ . 
\end{eqnarray*}
\end{rem}
This remark leads to the construction of the corresponding factorizing $F$-matrices acting in the tensor product $V_1 (\lambda_1) \otimes\ \dots\ \otimes V_n (\lambda_n)$ as,
\begin{prop}
Let the partial $F$-matrices be given as ,
\begin{eqnarray}
F_{1,\ 2\ \dots\ n}\ =\ e_1^{(11)}\ +\ e_1^{(22)}\ R_{1,\ 2\ \dots\ n}\ ,  \label{eq:frn} \\
F_{2\ \dots\ n,\ 1}\ =\ e_1^{(22)}\ +\ e_1^{(11)}\ R_{2\ \dots\ n,\ 1}\ .
\label{eq:fln}
\end{eqnarray}
They define a set of factorizing $F$-matrices in the sense of Definition (\ref{defi:factfm}) for the $R$-matrix (\ref{eq:vi}) provided the $\lambda_i$'s  in $V_1 (\lambda_1) \otimes\ \dots\ \otimes V_n (\lambda_n)$ are such that, 
$b\ ( \lambda_i,\ \lambda_j)$ and $c\ (\lambda_i,\ \lambda_j)$ are finite, $b\ ( \lambda_i,\ \lambda_j) \ne 0\ ,\ \forall i\ \ne\ j $ and $b\ ( \lambda_i,\ \lambda_j) \ne \pm c\ (\lambda_i,\ \lambda_j)\ ,\ \forall i\ \ne\ j $.
Moreover in this case, $F_{1\ \dots\ n}\ (\lambda_1,\ \dots\ \lambda_n)$ is a lower triangular matrix with determinant equal to $\prod_{i < j}\ b\ ( \lambda_i,\ \lambda_j)$ and is left invariant by the uniform shift $\lambda_i\ \rightarrow \lambda_i\ +\ \delta$ for any $\delta \in {\mathbf C}$.\\
Similarly, $F_{n\ \dots\ 1}$ is an upper triangular matrix.
\label{prop:exf}
\end{prop}
\Proof
We apply Proposition (\ref{prop:factfz}). First we have that the two sets of partial $F$-matrices coincide for $F_{12}$ by construction.\\
Then, let us prove the cocycle relation. We have,
\begin{eqnarray*}
F_{2\ \dots\  n-1,\ n }\ F_{1,\ 2\ \dots\ n} & =&(e_{n}^{(22)}\ +\ e_{n}^{(11)}\  
R_{2\ \dots\ n-1,\ n})\  
(e_{1}^{(11)} + e_{1}^{(22)} R_{1,\ 2\ \dots\ n})\ = \\ 
=\ e_{1}^{(11)}e_{n}^{(22)} &+& 
e_{1}^{(22)} \ e_{n}^{(22)}\ R_{1,\ 2\ \dots\ n}\ +\ 
e_{1}^{(11)} \ e_{n}^{(11)}\ R_{2\ \dots\ n-1,\ n}\ +\\
&+&\ e_{1}^{(22)}\ e_{n}^{(11)}\ R_{2\ \dots\ n-1,\ n}\ R_{1,\ 2\ \dots\ n}\ .  
\end{eqnarray*}
Using (\ref{eq:ree}) in the second and third terms, and Propositin (\ref{prop:cocyr}) in the last one this is rewritten as
\begin{eqnarray*}
e_{1}^{(11)}e_{n}^{(22)}\ &+&\
e_{1}^{(22)}\ e_{n}^{(22)}\ R_{1,2 \ldots n-1}\ +\ 
e_{1}^{(11)}\ e_{n}^{(11)}\ R_{1\ldots n-1,n}\ +\\
+\ e_{1}^{(22)}\ e_{n}^{(11)}\ R_{1,2 \ldots n-1}\ R_{1\ldots n-1,n}\ &=&\ 
(e_{1}^{(11)}\ +\ e_{1}^{(22)}\ R_{1,2 \ldots n-1})\ 
(e_{n}^{(22)}\ +\ e_{n}^{(11)}\   
R_{1\ldots n-1,n})\ , \\ 
&=&\ F_{1, 2 \ldots n-1}\  F_{1\ldots n-1, n}\ ,
\end{eqnarray*}
completing the proof of the cocycle relation for $F$.\\
Now, the symmetry relations (\ref{eq:rpfaz}, \ref{eq:rpfbz}) are direct consequences of the symmetry relations (\ref{eq:r0qsa}, \ref{eq:r0qsb}) for the $R$ matrix itself. For example, for any $\sigma \in S_n$,
\begin{eqnarray*}
R_{1\ \dots\ n}^{\sigma}\ F_{0,\ 1\ \dots\ n}\ &=&\ R_{1\ \dots\ n}^{\sigma}\ (e_0^{(11)}\ +\ e_0^{(22)}\ R_{0,\ 1\ \dots\ n})\ , \\
&=&\ (e_0^{(11)}\ +\ e_0^{(22)}\ R_{0,\ \sigma(1)\ \dots\ \sigma(n)})\ R_{1\ \dots\ n}^{\sigma}\ , \\
&=&\ F_{0,\ \sigma(1)\ \dots\ \sigma(n)}\ R_{1\ \dots\ n}^{\sigma}\ .
\end{eqnarray*}
Finally, we have to prove that these $F$-matrices indeed factorize the $R$-matrices. It is sufficient to prove eq. (\ref{eq:rpfcz}). Using unitarity of the $R$-matrix, we have,
\begin{eqnarray*}
F_{2\ \dots\ n,\ 1}\ R_{1,\ 2\ \dots\ n}\ &=&\  (e_1^{(22)}\ +\ e_1^{(11)}\ R_{2\ \dots\ n,\ 1})\ R_{1,\ 2\ \dots\ n}\ ,\\
&=&\ e_1^{(22)}\ R_{1,\ 2\ \dots\ n}\ +\ e_1^{(11)}\ , \\
&=&\ F_{1,\ 2\ \dots\ n}\ .
\end{eqnarray*}
What remains to be proven is the invertibility of $F$ and its lower triangularity. Writting $F_{0,\ 1\ \dots\ n}$ as a $2 \times 2$ matrix in space zero, we obtain,
\begin{equation}
F_{0,\ 1\ \dots\ n}\ (\lambda_0;\ \lambda_1,\ \dots\ ,\ \lambda_n)\ =\  \left(
\begin{array}{cc}
{\mathbf 1} & {\mathbf 0} \\
C_{1\ \dots\ n} (\lambda_0\ ;\ \lambda_1,\ \dots\ ,\ \lambda_n)  & D_{1\ \dots\ n}\ (\lambda_0\ ;\ \lambda_1,\ \dots\ ,\ \lambda_n)
\end{array}
\right)_{[0]}\ .\\
\label{eq:22f}
\end{equation}
Hence the determinant of $F_{0,\ 1\ \dots\ n}$ is equal to the determinant of $D_{1\ \dots\ n}\ (\lambda_1,\ \dots\ ,\ \lambda_n)$, namely is equal to,
\begin{equation}
det\ F_{0,\ 1\ \dots\ n}\ =\ \prod_{i\ =\ 1}^n\ \ b\ ( \lambda_0,\ \lambda_i)\ , \label{eq:detf}
\end{equation}
which leads to the desired result for the total $F$-matrix $F_{1\ \dots\ n}$,
\begin{equation}
det\ F_{1\ \dots\ n}\ =\ \prod_{i < j}\ b\ ( \lambda_i,\ \lambda_j)\ .
\label{eq:detft}
\end{equation}
It is non zero if and only if all $b\ ( \lambda_i,\ \lambda_j)$ are non zero, namely, if and only if all $\lambda_i$'s are different.\\
Being a product of lower triangular matrices (which follows from eq. (\ref{eq:22f}) and Proposition (\ref{prop:triabcd})), $F_{12\ \dots\ n}$ is also lower triangular. The fact that $F_{n\ \dots\ 21}$ is upper triangular follows from the same arguments.\\
By construction the total $F$-matrix is depending on the $\lambda_i$'s through the $R$-matrices $R_{j,\ j+1\ \dots\ n}$. These $R$-matrices are depending on the $\lambda_k$ for $k \geq j$'s only through the differences $\lambda_j\ -\ \lambda_k$. This means that the total $F$-matrix is depending on the $\lambda_i$'s only through the differences $\lambda_i\ -\ \lambda_j$ for $i < j$, leading to the invariance of $F_{1\ \dots\ n}\ (\lambda_1\ ,\ \dots\ ,\ \lambda_n)$ under a uniform arbirary shift $\lambda_i \rightarrow \lambda_i\ +\ \delta$ for any $\delta \in {\mathbf C}$.
\QED
\begin{prop}
Let $\overline{F}_{1\ \dots\ n}\ =\ {\cal C}_{1\ \dots\ n}\ F_{1\ \dots\ n}\ {\cal C}_{1\ \dots\ n}$. It defines a set of upper-triangular  factorizing $F$-matrices. In particular, we have,
\begin{eqnarray*}
\overline{F}_{1,\ 2 \dots\ n}\ &=&\ {\cal C}_{1\ \dots\ n}\ F_{1,\ 2 \dots\ n}\ {\cal C}_{1\ \dots\ n}\ , \\
&=&\ e_1^{(22)}\ +\ e_1^{(11)}\ R_{1,\ 2 \dots\ n}\ , \\ 
&=&\ \left(
\begin{array}{cc}
A_{2\ \dots\ n}\ (\lambda_1\ ;\ \lambda_2 \dots\ ,\ \lambda_n) & B_{2\ \dots\ n}\ (\lambda_1\ ;\ \lambda_2 \dots\ ,\ \lambda_n) \\
{\mathbf 0}  & {\mathbf 1}
\end{array}
\right)_{[1]}\ .
\end{eqnarray*}
\end{prop}
\Proof
Straightforward from the invariance of the $R$-matrix under the action of ${\cal C}$.
\QED
\begin{prop}
$G_{1\ \dots\ n}\ =\ (F_{n\ \dots\ 1}^{t_n\ \dots\ t_1})^{-1}$ define a set of factorizing $F$-matrices.
\label{prop:gn}
\end{prop}
\Proof
We use Proposition (\ref{prop:factfz}).\\
The invertibility of $G$ follows from the one of $F$. Then we first prove that $G$'s indeed factorize the $R$-matrices. Let us first note that,
\begin{eqnarray*}
G_{1,\ 2\ \dots\ n}\ &=&\ (F_{n\ \dots\ 2,\ 1}^{t_n\ \dots\ t_1})^{-1}\ , \\
G_{2\ \dots\ n,\ 1}\ &=&\ (F_{1,\ n\ \dots\ 2}^{t_n\ \dots\ t_1})^{-1}\ .
\end{eqnarray*}
For the decomposition of $R_{1,\ 2\ \dots\ n}$, we have,
\begin{eqnarray*}
R_{1,\ 2\ \dots\ n}\ G_{1,\ 2\ \dots\ n}^{-1}\ &=&\ R_{1,\ 2\ \dots\ n}\ F_{n\ \dots\ 2,\ 1}^{t_n\ \dots\ t_1}\ , \\
&=&\ (R_{n\ \dots\ 2,\ 1}^{t_n\ \dots\ t_1})^{-1}\ F_{n\ \dots\ 2,\ 1}^{t_n\ \dots\ t_1}\ , \\
&=&\ (F_{n\ \dots\ 2,\ 1}\ R_{n\ \dots\ 2,\ 1}^{-1} )^{t_n\ \dots\ t_1}\ ,\\
&=&\ F_{1,\ n\ \dots\ 2}^{t_n\ \dots\ t_1}\ ,\\
&=&\ G_{2\ \dots\ n,\ 1}^{-1}\ ,
\end{eqnarray*}
where we used the relation,
\begin{eqnarray*}
R_{n\ \dots\ 2,\ 1}^{t_n\ \dots\ t_1}\ &=&\ (R_{n1}\ \dots\ R_{21})^{t_n\ \dots\ t_1}\ , \\
&=&\ R_{21}\ \dots\ R_{n1}\ , \\ 
&=&\ R_{2\ \dots\ n,\ 1}\ =\  R_{1,\ 2\ \dots\ n}^{-1}\ .
\end{eqnarray*}
The symmetry properties follow from the same argument. For the cyclic permutation $\sigma_c$, we have,
\begin{eqnarray*}
R_{1,\ 2\ \dots\ n}\ G_{0,\ 1\ \dots\ n}\ &=&\ R_{1,\ n\ \dots\ 2}^{t_1\ \dots\ t_n}\ (F_{n\ \dots\ 1,\ 0}^{t_0\ t_1\ \dots\ t_n})^{-1}\ , \\
&=&\ [(R_{n\ \dots\ 2,\ 1}\ F_{n\ \dots\ 1,\ 0})^{t_0\ \dots\ t_n}]^{-1}\ , \\
&=&\ (F_{1n\ \dots\ 2,\ 0}^{t_1\ \dots\ t_n})^{-1}\ , \\
&=&\ G_{0,\ 2\ \dots\ n1}\ ,
\end{eqnarray*}
as it should be and similarly for the permutation of the two first elements $\sigma_{12}$,
\begin{eqnarray*}
R_{12}\ G_{0,\ 1\ \dots\ n}\ &=&\ R_{12}^{t_1 t_2}\ (F_{n\ \dots\ 1,\ 0}^{t_0\ t_1\ \dots\ t_n})^{-1}\ , \\
&=&\ [(R_{21}\ F_{n\ \dots\ 1,\ 0})^{t_0\ \dots\ t_n}]^{-1}\ , \\
&=&\ [F_{n\ \dots\ 312,\ 0}^{t_0\ \dots\ t_n}]^{-1}\ , \\
&=&\ G_{0,\ 213\ \dots\ n}\ .
\end{eqnarray*}
It remains to prove the cocycle relation for $G$. It follows from the one for $F$ by taking its total transposition and then its inverse.
\QED

%% file: sect4
\section{Factorizing $F$-matrices and Algebraic Bethe Ansatz}

The purpose of this section is to use the factorizing $F$-matrices constructed in section 3 in the context of Algebraic Bethe Ansatz for the $XXX$ and $XXZ$ spin-$1 \over 2$ quantum chains of lenght $N$. As first application, we will show that the generating matrix of scalar products of quantum states, namely, $<0\ |T_1 (\lambda_1)\ \dots\ T_n (\lambda_n)|\ 0>$, where $|\ 0>$ is the (ferromagnetic) pseudo-vacuum reference state, and $T_i (\lambda_i)$ the quantum monodromy matrix associated to the (auxiliary) space $V_i (\lambda_i)$,  can be diagonalized by the total $F$-matrix $F_{1\ \dots\ n}$ defined in (\ref{eq:frn}, \ref{eq:fln}). Moreover, the (operator) diagonal entries of the quantum monodromy matrix, namely the operators $A_{1\ \dots\ N}\ (\lambda_0)$ and $D_{1\ \dots\ N}\ (\lambda_0)$ can also be diagonalized by means of $\overline{F}_{1\ \dots\ N}$ and $F_{1\ \dots\ N}$ respectively, independently of $\lambda_0$. This is described in section 4.1. To make these properties effective in applications, additional properties of the $F$-matrices $\overline{F}_{1\ \dots\ N}$ and $F_{1\ \dots\ N}$ are also derived in section 4.2. This include in particular the behaviour of $F$-matrices by inversion and crossing transformations, leading to simple difference equations for the partial $F$-matrices.

\subsection{Factorizing $F$-matrices and the Quantum Monodromy Matrix}

Let us first give two general properties of factorizing $F$-matrices when applied to the quantum monodromy matrix.
\begin{prop}
Let $T_i (\lambda_i) \in End\ (V_i (\lambda_i)) \otimes {\cal A}_R$ such that,
\begin{equation}
R_{ij}\ (\lambda_i,\ \lambda_j)\ T_i (\lambda_i)\ T_j (\lambda_j)\ =\ T_j (\lambda_j)\ T_i (\lambda_i)\ R_{ij}\ (\lambda_i,\ \lambda_j)\ , 
\label{eq:rttij}
\end{equation}
where, $R_{ij} \ (\lambda_i,\ \lambda_j) \in End\ (V_i (\lambda_i) \otimes V_j (\lambda_j))$ is a unitary $R$-matrix satisfying the Yang Baxter equation. Then for any permutation $\sigma \in S_n$, we have,
\begin{equation}
R_{1\ \dots\ n}^{\sigma}\ (\lambda_1,\ \dots\ ,\ \lambda_n)\ T_1 (\lambda_1)\ \dots\ T_n (\lambda_n)\ =\ T_{\sigma (1)} (\lambda_{\sigma (1)})\ \dots\ T_{\sigma (n)} (\lambda_{\sigma (n)})\ R_{1\ \dots\ n}^{\sigma}\ (\lambda_1,\ \dots\ ,\ \lambda_n)\ .
\end{equation}
We  define,
\begin{eqnarray*}
T_{1\ \dots\ n}\ =\ T_{1\ \dots\ n}\ (\lambda_1,\ \dots\ ,\ \lambda_n)\ =\ T_1 (\lambda_1)\ \dots\ T_n (\lambda_n)\ .
\end{eqnarray*}
Let $F_{12\ \dots\ n}$ be a set of factorizing $F$-matrices associated to the $R$-matrix. Define, 
\begin{eqnarray*}
\tilde{T}_{1\ \dots\ n}\ =\ F_{1\ \dots\ n}\ T_{1\ \dots\ n}\ F_{1\ \dots\ n}^{-1}\ , 
\end{eqnarray*}
it verifies for any element $\sigma \in S_n$,
\begin{eqnarray*}
\tilde{T}_{1\ \dots\ n}\ =\ \tilde{T}_{\sigma (1)\ \dots\ \sigma(n)}\ ,
\end{eqnarray*}
hence giving a completly symmetric presentation of the operators entries of the product of monodromy matrices $T_1 (\lambda_1)\ \dots\ T_n (\lambda_n)$.
\label{prop:ftf}
\end{prop}
\Proof
It is sufficient to remark that $\hat{\cal R}_{1\ \dots\ n}^{\sigma}$ commutes with the product $T_1 (\lambda_1)\ \dots\ T_n (\lambda_n)$, leading to,
\begin{eqnarray*}
{\cal P}_{1\ \dots\ n}^{\sigma^{-1}}\ R_{1\ \dots\ n}^{\sigma}\ T_{1\ \dots\ n}\ &=&\ T_{1\ \dots\ n}\ {\cal P}_{1\ \dots\ n}^{\sigma^{-1}}\ R_{1\ \dots\ n}^{\sigma}\ , 
\end{eqnarray*}
hence leading to the relation,
\begin{eqnarray*}
R_{1\ \dots\ n}^{\sigma}\ T_{1\ \dots\ n}\ &=&\ T_{\sigma(1)\ \dots\ \sigma(n)}\ R_{1\ \dots\ n}^{\sigma}\ , 
\end{eqnarray*}
which, using the factorization of $R_{1\ \dots\ n}^{\sigma}$ in terms of $F_{1\ \dots\ n}$, can be rewritten as,
\begin{eqnarray*}
F_{1\ \dots\ n}\ T_{1\ \dots\ n}\ F_{1\ \dots\ n}^{-1}\ =\ F_{\sigma(1)\ \dots\ \sigma(n)}\ T_{\sigma(1)\ \dots\ \sigma(n)}\ F_{\sigma(1)\ \dots\ \sigma(n)}^{-1}\ .
\end{eqnarray*}
\QED
\begin{prop}
Let us consider quantum $L$ -operators $L_{0i}\ (\lambda_0,\ \lambda_i) \in End\ (V_0 (\lambda_0) \otimes V_i (\lambda_i)$ verifying,
\begin{eqnarray*}
R_{ij}\ (\lambda_i,\ \lambda_j)\ L_{ik}\ (\lambda_i,\ \lambda_k)\ L_{jk}\  (\lambda_j,\ \lambda_k)\ &=&\ L_{jk}\ (\lambda_j,\ \lambda_k)\ L_{ik}\  (\lambda_i,\ \lambda_k)\ R_{ij}\ (\lambda_i,\ \lambda_j)\ , \\
R_{jk}\ (\lambda_j,\ \lambda_k)\ L_{ik}\ (\lambda_i,\ \lambda_k)\ L_{ij}\ (\lambda_i,\ \lambda_j)\ &=&\ L_{ij}\ (\lambda_i,\ \lambda_j)\ L_{ik}\ (\lambda_i,\ \lambda_k)\ R_{jk}\ (\lambda_j,\ \lambda_k)\ .
\end{eqnarray*}
The associated monodromy matrix is defined as,
\begin{equation}
T_{0,\ 1\ \dots\ N}\ (\lambda_0\ ;\ \lambda_1,\ \dots\ ,\ \lambda_N)\ =\ L_{0N}\ (\lambda_0,\ \lambda_N)\ \dots\ L_{01}\ (\lambda_0,\ \lambda_1)\ .
\end{equation}
Let,
\begin{eqnarray*}
\tilde{T}_{0,\ 1\ \dots\ N}(\lambda_0\ ;\ \lambda_1,\ \dots\ ,\ \lambda_N)\ =\ F_{1\ \dots\ N}\ T_{0,\ 1\ \dots\ N}\ (\lambda_0\ ;\ \lambda_1,\ \dots\ ,\ \lambda_N)\ F_{1\ \dots\ N}^{-1}\ .
\end{eqnarray*}
It satisfies,
\begin{equation}
\tilde{T}_{0,\ 1\ \dots\ N}(\lambda_0\ ;\ \lambda_1,\ \dots\ ,\ \lambda_N)\ =\ \tilde{T}_{0,\ \sigma(1)\ \dots\ \sigma(N)}(\lambda_0\ ;\ \lambda_{\sigma(1)},\ \dots\ ,\ \lambda_{\sigma(N)})\ .
\end{equation}
\label{prop:ftlf}
\end{prop}
\Proof
Straightforward from the definitions
\QED
As it is well known (see for example \cite{IK1}), for the two-dimensional quantum integrable models associated to the $R$-matrix (\ref{eq:vi}), the quantum monodromy matrix is given as the following $2 \times 2$ matrix of operators,
\begin{eqnarray*}
T_i\ (\lambda_i)\ =\ 
\left( 
\begin{array}{cc}
A\ (\lambda_i) & B\ (\lambda_i) \\
C\ (\lambda_i) & D\ (\lambda_i)
\end{array}
\right)_{[i]}\ .
\end{eqnarray*}
It satisfies in particular Proposition (\ref{prop:ftf}). One of the crucial ingredient of the Algebraic Bethe Ansatz method is the existence of a reference state (or pseudo-vacuum state) denoted $|\ 0>$, and its dual $<0\ |$, such that we have the following action of the operator entries of the quantum monodromy matrix,  
\begin{eqnarray*}
\begin{array}{cc}
A\ (\lambda)\ |\ 0>\ =\ a\ (\lambda)\ |\ 0>\ , & <0\ |\ A\ (\lambda)\ =\ a\ (\lambda)\ <0\ | \ , \\
B\ (\lambda)\ |\ 0>\ \ne\ 0\ , & <0\ |\ B\ (\lambda)\ =\ 0 \ ,\\
C\ (\lambda)\ |\ 0>\ =\ 0 \ , & <0\ |\ C\ (\lambda)\  \ne\ 0\ ,\\
D\ (\lambda)\ |\ 0>\ =\ d\ (\lambda)\ |\ 0>\ , & <0\ |\ D\ (\lambda)\ =\ d\ (\lambda)\ <0\ |\ . 
\end{array}
\end{eqnarray*}

As a first application of factorizing $F$-matrices constructed in the last section, we have the following proposition.
\begin{prop}
Let $F_{1\ \dots\ n}\ (\lambda_1,\ \dots\ ,\ \lambda_n)$ be defined as in Proposition (\ref{prop:exf}). We have,
\begin{eqnarray*}
& &<0\ |\ T_1\ (\lambda_1)\ \dots\ T_n\ (\lambda_n)\ |\ 0>\ =\\ 
&=&\ F_{1\ \dots\ n}^{-1}\ (\lambda_1\ \dots\ \lambda_n)\ \otimes_{i\ =\ 1}^n\ 
<0\ |\ T_i\ (\lambda_i)\ |\ 0>\ F_{1\ \dots\ n}\ (\lambda_1\ \dots\ \lambda_n)\ ,
\end{eqnarray*}
with,
\begin{eqnarray*}
<0\ |\ T_i\ (\lambda_i)\ |\ 0>\ &=&\ \left( 
\begin{array}{cc}
a\ (\lambda_i) & 0\\
0 & d\ (\lambda_i)
\end{array}
\right)_{[i]}\ .
\end{eqnarray*}
Hence, the factorizing $F$-matrices defined in Proposition (\ref{prop:exf}) diagonalize the matrix generating all scalar products of quantum states.
\label{prop:diagtq}
\end{prop}
\Proof
Let us first remark that the $2^n \times 2^n$ matrix $<0\ |T_1\ (\lambda_1)\ \dots\ T_n\ (\lambda_n)|\ 0>$ is a lower triangular matrix while the matrix $<0\ |T_n (\lambda_n)\ \dots\ T_1\ (\lambda_1)|\ 0>$ is an upper triangular matrix. These two properties are easily proved by induction on $n$.\\
Then, $F_{1\ \dots\ n}$ being a lower triangular matrix, it follows that the matrix,
\begin{eqnarray*}
F_{1\ \dots\ n}\ (\lambda_1\ \dots\ \lambda_n)\ <0\ |\ T_1\ (\lambda_1)\ \dots\ T_n\ (\lambda_n)\ |\ 0>\  F_{1\ \dots\ n}^{-1}\ (\lambda_1\ \dots\ \lambda_n)\ , 
\end{eqnarray*}
is also a lower triangular matrix which is invariant under any simultaneous permutation of the spaces $(1\ \dots\ n)$ and the corresponding spectral parameters, due to Propositions (\ref{prop:exf}, \ref{prop:ftf}). Hence, in particular, it is equal to the matrix,
\begin{eqnarray*}
F_{n\ \dots\ 1}\ (\lambda_n\ \dots\ \lambda_1)\ <0\ |\ T_n\ (\lambda_n)\ \dots\ T_1\ (\lambda_1)\ |\ 0>\  F_{n\ \dots\ 1}^{-1}\ (\lambda_n\ \dots\ \lambda_1)\ .
\end{eqnarray*}
But this matrix is a product of three upper triangular matrices, and therefore is also upper triangular. 
We conclude from the two assertions that it is a diagonal matrix. Moreover, it has then to be equal to the diagonal part of the matrix $<0\ |T_1 (\lambda_1)\ \dots\ T_n (\lambda_n)|\ 0>$, leading to the result.
\QED
\begin{rem}
The matrices $<0\ |T_1\ (\lambda_1)\ \dots\ T_n\ (\lambda_n)|\ 0>$ generate all scalar products of quantum states. In particular we have,
\begin{eqnarray*}
& &<0\ |\ C\ (\lambda_1)\ \dots\ C\ (\lambda_n)\ B\ (\lambda_{n+1})\ \dots\ B\ (\lambda_{2n})\ |\ 0>\ =\\ 
&=&\ tr_{1\ \dots\ 2n}\ (\Lambda_{1\ \dots\ 2n}^{+\ -}\ <0\ |\ T_1 (\lambda_1)\ \dots\ T_{2n} (\lambda_{2n})\ |\ 0>)\ , \\
&=&\ tr_{1\ \dots\ 2n}\ ( F_{1\ \dots\ 2n}\ (\lambda_1\ \dots\ \lambda_{2n})\ \Lambda_{1\ \dots\ 2n}^{+\ -}\ F_{1\ \dots\ 2n}^{-1}\ (\lambda_1\ \dots\ \lambda_{2n})\ {\cal D}_{1\ \dots\ 2n})\ , 
\end{eqnarray*}
where we have defined the matrices,
\begin{eqnarray*}
\Lambda_{1\ \dots\ 2n}^{+\ -}\ &=&\ \sigma_1^+\ \dots\ \sigma_n^+\ \sigma_{n+1}^-\ \dots\ \sigma_{2n}^-\ , \\
{\cal D}_{1\ \dots\ 2n}\ &=&\ \otimes_{i\ =\ 1}^n\ \left(  
\begin{array}{cc}
a\ (\lambda_i) & 0\\
0 & d\ (\lambda_i)
\end{array}
\right)_{[i]}\ .
\end{eqnarray*}
Thus above diagonalization separates  the contributions to the scalar products of quantum states coming on the one hand from the non trivial commutation relations of the quantum monodromy matrix algebra (contained in $F$) and on the other hand from the specific model (contained in the $<0\ |\ T_i (\lambda_i)\ |\ 0>$'s). This should be useful in the computation of correlation functions in the framework developped by Izergin and Korepin \cite{IK1}. In particular, let us note that the space dependence of the correlators is entirely contained in the pseudo-vacuum expectation values of individual $T_i^{x\ y}\ (\lambda_i)$'s  corresponding to the monodromy operator from site $x$ to site $y$ on the lattice. In our formula, these contributions appear only through diagonal matrices.
\end{rem}
\begin{rem}
Note that the above Proposition applies not only to the Heisenberg quantum spin chains but also to other integrable models associated to the same $R$-matrix and solvable by the Algebraic Bethe Ansatz, like for example the Sine-Gordon model.
\end{rem}

The next application of the factorizing $F$-matrices $F_{1\ \dots\ N}$ is specific to the $XXX$ and $XXZ$ spin-$1 \over 2$ inhomogeneous  quantum chains of lenght $N$. From the point of view of the Quantum inverse Scattering Method, these models are {\em fundamental} models associated to the $R$-matrices defined in (\ref{eq:vi}). Indeed, to each site $i$ of the lattice, $i\ =\ 1,\ \dots\ N$, we can associate a quantum $L$-operator given simply in terms of the corresponding $R$-matrix as,
\begin{equation}
L_{0i}\ (\lambda_0,\ \xi_i)\ =\ R_{0i}\ (\lambda_0,\ \xi_i)\ .
\label{eq:LR}
\end{equation}
We will refer to the parameter $\lambda_0$ as the spectral parameter and to the parameters $\xi_i$ as the inhomogeneity parameters.\\
The quantum monodromy operator is then,
\begin{eqnarray*}
T_{0,\ 1\ \dots\ N}\ (\lambda_0\ ;\ \xi_1,\ \dots\ ,\ \xi_N)\ &=&\ L_{0N}\ (\lambda_0,\ \xi_N)\ \dots\ L_{01}\ (\lambda_0,\ \xi_1)\ , \\
&=&\ R_{0,\ 1\ \dots\ N}\ (\lambda_0\ ;\ \xi_1,\ \dots\ ,\ \xi_N)
\end{eqnarray*}
Hence the quantum space of states of the model is ${\cal H}_N\ \simeq\ \otimes_{i\ =\ 1}^{N}\ {\cal H}_i$, where ${\cal H}_i\ \simeq {\mathbf C}^2$ for all $i$.\\
For later use we  define the ferromagnetic reference state $|\ 0>$ (hereafter also called the pseudo-vacuum) as the tensor product,
\begin{eqnarray*}
|\ 0>\ &=&\ \otimes_{i\ =\ 1}^N\ |\ 0>_i \ ,\\
|\ 0>_i\ &=&\ \left( 
\begin{array}{c}
1 \\
0 
\end{array}
\right)_{[i]}\ .
\end{eqnarray*}
For these models we have,
\begin{eqnarray*}
\begin{array}{cc}
a\ (\lambda)\ =\ \prod_{i\ =\ 1}^{N}\ a_i\ (\lambda,\ \xi_i)\ , & a_i\ (\lambda,\ \xi_i)\ =\ 1\ ,\\
d\ (\lambda)\ =\ \prod_{i\ =\ 1}^{N}\ d_i\ (\lambda,\ \xi_i)\ , & d_i\ (\lambda,\ \xi_i)\ =\ b\ (\lambda,\ \xi_i)\ .
\end{array}
\end{eqnarray*}
Using the notations introduced in the last section, we have :
\begin{prop}
Let $F_{1\ \dots\ N}\ (\xi_1,\  \dots\ ,\ \xi_N)$ be the total $F$ operator acting in ${\cal H}_N$. Let its conjugated matrix by ${\cal C}_{1\ \dots\ N}\ =\ \otimes_{i\ =\ 1}^N\ \sigma_i^x$ be denoted by  $\overline{F}_{1\ \dots\ N}\ (\xi_1,\  \dots\ ,\ \xi_N)$ . They satisfy for any spectral parameter $\lambda$,
\begin{eqnarray*}
F_{1\ \dots\ N}\ (\xi_1,\  \dots\ ,\ \xi_N)\ D_{1\ \dots\ N}\ (\lambda\ ;\ \xi_1,\  \dots\ ,\ \xi_N)\ F_{1\ \dots\ N}^{-1}\ (\xi_1,\  \dots\ ,\ \xi_N)\ =\\ =\ \otimes_{i\ =\ 1}^{N}\ 
\left( 
\begin{array}{cc}
b\ (\lambda,\ \xi_i) & 0\\
0 & 1
\end{array}
\right)_{[i]}\ ,
\end{eqnarray*}
\begin{eqnarray*}
\overline{F}_{1\ \dots\ N}\ (\xi_1,\  \dots\ ,\ \xi_N)\ A_{1\ \dots\ N}\ (\lambda\ ;\ \xi_1,\  \dots\ ,\ \xi_N)\ \overline{F}_{1\ \dots\ N}^{-1}\ (\xi_1,\  \dots\ ,\ \xi_N)\ =\\ =\ \otimes_{i\ =\ 1}^{N}\ 
\left( 
\begin{array}{cc}
1 & 0\\
0 & b\ (\lambda,\ \xi_i)
\end{array}
\right)_{[i]}\ , 
\end{eqnarray*}
and similarly,
\begin{eqnarray*}
F_{1\ \dots\ N}\ (\xi_1,\  \dots\ ,\ \xi_N)\ A_{1\ \dots\ N}^{\dagger}\ (\lambda\ ;\ \xi_1,\  \dots\ ,\ \xi_N)\ F_{1\ \dots\ N}^{-1}\ (\xi_1,\  \dots\ ,\ \xi_N)\ =\\ =\ \otimes_{i\ =\ 1}^{N}\ 
\left( 
\begin{array}{cc}
1 & 0\\
0 & b\ (\xi_i,\ \lambda)
\end{array}
\right)_{[i]}\ , 
\end{eqnarray*}
\begin{eqnarray*}
\overline{F}_{1\ \dots\ N}\ (\xi_1,\  \dots\ ,\ \xi_N)\ D_{1\ \dots\ N}^{\dagger}\ (\lambda\ ;\ \xi_1,\  \dots\ ,\ \xi_N)\ \overline{F}_{1\ \dots\ N}^{-1}\ (\xi_1,\  \dots\ ,\ \xi_N)\ =\\ =\ \otimes_{i\ =\ 1}^{N}\ 
\left( 
\begin{array}{cc}
b\ (\xi_i,\ \lambda) & 0\\
0 & 1
\end{array}
\right)_{[i]}\ .
\end{eqnarray*}
\label{prop:fDA+}
\end{prop}
\Proof
It is enough to prove the above relation for $D_{1\ \dots\ N}$. We will do the proof by induction on N. The relation is true for $N\ =\  2$ as can be checked by direct computations using the matrix $F_{12}$ given in (\ref{eq:f12}),  the formula,
\begin{eqnarray*}
D_{12}\ (\lambda_0\ ;\ \xi_1,\ \xi_2)\ =\ \left(
\begin{array}{cccc}
b\ (\lambda_0\ ,\ \xi_1) \ b\ (\lambda_0\ ,\ \xi_2) & 0 & 0 & 0\\
0 & b\ (\lambda_0\ ,\ \xi_1) & 0 & 0\\
0 & c\ (\lambda_0\ ,\ \xi_1)\ c\ (\lambda_0\ ,\ \xi_2) & b\ (\lambda_0\ ,\ \xi_2) & 0\\
0 & 0 & 0 & 1
\end{array}
\right)_{[12]}\ , 
\end{eqnarray*}
and the fact that,
\begin{eqnarray*}
b\ (\lambda_0\ ,\ \xi_1)\ c\ (\xi_1,\ \xi_2)\ +\ c\ (\lambda_0\ ,\ \xi_1)\ c\ (\lambda_0\ ,\ \xi_2)\ b\ (\xi_1,\ \xi_2)\ =\ b\ (\lambda_0\ ,\ \xi_2)\ c\ (\xi_1,\ \xi_2)\ , 
\end{eqnarray*}
due to the Yang-Baxter equation for the $R$-matrix.\\
Let the property be verified for $n-1$ spaces. Taking into account the decomposition of the total $F$-matrix as $F_{1\ \dots\ n}\ =\ F_{2\ \dots\ n}\ F_{1,\ 2\ \dots\ n}$, we want to compute $F_{1\ \dots\ n}\ D_{1\ \dots\ n}\ (\lambda_0)$. Let us consider this matrix product as a $2 \times 2$ matrix in space one whose  matrix entries are $2^{n-1} \times 2^{n-1}$ matrices. We get, using simplified notations, namely not writing explicitly the spectral parameter dependence $(\lambda_1,\ \dots\ ,\ \lambda_n)$,
\begin{eqnarray*}
 F_{1\ \dots\ n}\ D_{1\ \dots\ n}\ (\lambda_0)\ &=&\ F_{2\ \dots\ n}\ \left( \begin{array}{cc}
{\mathbf 1} & {\mathbf 0} \\
C_{2\ \dots\ n} (\lambda_1)  & D_{2\ \dots\ n}\ (\lambda_1)
\end{array}
\right)_{[1]}\ D_{1\ \dots\ n}\ (\lambda_0)\ ,\\
&=&\  F_{2\ \dots\ n}\ \left(
\begin{array}{cc}
X_{2\ \dots\ n} & {\mathbf 0} \\
Y_{2\ \dots\ n}  & Z_{2\ \dots\ n}
\end{array}
\right)_{[1]}\ ,
\end{eqnarray*}
where,
\begin{eqnarray*}
X_{2\ \dots\ n}\ &=&\ b\ (\lambda_0,\ \lambda_1)\ D_{2\ \dots\ n}\ (\lambda_0)\ ,\\
Y_{2\ \dots\ n}\ &=&\ b\ (\lambda_0,\ \lambda_1)\ C_{2\ \dots\ n}\ (\lambda_1)\ D_{2\ \dots\ n}\ (\lambda_0)\ +\ c\ (\lambda_0,\ \lambda_1)\ D_{2\ \dots\ n}\ (\lambda_1)\ C_{2\ \dots\ n}\ (\lambda_0)\ ,\\
Z_{2\ \dots\ n}\ &=&\ D_{2\ \dots\ n}\ (\lambda_1)\ D_{2\ \dots\ n}\ (\lambda_0)\ .
\end{eqnarray*}
Now using the Yang-Baxter equation written in simplified notations as,
\begin{eqnarray*}
R_{01}\ (\lambda_0,\ \lambda_1)\ R_{0,\ 2\ \dots\ n}\ (\lambda_0)\ R_{1,\ 2\ \dots\ n}\ (\lambda_1)\ = \ R_{1,\ 2\ \dots\ n}\ (\lambda_1)\  R_{0,\ 2\ \dots\ n}\ R_{01}\ (\lambda_0,\ \lambda_1)\ , 
\end{eqnarray*}
we obtain the following relation for $D$ and $C$ operators,
\begin{eqnarray*} 
b\ (\lambda_0,\ \lambda_1)\ C_{2\ \dots\ n}\ (\lambda_1)\ D_{2\ \dots\ n}\ (\lambda_0)\ +\ c\ (\lambda_0,\ \lambda_1)\ D_{2\ \dots\ n}\ (\lambda_1)\ C_{2\ \dots\ n}\ (\lambda_0)\ = \\ 
=\ D_{2\ \dots\ n}\ (\lambda_0)\  C_{2\ \dots\ n}\ (\lambda_1)\ ,
\end{eqnarray*}
and,
\begin{eqnarray*}
D_{2\ \dots\ n}\ (\lambda_0)\  D_{2\ \dots\ n}\ (\lambda_1)\  =\ D_{2\ \dots\ n}\ (\lambda_1)\  D_{2\ \dots\ n}\ (\lambda_0)\ .
\end{eqnarray*}
This leads to the relation,
\begin{eqnarray*}
F_{1\ \dots\ n}\ D_{1\ \dots\ n}\ (\lambda_0)\ &=&\ F_{2\ \dots\ n}\ \left( \begin{array}{cc} 
b\ (\lambda_0,\ \lambda_1) & 0 \\
0  & 1
\end{array}
\right)_{[1]}\ D_{2\ \dots\ n}\ (\lambda_0)\ F_{1,\ 2\ \dots\ n}\ ,\\
&=&\ \left( \begin{array}{cc}
b\ (\lambda_0,\ \lambda_1) & 0 \\
0  & 1
\end{array}
\right)_{[1]}\ F_{2\ \dots\ n}\ D_{2\ \dots\ n}\ (\lambda_0)\ F_{2\ \dots\ n}^{-1}\ F_{1\ \dots\ n}\ .
\end{eqnarray*}
Applying the recursion relation to $D_{2\ \dots\ n}\ (\lambda_0)$ leads to the desired result, namely,
\begin{eqnarray*}
F_{1\ \dots\ N}\ D_{1\ \dots\ N}\ (\lambda_0)\ &=&\ \otimes_{i\ =\ 1}^{N}\ \left( \begin{array}{cc}
b\ (\lambda_0,\ \lambda_i) & 0 \\
0  & 1
\end{array}
\right)_{[i]}\ F_{1\ \dots\ N}\ .
\end{eqnarray*}
The relation for $A_{1\ \dots\ N}$ follows from the above one by conjugaison by the operator ${\cal C}_{1\ \dots\ N}$. For $A^{\dagger}$ and $D^{\dagger}$ we use the Proposition (\ref{prop:crossrn}) from which we obtain that these operators are proportional to $D$ and $A$ respectively but with shifted spectral parameters. 
\QED
\begin{rem}
The above result shows that there exists a factorizing $F$-matrix diagonalizing  the operator $A$ and another (different) one diagonalizing  the operator $D$. A natural question is whether there exists a factorizing $F$-matrix diagonalizing the operator $A\ +\ D$ which contains the Hamiltonian and the series of conserved charges. We conjecture this question to have a positive answer.
\end{rem}
\begin{rem}
The effective diagonalization of the operator $D$ and $A$ would require that we are in fact able to obtain the inverse matrix $F_{1\ \dots\ N}^{-1}$. This is a $2^N \times 2^N$ matrix, where $N$ is the size of the lattice. Hence, even if it is a  lower triangular matrix, this is a non trivial problem. It will be solved in the next section.
\end{rem}

\subsection{Algebraic Properties of Factorizing $F$-matrices}

First we recall that following the Proposition (\ref{prop:daggerr}), we have,
\begin{eqnarray*}
F_{1\ \dots\ n}^\dagger\ (\lambda_1,\ \dots\ ,\ \lambda_n)\ =\ 
F_{1\ \dots\ n}^{t_1\ \dots\ t_n}\ (-\lambda_1,\ \dots\ ,\ -\lambda_n)\ .
\end{eqnarray*}
Let us also define the diagonal matrix,
\begin{eqnarray}
\Delta_{1\ \dots\ n}\ &=&\ \Delta_{1\ \dots n}\ (\lambda_1,\ \dots\ \lambda_n)\ , \nonumber\\
&=&\ \prod_{i<j}\ \Delta_{ij}(\lambda_i,\ \lambda_j)\ ,
\label{eq:deltan}
\end{eqnarray}
where,
\begin{eqnarray*}
\Delta_{ij}\ (\lambda_i,\ \lambda_j)\ =\ 
\left(
\begin{array}{cccc}
1 &  &   &  \\
& b^{-1}(\lambda_j , \lambda_i) &  &  \\
 &  &b^{-1}(\lambda_i, \lambda_j) &  \\
 &  &  & 1
\end{array}
\right)_{[ij]}\ .
\end{eqnarray*}
We have,
\begin{lem}
\begin{eqnarray*}
\Delta_{1\ \dots\ n}^{\dagger}\ &=&\ {\cal C}_{1\ \dots\ n}\ \Delta_{1\ \dots\ n}\ {\cal C}_{1\ \dots\ n}\ ,\\
\Delta_{\sigma(1)\ \dots\ \sigma(n)}\ &=&\ \Delta_{1\ \dots\ n}\ ,\ \forall \sigma \in S_n\ .
\end{eqnarray*}
\end{lem}
\Proof
It follows trivially from the same properties for the elementary matrix $\Delta_{12}$ which can be checked by direct computation.
\QED
\begin{prop}
 The following formula for the inverse of $F_{1\ \dots\ n}$ defined in Proposition (\ref{prop:exf}) holds, 
\begin{equation}
F_{1\ \dots\ n}^{-1}\ =\  {\cal C}_{1\ \dots\ n}\ F_{1\ \dots\ n}^{\dagger}\ {\cal C}_{1\ \dots\ n}\ \Delta_{1\ \dots\ n}\ .
\label{eq:ff+}
\end{equation}
\label{prop:f-1}
\end{prop}
\Proof
By induction on $n$. The property is true for $n\ =\ 2$ as can be checked by direct computation. Suppose it is true for $n$. Then we want to prove that,
\begin{eqnarray*}
F_{01\ \dots\ n}^{-1}\ =\ {\cal C}_{01\ \dots\ n}\ F_{01\ \dots\ n}^{\dagger}\ {\cal C}_{01\ \dots\ n}\ \Delta_{01\ \dots\ n}\ .
\end{eqnarray*}
Factorizing $F_{01\ \dots\ n}$ as $F_{1\ \dots\ n}\ F_{0,\ 1\ \dots\ n}$ and assuming the property to be true for $n$, we can see it is enough to prove that,
\begin{equation}
F_{0,\ 1\ \dots\ n}\ {\cal C}_{01\ \dots\ n}\ F_{0,\ 1\ \dots\ n}^{\dagger}\ {\cal C}_{01\ \dots\ n}\ =\ F_{1\ \dots\ n}^{-1}\ \Delta_{0,\ 1\ \dots\ n}^{-1}\ F_{1\ \dots\ n}\ ,
\label{eq:pf-1}
\end{equation}
where we defined,
\begin{eqnarray*}
\Delta_{0,\ 1\ \dots\ n}\ =\ \Delta_{0 1\ \dots\ n}\ \Delta_{1\ \dots\ n}^{-1}\ .
\end{eqnarray*}
This equality can be proven by writing every object as a $2 \times 2$ matrix in space $0$. We have,
\begin{eqnarray*}
F_{0,\ 1\ \dots\ n}\ (\lambda_0)\ &=&\  \left(
\begin{array}{cc}
{\mathbf 1} & {\mathbf 0} \\
C\ (\lambda_0)  & D\ (\lambda_0)
\end{array}
\right)_{[0]}\ ,\\
F_{0,\ 1\ \dots\ n}^{\dagger}\ (\lambda_0)\ &=&\  \left(
\begin{array}{cc}
{\mathbf 1} & C^{\dagger}\ (\lambda_0) \\
{\mathbf 0} & D^{\dagger}\ (\lambda_0)
\end{array}
\right)_{[0]}\ ,
\end{eqnarray*}
with $C\ (\lambda_0)\ =\ C_{1\ \dots\ n}\ (\lambda_0\ ;\ \xi_1,\ \dots\ ,\ \xi_n)$ and similarly for $D\ (\lambda_0)$. Then using the unitarity of the $R$-matrix, we have, 
\begin{eqnarray*}
C\ ({\lambda}_0)\ A^{\dagger}\ ({\lambda}_0)\ +\ D\ ({\lambda}_0)\ B^{\dagger}\ ({\lambda}_0)\ =\ 0\ ,
\end{eqnarray*}
and hence we get,
\begin{eqnarray*}
F_{0,\ 1\ \dots\ n}\ {\cal C}_{01\ \dots\ n}\ F_{0,\ 1\ \dots\ n}^{\dagger}\ {\cal C}_{01\ \dots\ n}\ =\ \left(
\begin{array}{cc}
A^{\dagger}\ (\lambda_0) & {\mathbf 0} \\
{\mathbf 0} & D\ (\lambda_0)
\end{array}
\right)_{[0]}\ .
\end{eqnarray*}
In order to compute the right hand side of the equation (\ref{eq:pf-1}), let us first note that,
\begin{eqnarray*}
\Delta_{0,\ 1\ \dots\ n}^{-1}\ =\ \left(
\begin{array}{cc}
{\cal C}\ \tilde{D}^{\dagger}\ (\lambda_0)\ {\cal C}\  & {\mathbf 0} \\
{\mathbf 0} & \tilde{D}\ (\lambda_0)
\end{array}
\right)_{[0]}\ .
\end{eqnarray*}
Then it is sufficient to use formulas given in Proposition (\ref{prop:fDA+}) to obtain the desired result.
\QED
\begin{prop}
For the factorizing $F$-matrix defined in Proposition (\ref{prop:exf}), we have the following identity,
\begin{eqnarray*}
F_{1\ \dots\ n}\ F_{n\ \dots\ 1}^{t_1\ \dots\ t_n}\ =\ \Delta_{1\ \dots\ n}^{-1}\ . 
\end{eqnarray*}
It follows that the two factorizing $F$-matrices defined in the last section, namely $F_{1\ \dots\ n}$ in Proposition (\ref{prop:exf}) and $G_{1\ \dots\ n}$ in Proposition (\ref{prop:gn}), are related by,
\begin{eqnarray*}
G_{1\ \dots\ n}\ =\ \Delta_{1\ \dots\ n}\ F_{1\ \dots\ n}\ .
\end{eqnarray*}
\label{prop:fftd}
\end{prop}
\Proof
This follows from the fact that,
\begin{eqnarray*}
F_{1\ \dots\ n}^{\dagger}\ =\ {\cal C}_{1\ \dots\ n}\ F_{n\ \dots\ 1}^{t_1\ \dots\ t_n}\ {\cal C}_{1\ \dots\ n}\ , 
\end{eqnarray*}
which is proven by induction on $n$. Let this property be true for $n-1$. Then, to prove it for $n$ we just need to prove the identity,
\begin{eqnarray*}
F_{1,\ 2\ \dots\ n}^{\dagger}\ =\ {\cal C}_{1\ \dots\ n}\  F_{n\ \dots\ 2,\  1}^{t_1\ \dots\ t_n}\ {\cal C}_{1\ \dots\ n}\ .
\end{eqnarray*} 
We have,
\begin{eqnarray*}
F_{1,\ 2\ \dots\ n}^{\dagger}\ &=&\ (e_1^{(11)}\ +\ e_1^{(22)}\ R_{1,\ 2\ \dots\ n})^{\dagger}\ ,\\
&=&\ e_1^{(11)}\ + R_{1,\ 2\ \dots\ n}^{\dagger}\ e_1^{(22)}\ .
\end{eqnarray*}
Similarly,
\begin{eqnarray*}
{\cal C}_{1\ \dots\ n}\  F_{n\ \dots\ 2,\  1}^{t_1\ \dots\ t_n}\ {\cal C}_{1\ \dots\ n}\ &=&\ {\cal C}_{1\ \dots\ n}\  (e_1^{(22)}\ +\ e_1^{(11)}\ R_{n\ \dots\ 2,\ 1})^{t_1\ \dots\ t_n}\ {\cal C}_{1\ \dots\ n}\ ,\\
&=&\ {\cal C}_{1\ \dots\ n}\  (e_1^{(22)}\ +\ R_{n\ \dots\ 2,\ 1}^{t_1\ \dots\ t_n}\ e_1^{(11)})\ {\cal C}_{1\ \dots\ n}\ ,\\
&=&\ e_1^{(11)}\ +\ R_{2\ \dots\ n,\ 1}\ e_1^{(22)}\ ,
\end{eqnarray*}
leading to the result. 
\QED
We are now interested in the behaviour of the $F$-matrices given in Proposition (\ref{prop:exf}) under crossing symmetry. We have the following Proposition.
\begin{prop}
Under the transformation $\lambda_0 \rightarrow \lambda_0^s$ the partial $F$-matrix $\tilde{F}_{0,\ 1\ \dots\ n}$ following from Proposition (\ref{prop:exf}) verifies the relation,
\begin{equation}
\tilde{F}_{0,\ 1\ \dots\ n}\ (\lambda_0)\ \gamma_0\ \tilde{F}_{0,\ 1\ \dots\ n}^{t_0}\ (\lambda_0^s)\ \gamma_0\ =\ \Theta_0\ \Delta_{0,\ 1\ \dots\ n}^{-1}\ ,
\end{equation}
where, $\tilde{F}_{0,\ 1\ \dots\ n}\ =\ F_{1\ \dots\ n}\ F_{0,\ 1\ \dots\ n}\ F_{1\ \dots\ n}^{-1}$ and $\Theta_0$ is the following  $2 \times 2$ matrix acting in space $0$,
\begin{eqnarray*}
\Theta_0\ =\ \left(
\begin{array}{cc}
\rho\ (\lambda_0;\ \lambda_1\ \dots\ \lambda_n) &  0 \\
0 & 1 
\end{array}
\right)_{[0]} \ .
\end{eqnarray*}
\label{prop:fcross}
\end{prop}
\Proof
We have, using simplified notations as $q \equiv (1\ \dots\ n)$,
\begin{eqnarray*}
F_{0,\ q}\ (\lambda_0^s)\ &=&\ e_0^{(11)}\ +\ e_0^{(22)}\ R_{0,\ q}\ (\lambda_0^s) \\
&=&\ e_0^{(11)}\ +\ \rho\ (\lambda_0)\ e_0^{(22)}\ \gamma_0\ R_{q,\ 0}^{t_0}\ (\lambda_0)\ \gamma_0\ .
\end{eqnarray*}
Then using the fact that $\gamma_0\ e_0^{(11)}\ \gamma_0\ =\ e_0^{(22)}$, we get,
\begin{eqnarray*}
F_{0,\ 1\ \dots\ n}\ (\lambda_0)\ \gamma_0\ F_{0,\ 1\ \dots\ n}^{t_0}\ (\lambda_0^s)\ \gamma_0\ &=&\ \rho\ (\lambda_0)\ e_0^{(11)}\ R_{q,\ 0}\ (\lambda_0)\ e_0^{(11)}\ +\ e_0^{(22)}\ R_{0,\ q}\ (\lambda_0)\ e_0^{(22)}\ ,\\
&=&\ \Theta_0\ \left(
\begin{array}{cc}
A^{\dagger}\ (\lambda_0) & {\mathbf 0} \\
{\mathbf 0} & D\ (\lambda_0)
\end{array}
\right)_{[0]}\ ,
\end{eqnarray*}
leading to the result.
\QED
We also have the following difference equation for the $F$-matrices, involving twice the shift associated to crossing and denoted by  $\lambda_0 \rightarrow \lambda_0^{ss}$.
\begin{prop}
\begin{eqnarray*}
R_{0,\ q}^{t_0}\ (\lambda_0^{ss})\ F_{q,\ 0}^{t_0}\ (\lambda_0)\ =\ F_{0,\ q}^{t_0}\ (\lambda_0^{ss})\ ( \rho\ (\lambda_0^{s})\ \rho^{-1}\ (\lambda_0)\ e_0^{(11)}\ +\ e_0^{(22)})\ ,
\end{eqnarray*}
where we used the simplified notations with $q \equiv (1\ \dots\ n)$.
\end{prop}
\Proof
We have from the crossing symmetry of the $R$-matrix,
\begin{eqnarray*}
R_{0,\ q}\ (\lambda_0^s)\ =\ \rho\ (\lambda_0)\ \gamma_0\ (R_{0,\ q}^{-1}\ (\lambda_0))^{t_0}\ \gamma_0\ ,
\end{eqnarray*}
leading to the relation for double shifts,
\begin{eqnarray*}
(R_{0,\ q}^{t_0}\ (\lambda_0^{ss}))^{-1}\ =\ \rho^{-1}\ (\lambda_0^{s})\ \rho\ (\lambda_0)\ (R_{0,\ q}^{-1}\ (\lambda_0))^{t_0}\ .
\end{eqnarray*}
We have also from the definition of the $F$-matrices,
\begin{eqnarray*}
F_{q,\ 0}^{t_0}\ (\lambda_0)\ &=&\ e_0^{(22)}\ +\ (R_{0,\ q}^{-1}\ (\lambda_0))^{t_0}\ e_0^{(11)}\ ,\\
&=&\ e_0^{(22)}\ +\ \rho\ (\lambda_0^{s})\ \rho^{-1}\ (\lambda_0)\ (R_{0,\ q}^{t_0}\ (\lambda_0^{ss}))^{-1}\ e_0^{(11)}\ ,
\end{eqnarray*}
leading to the result.
\QED

%% file: sect5
\section{$F$-basis for the $XXX$-$1 \over 2$ quantum chain}

The purpose of this section is to study the operator algebra generated by the quantum monodromy matrix entries in the basis generated by the column vectors of the inverse of the $F$-matrix $F_{1\ \dots\ N}$ defined in Proposition (\ref{prop:exf}) acting in ${\cal H}_N$.\\
The idea we would like to advertize, here in this very simple example, is that this basis is particularly convenient in the sense that the explicit expressions for the operator entries of the quantum monodromy matrix are far more simple in this basis than in the original one. This will also results in a more explicit formula for these $F$-matrices themselves. In all this section, $F_{1\ \dots\ n}$ stands for the factorizing $F$-matrices defined in Proposition (\ref{prop:exf}).

\subsection{The $F$-basis and Algebraic Bethe Ansatz : Outlook}

The general philosophy of the Algebraic Bethe Ansatz is that one concentrate on the (quadratic) commutation relations of the operator entries of the quantum monodromy matrix rather than on their explicit action in the space of quantum state ${\cal H}_N$, which, in some sense, is obtain as a byproduct of the representation theory of this Yang-Baxter algebra (see for example \cite{KBI}). In fact the explicit expressions of the operator entries of the quantum monodromy matrix are too complicated to be used efficiently in this context. More precisely, in the case of the $XXX$-$1 \over 2$ model, the  quantum monodromy operator is a $2 \times 2$ matrix with entries $A,\ B,\ C,\ D$ which are  obtained as a sums of $2^{N-1}$ operators which themselves are products of $N$ local operators on the quantum chain. As an example, the $B$ operator is given as,
\begin{equation}
B_{1\ \dots\ N}\ (\lambda)\ =\ \sum_{i\ =\ 1}^{N}\ S_i^-\ \Omega_i\ +\ \sum_{i\ \ne j\ \ne\ k}\ S_i^-\ (S_j^-\ S_k^+)\ \Omega_{ijk}\ +\ higher\ terms\ ,
\label{eq:B2N}
\end{equation}
where the matrices $\Omega_i$, $\Omega_{ijk}$, are diagonal operators acting respectively on all sites but $i$, on all sites but $i, j, k$, and  the higher order terms involves more and more exchange spin terms like $S_j^-\ S_k^+$. It means that the $B$ operator return one spin somewhere on the chain, this operation being however dressed non-locally and with non-diagonal operators by multiple exchange terms of the type $S_j^-\ S_k^+$.\\
We first would like to argue that in the $F$-basis, the $B$ and the $C$ operators of the $XXX$-$1 \over 2$ model should have a much simpler expression than in the original basis. This is already true for the operator $D$ which becomes diagonal in the $F$-basis  since we have from the last section,
\begin{eqnarray*}
F_{1\ \dots\ N}\ (\xi_1,\  \dots\ ,\ \xi_N)\ D_{1\ \dots\ N}\ (\lambda\ ;\ \xi_1,\  \dots\ ,\ \xi_N)\ F_{1\ \dots\ N}^{-1}\ (\xi_1,\  \dots\ ,\ \xi_N)\ =\\ =\ \tilde{D}_{1\ \dots\ N}\ (\lambda\ ;\ \xi_1,\  \dots\ ,\ \xi_N)\ =\ \otimes_{i\ =\ 1}^{N}\ 
\left( 
\begin{array}{cc}
b\ (\lambda,\ \xi_i) & 0\\
0 & 1
\end{array}
\right)_{[i]}\ .
\end{eqnarray*}
Note here that for convenience  we will take the following conventions for the spectral parameter dependence : parameters associated to lattice sites, namely to the quantum space of states, will be denoted by $\xi_i$ while spectral parameters associated to the (so called auxiliary) matrix spaces will be denoted by parameters $\lambda_i$.
Using the unitarity and the crossing symmetries of the $R$-matrix we have also obtained that,
\begin{equation}
A^{\dagger}\ ({\lambda}_0)\ B\ ({\lambda}_0)\ +\ C^{\dagger}\ ({\lambda}_0)\ D\ ({\lambda}_0)\ =\ 0\ ,
\end{equation}
together with
\begin{equation}
C^{\dagger}\ ({\lambda}_0)\ =\  -\ \rho^{-1}\ (\lambda_0)\ B\ ({\lambda}_0^s)
\end{equation}
which leads to the fact that 
\begin{equation}
A^{\dagger}\ ({\lambda}_0)\ B\ ({\lambda}_0)\ =\ \rho^{-1}\ (\lambda_0)\ B\ ({\lambda}_0^s)\ D\ ({\lambda}_0)\ .
\end{equation}
Going to the $F$-basis, $A^{\dagger}$ and $D$ becomes diagonal matrices, leading to the relation,
\begin{equation}
\tilde{B}\ ({\lambda}_0)\ =\ \rho^{-1}\ (\lambda_0)\ {\cal C}\ (\tilde{D}^{\dagger})^{-1}\ (\lambda_0)\ {\cal C}\ \tilde{B}\ ({\lambda}_0^s)\ \tilde{D}\ (\lambda_0)\ .
\end{equation}
Hence the behaviour of the $B$ operator under the shift $\lambda_0 \rightarrow \lambda_0^s$ is given in the $F$-basis by the action of diagonal matrices. The same is true for the operator $C$. This means that their explicit expressions should simplify in this basis. 
Indeed, due to the fact that the total spin operator commutes with the $R$-matrix, we have the following formulas for the operators $A,\ B,\ C$ in terms of the operator $D$ and the total spin operator,
\begin{eqnarray}
C_{1\ \dots\ N}\ (\lambda)\ &=&\ \left[\ S_{1\ \dots\ N}^+\ ,\ D_{1\ \dots\ N}\ (\lambda)\ \right]\ , \nonumber\\
B_{1\ \dots\ N}\ (\lambda)\ &=&\ \left[\ D_{1\ \dots\ N}\ (\lambda)\ ,\ S_{1\ \dots\ N}^- \ \right]\ , \nonumber\\
A_{1\ \dots\ N}\ (\lambda)\ &=&\ D_{1\ \dots\ N}\ (\lambda)\ -\ \left[\ S_{1\ \dots\ N}^+\ ,\ \left[\ S_{1\ \dots\ N}^-\ ,\ D_{1\ \dots\ N}\ (\lambda)\  \right]\  \right]\ .
\label{eq:ABCDr}
\end{eqnarray}
This means that in the $F$-basis all the above operators are given in terms of a diagonal matrix $\tilde{D}$ and of the total spin operator in this basis. So, if the expression of the total spin operator in the $F$-basis is simple so will be the expressions of $\tilde{A},\ \tilde{B},\ \tilde{C}$. Since the $F$-matrices constructed in section 3 are functions of the $R$ matrix it is possible to compute explicitly $\tilde{S}$ the total spin operator in the basis given by $F$. This is done in the next section.

\subsection{The Total Spin Operator in the $F$-basis}

The total spin operator $S_{1\ \dots\ N}^{\alpha}$, $\alpha\ =\ 1,\ 2,\ 3$ is defined as,
\begin{equation}
S_{1\ \dots\ N}^{\alpha}\ =\ \sum_{i\ =\ 1}^{N}\ S_i^{\alpha}\ ,
\label{eq:Sn}
\end{equation}
where $S_i^{\alpha}$ is the $\alpha$ component of the spin operator at site $i$. The value of the spin operator in the $F$-basis is given as, 
\begin{equation}
\tilde{S}_{1\ \dots\ N}^{\alpha}\ =\ F_{1\ \dots\ N}\ S_{1\ \dots\ N}^{\alpha}\ F_{1\ \dots\ N}^{-1}\ ,
\label{eq:tSn}
\end{equation}
where, $F_{1\ \dots\ N}\ \equiv\ F_{1\ \dots\ N}\ (\xi_1\ \dots\ \xi_N)$. 
We have the following instructive result for $N\ =\ 2$ :
\begin{lem}
We have, $\tilde{S}_{12}^{\pm}\ =\ S_{12}^{\pm}\ \Delta_{12}$, where $\Delta_{12}\ (\xi_1,\ \xi_2)$ is the diagonal matrix defined in the last section as,
\begin{eqnarray*}
\Delta_{12}\ (\xi_1,\ \xi_2)\ =\ 
\left(
\begin{array}{cccc}
1 &  &   &  \\
& b^{-1}(\xi_2 , \xi_1) &  &  \\
 &  &b^{-1}(\xi_1, \xi_2) &  \\
 &  &  & 1
\end{array}
\right)_{[12]}\ ,
\end{eqnarray*}
and such that $\Delta_{12}\ (\xi_1,\ \xi_2)\ =\ \Delta_{21}\ (\xi_2,\ \xi_1)$.
\label{lem:recS}
\end{lem}
\Proof
By direct computation of products of $4 \times 4$ matrices.
\QED
In order to give the result for $\tilde{S}_{1\ \dots\ N}^{\alpha}$ in a closed form, let us first prove the following Lemma.
\begin{lem}
\begin{equation}
\tilde{S}_{1\ \dots\ N}^{-}\ =\ \left( S_1^{-}\ \ +\ e_1^{(11)}\ \tilde{S}_{2\ \dots\ N}^{-}\ {\cal C}\ \tilde{D}_{2\ \dots\ N}^{\dagger}\ (\xi_1)\ {\cal C}\  
+\ e_1^{(22)}\ \tilde{D}_{2\ \dots\ N}\ (\xi_1)\ \tilde{S}_{2\ \dots\ N}^{-}\ \right)\ \Delta_{1,\ 2\ \dots\ N}\ ,
\\
\end{equation}
where the matrices $\Delta_{1,\ 2\ \dots\ N}$, ${\cal C}\ \tilde{D}_{2\ \dots\ N}^{\dagger}\ (\xi_1)\ {\cal C}$ and $\tilde{D}_{2\ \dots\ N}\ (\xi_1)$ are the following diagonal matrices,
\begin{eqnarray*}
\Delta_{1,\ 2\ \dots\ N}\ &=&\ \Delta_{12\ \dots\ N}\ \Delta_{2\ \dots\ N}^{-1}\ ,  \\
{\cal C}\ \tilde{D}_{2\ \dots\ N}^{\dagger}\ (\xi_1)\ {\cal C}\ &=&\ \otimes_{i\ =\ 2}^{N}\ \left(
\begin{array}{cc}
1 & 0  \\
0 & b\ (\xi_i,\ \xi_1)
\end{array}
\right)_{[i]}\ ,
 \\
\tilde{D}_{2\ \dots\ N}\ (\xi_1)\ &=&\  \otimes_{i\ =\ 2}^{N}\ \left(
\begin{array}{cc}
b\ (\xi_1,\ \xi_i) & 0  \\
0 & 1
\end{array}
\right)_{[i]}\ .
\end{eqnarray*}
\end{lem}
\Proof 
We first recall the following property of $F$-matrices taken from Proposition (\ref{prop:fftd}), 
\begin{eqnarray*}
F_{1\ \dots\ n}^{-1}\ =\ F_{n\ \dots\ 1}^{t_1\ \dots\ t_n}\ \Delta_{1\ \dots\ n}\ .
\end{eqnarray*}
Hence,
\begin{eqnarray*}
\tilde{S}_{1\ \dots\ n}^{-}\ &=&\ F_{1\ \dots\ n}\ S_{1\ \dots\ n}^{-}\ F_{n\ \dots\ 1}^{t_1\ \dots\ t_n}\ \Delta_{1\ \dots\ n}\ , \\
&=&\ F_{2\ \dots\ n}\ F_{1,\ 2\ \dots\ n}\ S_{1\ \dots\ n}^{-}\ F_{n\ \dots\ 2,\ 1}^{t_1\ \dots\ t_n}\ F_{2\ \dots\ n}^{-1}\ \Delta_{1,\ 2\ \dots\ n}\ .
\end{eqnarray*}
We have,
\begin{equation}
F_{n\ \dots\ 2,\ 1}\ =\ e_1^{(22)}\ +\ e_1^{(11)}\ R_{n\ \dots\ 2,\ 1}\ . \nonumber
\end{equation}
Therefore,
\begin{eqnarray*}
F_{n\ \dots\ 2,\ 1}^{t_1\ \dots\ t_n}\ &=&\ e_1^{(22)}\ +\ R_{n\ \dots\ 2,\ 1}^{t_1\ \dots\ t_n}\ e_1^{(11)}\ , \\
&=&\ e_1^{(22)}\ +\ R_{2\ \dots\ n,\ 1}\ e_1^{(11)}\ .
\end{eqnarray*}
Hence we obtain,
\begin{eqnarray*}
& &\tilde{S}_{1\ \dots\ n}^{-}\ =\\
&=&\ F_{2\ \dots\ n}\ (e_1^{(11)}\ +\ e_1^{(22)}\ R_{1,\ 2\ \dots\ n})\ S_{1\ \dots\ n}^{-}\ (e_1^{(22)}\ +\ R_{2\ \dots\ n,\ 1}\ e_1^{(11)})\ F_{2\ \dots\ n}^{-1}\ \Delta_{1,\ 2\ \dots\ n}\ ,\\
&=&\ F_{2\ \dots\ n}\ \left[\ S_1^{-}\ +\ e_1^{(22)}\ R_{1,\ 2\ \dots\ n}\ S_{1\ \dots\ n}^{-}\ e_1^{(22)}\ +\ e_1^{(11)}\ S_{1\ \dots\ n}^{-}\ R_{2\ \dots\ n,\ 1}\ e_1^{(11)}\ \right]\ \cdot\\ 
& &\cdot\ F_{2\ \dots\ n}^{-1}\ \Delta_{1,\ 2\ \dots\ n}\ ,
\end{eqnarray*}
where we used the commutation of $S_{1\ \dots\ n}^{-}$ with the $R$-matrix $R_{1,\ 2\ \dots\ n}$, the unitarity of the $R$-matrix and the elementary relations,
\begin{equation}
e_1^{(11)}\ S_1^-\ =\ 0\ \ \ S_1^-\ e_1^{(11)}\ =\ S_1^-\ \ \  S_1^-\ e_1^{(22)}\ =\ 0\ \ \ e_1^{(22)}\ S_1^-\ =\ S_1^-\ .
\end{equation}
Therefore we have,
\begin{eqnarray*}
\tilde{S}_{1\ \dots\ n}^{-}\ &=&\  F_{2\ \dots\ n}\ \left[\ S_1^{-}\ +\ e_1^{(22)}\ R_{1,\ 2\ \dots\ n})\ S_{2\ \dots\ n}^{-}\ e_1^{(22)}\ +\ e_1^{(11)}\ S_{2\ \dots\ n}^{-}\ R_{2\ \dots\ n,\ 1}\ e_1^{(11)}\ \right]\ \cdot\\ 
& &\cdot\ F_{2\ \dots\ n}^{-1}\ \Delta_{1,\ 2\ \dots\ n}\ .
\end{eqnarray*}
Now recall that $R_{2\ \dots\ n,\ 1}\ =\ R_{1,\ 2\ \dots\ n}^{\dagger}$ to get,
\begin{eqnarray*}
\tilde{S}_{1\ \dots\ n}^{-}\ =\ \left[\ S_1^{-}\ +\ e_1^{(22)}\ \tilde{D}_{2\ \dots\ n}\ (\xi_1)\ \tilde{S}_{2\ \dots\ n}^{-}\  +\ e_1^{(11)}\ \tilde{S}_{2\ \dots\ n}^{-}\ {\cal C}\ \tilde{D}_{2\ \dots\ n}^{\dagger}\ (\xi_1)\ {\cal C}\ \right]\  \Delta_{1,\ 2\ \dots\ n}\ ,
\end{eqnarray*}
where we used the diagonalization of $D_{2\ \dots\ n}$ and $A_{2\ \dots\ n}^{\dagger}$ by $F_{2\ \dots\ n}$ from Proposition (\ref{prop:fDA+}).
\QED
\begin{rem}
The only term in $\tilde{S}_{1\ \dots\ n}^{-}$ containing $S_1^-$ is given as $S_1^-\ \Delta_{1,\ 2\ \dots\ n}$. But we also know that for any permutation $\sigma \in S_n$,
\begin{eqnarray*}
R_{1\ \dots\ n}^{\sigma}\ S_{1\ \dots\ n}^{\alpha}\ &=&\ S_{1\ \dots\ n}^{\alpha}\ R_{1\ \dots\ n}^{\sigma}\ , \\
&=&\ S_{\sigma(1)\ \dots\ \sigma(n)}^{\alpha}\ R_{1\ \dots\ n}^{\sigma}\ ,
\end{eqnarray*}
since the total spin operator commutes with the $R$-matrices and is also completly symmetric under permutations of $S_n$. Hence for any element $\sigma \in S_n$,
\begin{equation}
\tilde{S}_{1\ \dots\ n}^{\alpha}\ =\ \tilde{S}_{\sigma(1)\ \dots\ \sigma(n)}^{\alpha}\ .
\end{equation}
So, by an elementary symmetry argument, we should have the following Proposition.
\end{rem}
\begin{prop}
\begin{equation}
\tilde{S}_{1\ \dots\ N}^{-}\ =\ \sum_{i\ =\ 1}^{N}\ S_i^-\ \Delta_i^N\ ,
\label{eq:tildes-}
\end{equation}
where $\Delta_i^N\ =\ \prod_{j \ne i,\ j=1}^N\ \Delta_{ij}$. 
\label{prop:fs-}
\end{prop}
\Proof
By induction on $n$ using the last Lemma. Indeed, the proposition is true for $n\ =\ 2$ from the first Lemma. Let it be true for $n-1$. Then from the above Lemma, we obtain, applying the induction hypothesis to $\tilde{S}_{2\ \dots\ N}^{-}$, 
\begin{eqnarray*}
\tilde{S}_{1\ \dots\ n}^{-}\  &=&\ S_1^{-}\ \Delta_{1,\ 2\ \dots\ n}\ +\ e_1^{(22)}\ \tilde{D}_{2\ \dots\ n}\ (\xi_1)\ \sum_{i\ =\ 2}^{n}\ S_i^-\ (\prod_{j \ne i,\ j\ =\ 2}^{n}\ \Delta_{ij})\ \Delta_{1,\ 2\ \dots\ n}\  +\\ 
& &+ \ e_1^{(11)}\ \sum_{i\ =\ 2}^{n}\ S_i^-\ (\prod_{j \ne i,\ j\ =\ 2}^{n}\ \Delta_{ij})\ {\cal C}\ \tilde{D}_{2\ \dots\ n}^{\dagger}\ (\xi_1)\ {\cal C}\ \Delta_{1,\ 2\ \dots\ n}\ .
\end{eqnarray*}
Now we have the following equality for diagonal matrices involved in the above formula,
\begin{eqnarray*}
e_1^{(11)}\ \Delta_{1,\ 2\ \dots\ N}\ &=&\ {\cal C}\ (\tilde{D}_{2\ \dots\ n}^{\dagger}\ (\xi_1))^{-1}\ e_1^{(11)}\ {\cal C}\ ,\\
e_1^{(22)}\ \Delta_{1,\ 2\ \dots\ N}\ &=&\ \tilde{D}_{2\ \dots\ n}^{-1}\ (\xi_1)\ e_1^{(22)}\ .
\end{eqnarray*}
We also have,
\begin{eqnarray*}
S_i^-\ \prod_{j \ne i,\ j\ =\ 2}^{n}\ \Delta_{ij}\ =\ S_i^-\ \otimes_{j \ne i}\ \left( 
\begin{array}{cc}
1 & 0\\
0 & b^{-1}\ (\xi_j,\ \xi_i)
\end{array}
\right)_{[j]}\ .
\end{eqnarray*}
We finally obtain the result using $e_1^{(11)}\ +\ e_1^{(22)}\ =\ {\mathbf 1}_1$.
\QED
\begin{rem}
This formula defines the action of $S_{1\ \dots\ N}^-$ in the $F$-basis as a quasi-local operator since each $S_i^-$ in the above sum is ``dressed'' by a diagonal matrix acting on all sites but $i$.
\end{rem}
In the same way we obtain,
\begin{prop}
\begin{equation}
\tilde{S}_{1\ \dots\ N}^{+}\ =\ \sum_{i\ =\ 1}^{N}\ S_i^+\ \Delta_i^N\ ,
\label{eq:tildes+}
\end{equation}
and,
\begin{equation}
\tilde{S}_{1\ \dots\ N}^{z}\ =\ \sum_{i\ =\ 1}^{N}\ S_i^z\ = \ S_{1\ \dots\ N}^{z}\ .
\label{eq:tildesz}
\end{equation}
\label{prop:fs+z}
\end{prop}
\Proof
The formula for $\tilde{S}_{1\ \dots\ N}^{+}$ follows along the same arguments as the one for $\tilde{S}_{1\ \dots\ N}^{-}$. Hence we give here only the derivation for $\tilde{S}_{1\ \dots\ N}^{z}$. It is obtain from the fact that $S_{1\ \dots\ N}^{z}$ commutes with any partial $F$-matrix, namely we have,
\begin{eqnarray*}
F_{1,\ 2\ \dots\ n}\ S_{1\ \dots\ n}^{z}\ &=&\ (e_1^{(11)}\ +\ e_1^{(22)}\ R_{1,\ 2\ \dots\ n})\ S_{1\ \dots\ n}^{z}\ , \\
&=&\ S_{1\ \dots\ n}^{z}\ F_{1,\ 2\ \dots\ n}
\end{eqnarray*}
since $S_{1\ \dots\ N}^{z}$ commutes with the $R$-matrix of the $XXX$-$1 \over 2$ model and also with the diagonal matrices $e_1^{(11)}$ and $e_1^{(22)}$. Hence, $S_{1\ \dots\ N}^{z}$ commutes with the total $F$-matrix $F_{1\ \dots\ N}$.
\QED
\begin{rem}
The components of $\tilde{S}_{1\ \dots\ N}$ obey the same commutation relations as the one of $S_{1\ \dots\ N}$ which can be checked by direct computation using some identities of rational functions.
\end{rem}

\subsection{The Monodromy Matrix and the $F$-matrix in the $F$-basis}

Using the results of the last section, we can now compute the values of the operators $A,\ B,\ C,\ D$ and hence of the quantum monodromy matrix in the $F$-basis. First recall that,
\begin{eqnarray*}
F_{1\ \dots\ N}\ (\xi_1,\  \dots\ ,\ \xi_N)\ D_{1\ \dots\ N}\ (\lambda\ ;\ \xi_1,\  \dots\ ,\ \xi_N)\ F_{1\ \dots\ N}^{-1}\ (\xi_1,\  \dots\ ,\ \xi_N)\ =\\ =\ \tilde{D}_{1\ \dots\ N}\ (\lambda\ ;\ \xi_1,\  \dots\ ,\ \xi_N)\ =\  \otimes_{i\ =\ 1}^{N}\ 
\left( 
\begin{array}{cc}
b\ (\lambda,\ \xi_i) & 0\\
0 & 1
\end{array}
\right)_{[i]}\ .
\end{eqnarray*}
Then we have the following Proposition.
\begin{prop}
The operator $\tilde{B}$ representing the operator $B$ in the $F$-basis is given by,
\begin{eqnarray*}
\tilde{B}_{1\ \dots\ N}\ (\lambda)\ =\ \sum_{i\ =\ 1}^{N}\ S_i^-\ c\ (\lambda,\ \xi_i)\ \otimes_{j \ne i}\ 
\left( 
\begin{array}{cc}
b\ (\lambda,\ \xi_j) & 0\\
0 & b^{-1}\ (\xi_j,\ \xi_i)
\end{array}
\right)_{[j]}\ .
\end{eqnarray*}
Similarly we have for the operator $\tilde{C}$,
\begin{eqnarray*}
\tilde{C}_{1\ \dots\ N}\ (\lambda)\ =\ \sum_{i\ =\ 1}^{N}\ S_i^+\ c\ (\lambda,\ \xi_i)\ \otimes_{j \ne i}\ 
\left( 
\begin{array}{cc}
b\ (\lambda,\ \xi_j)\ b^{-1}\ (\xi_i,\ \xi_j) & 0\\
0 & 1
\end{array}
\right)_{[j]}\ ,
\end{eqnarray*}
and for the operator $\tilde{A}$,
\begin{eqnarray*}
\tilde{A}_{1\ \dots\ N}\ (\lambda)\ =\ \tilde{D}_{1\ \dots\ N}\ (\lambda)\ +\\
+\  \sum_{i\ =\ 1}^{N}\ S_i^z\ c\ (\lambda,\ \xi_i)\ \otimes_{j \ne i}\ 
\left( 
\begin{array}{cc}
b\ (\lambda,\ \xi_j)\ b^{-1}\ (\xi_i,\ \xi_j) & 0\\
0 & b^{-1}\ (\xi_j,\ \xi_i)
\end{array}
\right)_{[j]}\ +\ \sum_{i \ne j}^{N}\ S_i^+\ S_j^-\ \cdot \\
\cdot\ c\ (\lambda,\ \xi_i)\ c\ (\lambda,\ \xi_j) \ b^{-1}\ (\xi_i,\ \xi_j) \otimes_{k \ne i,\ k \ne j}\ 
\left( 
\begin{array}{cc}
b\ (\lambda,\ \xi_k)\ b^{-1}\ (\xi_i,\ \xi_k) & 0\\
0 & b^{-1}\ (\xi_k,\ \xi_j)
\end{array}
\right)_{[k]}\ .
\end{eqnarray*}
\label{prop:tildeABC}
\end{prop}
\Proof
We have from equations (\ref{eq:ABCDr}),
\begin{eqnarray*}
\tilde{B}_{1\ \dots\ N}\ (\lambda)\ &=&\ \left[\ \tilde{D}_{1\ \dots\ N}\ (\lambda)\ ,\ \tilde{S}_{1\ \dots\ N}^- \ \right]\ , \\
\tilde{C}_{1\ \dots\ N}\ (\lambda)\ &=&\ \left[\ \tilde{S}_{1\ \dots\ N}^+\ ,\ \tilde{D}_{1\ \dots\ N}\ (\lambda)\ \right]\ .
\end{eqnarray*}
This leads to,
\begin{eqnarray*}
\tilde{B}_{1\ \dots\ N}\ (\lambda)\ &=&\ \left[\ \tilde{D}_{1\ \dots\ N}\ (\lambda)\ ,\ \sum_{i\ =\ 1}^{N}\ S_i^-\ \Delta_i^N\ \right]\ =\ \sum_{i\ =\ 1}^{N}\ \left[\ \tilde{D}_{1\ \dots\ N}\ (\lambda)\ ,\ S_i^-\ \Delta_i^N\ \right]\ , \\
&=&\ \sum_{i\ =\ 1}^{N}\ S_i^-\ \left[\ \tilde{D}_{1\ \dots\ i-1,\ i+1,\ \dots\  N}\ (\lambda)\ -\ \tilde{D}_{1\ \dots\ N}\ (\lambda)\ \right]\ \Delta_i^N\ ,
\end{eqnarray*}
leading to the result since we also have,
\begin{eqnarray*}
\left( 
\begin{array}{cc}
0 & 0\\
1 & 0
\end{array}
\right)\ \cdot\ 
\left( 
\begin{array}{cc}
b & 0\\
0 & 1
\end{array}
\right)\ =\ b\ \left( 
\begin{array}{cc}
0 & 0\\
1 & 0
\end{array}
\right)\ .
\end{eqnarray*}
The computation for $\tilde{C}$ is similar. For $\tilde{A}$ we use the relation,
\begin{equation}
A\ (\lambda)\ =\ D\ (\lambda)\ +\ \left[\ S_{1\ \dots\ N}^+\ ,\ B_{1\ \dots\ N}\ (\lambda)\ \right]\ ,
\end{equation}
which rewritten in the $F$-basis leads to the result.
\QED
At this point several remarks are in order.
\begin{rem}
We wish to stress that the operators $\tilde{A},\ \tilde{B},\ \tilde{C},\ \tilde{D}$ satisfy the same quadratic commutation relations as the one satisfied by $A,\ B,\ C,\ D$. 
\end{rem}
\begin{rem}
The operators $\tilde{B}$ and $\tilde{C}$ are reduced to sums of only $N$ elementary operators which is to be compared to their expressions in the original basis where they are given as sums of $2^N$ terms involving much more complicated operators. Here, they are given essentially as sums on sites of the corresponding spin operator at each site {\em dressed diagonally}. In fact the action of $F$ essentially reduces the sum (\ref{eq:B2N}) to a term similar in structure to the first term of (\ref{eq:B2N}).\\
It really means that the factorizing $F$-matrices we have constructed solve the combinatorial problem induced by the non-trivial action of the permutation group $S_N$ given by the $R$-matrix. In the $F$-basis the action of the permutation group on the operators $\tilde{A}, \tilde{B}, \tilde{C}, \tilde{D}$ is trivial.  Moreover the operator $\tilde{A}\ +\ \tilde{D}$ which contains the Hamiltonian of the model together with the series of conserved quantities is now a quasi-bi-local operator.
\end{rem}
\begin{rem}
The operators $B$ and $C$ in the $F$-basis can be viewed as a diagonal dressing of the corresponding operators in the Gaudin model. The same is true for the Hamiltonian of the inhomogeneous $XXX$-$1 \over 2$ model in the $F$-basis  compared to the Hamiltonian of the Gaudin model. The ``linearization'' limit of the  $XXX$-$1 \over 2$ model, given by taking the parameter $\eta$ going to zero, indeed reproduce the Gaudin model from the above formula. In this limit, all diagonal ``dressing'' matrices are sent to the identity matrix.\\
Hence the $XXX$-$1 \over 2$ model  rewritten in the $F$-basis can be interpreted as a ``diagonal dressing'' of the Gaudin model.
\end{rem}
\begin{rem}
It can be shown that the pseudo-vacuum state is left invariant, namely, it is an eigenvector of the total $F$-matrix with eigenvalue $1$. Hence, in particular the Algebraic Bethe Ansatz can be carried out also in the $F$-basis. Moreover for the scalar products of the quantum states of the model, we have,
\begin{eqnarray*}
< 0\ |\ C\ (\lambda_1)\ \dots\ C\ (\lambda_n)\ B\ (\lambda_{n+1})\ \dots\ B\ (\lambda_{2n})\ |\ 0>\ = \\ 
=\ < 0\ |\ \tilde{C}\ (\lambda_1)\ \dots\ \tilde{C}\ (\lambda_n)\ \tilde{B}\ (\lambda_{n+1})\ \dots\ \tilde{B}\ (\lambda_{2n})\ |\ 0>\ .
\label{eq:sp}
\end{eqnarray*}
This should shed some new lights on the Gaudin formula for these scalar products here in the case of the $XXX$-$1 \over 2$ model.
\end{rem}

Using this remarks, an interesting exercise is to compute the first Bethe state given as the action of $\tilde{B}\ (\lambda)$ on the pseudo vacuum state. It can be shown directly that Bethe equations for this state are equivalent to the condition for this state to be highest weight with respect to $\tilde{S}$.\\
Finally, we wish to come back to the problem of the computation of the factorizing $F$-matrices defined in section 3. There we obtained factorizing $F$-matrices from an ordered product of partial $F$-matrices like $F_{1,\ 2\ \dots\ n}$. These object are given in term of the $R$-matrix $R_{1,\ 2\ \dots\ n}$. However this object is hightly non trivial to compute explicitly since it involves in fact sums of $2^{n-1}$ terms. In contrast, the partial $F$-matrices in the $F$-basis can be obtained explicitly, while they also lead to the construction of factorizing $F$-matrices $F_{12\ \dots\ n}$. We have,
\begin{eqnarray*}
F_{1\ \dots\ n}\ =\ \tilde{F}_{1,\ 2\ \dots\ n}\ F_{2\ \dots\ n}\ =\        \tilde{F}_{1,\ 2\ \dots\ n}\ \tilde{F}_{2,\ 3\ \dots\ n}\ \dots\ F_{n-1\ n}\ ,
\end{eqnarray*}
where,
\begin{eqnarray*}
\tilde{F}_{1,\ 2\ \dots\ n}\ &=&\ F_{2\ \dots\ n}\ F_{1,\ 2\ \dots\ n}\ F_{2\ \dots\ n}^{-1}\\
&=&\ \left(
\begin{array}{cc}
{\mathbf 1} & {\mathbf 0} \\
\tilde{C}_{2\ \dots\ n}\ (\xi_1) & \tilde{D}_{2\ \dots\ n}\ (\xi_1)
\end{array}
\right)_{[1]}\ ,
\end{eqnarray*}
is a very simple object to compute from the results of this section. Hence using the $F$-basis we have also obtained a more explicit and elementary formula for the $F$-matrix itself.\\
These results are relevant to the algebraic structure of the correlation functions of this model.

\par \vskip .5in \noindent
{\bf Aknowledgements.} The authors would like to thank their colleagues of the ENSLAPP Laboratory for Theoretical Physics in ENS Lyon for useful discussions and remarks. One of the authors (J. S. de S.) would like to thank the ENSLAPP Laboratory for Theoretical Physics in ENS Lyon for hospitality during which this work has been done. One of the authors (J.M.M.) would like to thank A. Izergin, N. Reshetikhin, E. Sklyanin and V. Tarassov for useful discussions.

%% file: F.bbl
\begin{thebibliography}{99}
\bibitem{Bax1}
Baxter, R.J.: Exactly solved models in statistical mechanics. London, Academic Press (1982) (and references therein)
\bibitem{CP1}
Chari, V. and Pressley, A.: Quantum affine algebras. Commun. Math. Phys. {\bf 142}, 261 (1991)
\bibitem{Drin1}
Drinfel'd, V.G.: Hopf algebras and the quantum Yang-Baxter equation. Soviet. Math. Dokl. Vol. {\bf 32}, 1060 (1985)
\bibitem{Drin2}
Drinfel'd, V.G.: Quantum Groups. in Proc. of the I. C. M. (Berkeley, 1986)
\bibitem{Drin3}
Drinfel'd, V.G.: Quasi-Hopf algebras. Leningrad Math. J. {\bf 1}, 1419 (1990)
\bibitem{Drin4}
Drinfel'd, V.G.: Constant quasi-classical solutions of the Yang-Baxter quantum equation. Sov. Math. Dokl. {\bf 28}, 667 (1983)
\bibitem{Drin5}
Drinfel'd, V.G.: A new realization of Yangians and quantum affine algebras. 
Sov. Math. Dokl. {\bf 36}, 212 (1988)
\bibitem{E1}
Enriquez, B.: Rational forms for twistings of envelopping algebras of simple Lie algebras. Lett. Math. Phys. {\bf 25}, 111 (1992)
\bibitem{EK1}
Etingof, P. and Kazhdan, D.: Quantization of Lie bialgebras, I. q-alg/9506005 
\bibitem{F1}
Faddeev, L.D.: Integrable models in $1 + 1$ dimensional quantum field theory. In Les Houches (1982), 561, Elsevier Science Publ. (1984)
\bibitem{FRT}
Faddeev, L.D., Reshetikhin N., and Takhtadzhyan, L.: Quantization of Lie groups and Lie algebras. Leningrad Math. J. {\bf 1}, 193 (1990)
\bibitem{FST}
Faddeev, L.D., Sklyanin, E.K., and Takhtajan, L.: The quantum inverse scattering method. Theor. Math. Phys. {\bf 40}, 688 (1979)
\bibitem{FR1}
Frenkel, I.B. and Reshetikhin, N.: Quantum affine algebras and holonomic difference equations. Commun. Math. Phys. {\bf 146}, 1 (1992)
\bibitem{Gaud1}
Gaudin, M.: La fonction d'onde de Bethe. Masson, Paris (1983)
\bibitem{IK1}
Izergin, A.G. and Korepin, V.E.: The quantum inverse scattering method approach to correlation functions. Commun. Math. Phys. {\bf 94}, 67 (1984)
\bibitem{Jim1}
Jimbo, M.: A q-difference analogue of ${\cal U}({\cal G})$ and the Yang-Baxter equation. Lett. Math. Phys. {\bf 10}, 63 (1985)
\bibitem{Jim2}
Jimbo, M.: A q-analogue of ${\cal U}(gl(N+1))$, Hecke algebra, and the Yang-Baxter equation. Lett. Math. Phys. {\bf 11}, 247 (1986) 
\bibitem{Jim3}
Jimbo, M.:  Quantum R-matrix for the generalized Toda system. Commun. Math. Phys. {\bf 102},  537 (1986)
\bibitem{JM1}
Jimbo, M. and Miwa, T.:Algebraic analysis of solvable lattice models. AMS, (1995)
\bibitem{KT1}
Khoroshkin, S.M. and Tolstoy, V.N.: Twisting of quantum (super) algebras. Connection of Drinfel'd and Cartan-Weyl realizations for quantum affine algebras. Preprint MPI-94-23, hepth/9404036
\bibitem{KZ}
Knizhnik, V.G. and Zamolodchikov, A.B.: Current algebra and Wess-Zumino model in two dimensions. Nucl. Phys. B {\bf 247}, 83 (1984)
\bibitem{KBI}
Korepin, V. E., Bogoliubov, N.M. and Izergin, A.G. : Quantum Inverse scattering method and correlation functions. Cambridge University Press, Cambridge (1993)
\bibitem{KR1}
Kulish, P.P. and Reshetikhin, N.: Quantum linear problem for the Sine-Gordon equation and higher representations. Zap. Nauch. semin. LOMI, {\bf 101}, 101 (1981)
\bibitem{KS1}
Kulish, P.P, and Sklyanin, E.K.: Quantum spectral transform method. Recent developments. In Lect. Notes in Phys. vol. {\bf 151}, 61, Springer verlag (1982)
\bibitem{JMM1}
Maillet, J.M.: Lax equations and quantum groups. Phys. Lett. B {\bf 245}, 480 (1990)
\bibitem{R1}
Reshetikhin, N.: Multiparameter quantum groups and twisted quasitriangular Hopf algebras. Lett. Math. Phys. {\bf 20}, 331 (1990)
\bibitem{Skl1}
Sklyanin, E.K.: On some algebraic structures connected with the Yang-Baxter equation. Funct. Anal. Appl. {\bf 16}, 263 (1983)
\bibitem{Skl2}
Sklyanin, E.K.: Quantum version of the method of inverse scattering problem. Zap. Nauchn. Semin. LOMI {\bf 95}, 55 (1980)







\end{thebibliography}
